\documentclass[twoside,11pt]{article}

\usepackage{blindtext}

\usepackage{jmlr2e}

\usepackage{lipsum}
\usepackage{amsfonts}
\usepackage{amsmath}
\usepackage{graphicx,subfigure}
\usepackage{epstopdf}
\usepackage{algorithm}
\usepackage{algorithmic}
\usepackage{booktabs}
\usepackage{multirow} 
\usepackage{caption}
\usepackage{subcaption}
\usepackage{bm}

\newtheorem{assumption}[theorem]{Assumption}

\usepackage{lastpage}
\jmlrheading{23}{2026}{1-\pageref{LastPage}}{\today}{}{21-0000}{Zhangyong Liang, Huanhuan Gao}

\ShortHeadings{Neural Multiscale Decomposition}{Zhangyong Liang, Huanhuan Gao}
\firstpageno{1}

\begin{document}

\title{Neural Multiscale Decomposition for the Nonlinear Klein-Gordon Equation with Time Oscillation}

\author{\name Zhangyong Liang 
       \email zyliang1994@tju.edu.cn \\
       \addr
       National Center for Applied Mathematics \\
       Tianjin University \\
       Tianjin, 300072, China
       \AND \name Huanhuan Gao\thanks{Corresponding author.}
       \email gao\_huanhuan@jlu.edu.cn\\
       \addr  School of Mechanical and Aerospace Engineering, Jilin University \\
       Jilin University \\
       Changchun, 130025, China
       }

\editor{With editor}

\maketitle

\begin{abstract}
In this paper, we propose a neural multiscale decomposition method (NeuralMD) for solving the nonlinear Klein--Gordon equation (NKGE) with a dimensionless parameter $\varepsilon\in(0,1]$ from the relativistic regime to the nonrelativistic limit regime. 
The solution of the NKGE propagates waves with wavelength at $O(1)$ and $O(\varepsilon^2)$ in space and time, respectively, which brings the oscillation in time.
Existing collocation-based methods for solving this equation lead to spectral bias and propagation failure.
To mitigate the spectral bias induced by high-frequency time oscillation, we employ a multiscale time integrator (MTI) to absorb the time oscillation into the phase. 
This decomposes the NKGE into a nonlinear Schr\"odinger equation with wave operator (NLSW) with well-prepared initial data and a remainder equation with small initial data. 
As $\varepsilon \to 0$, the NKGE converges to the NLSW at rate $O(\varepsilon^{2})$, and the contribution of the remainder equation becomes negligible. 
Furthermore, to alleviate propagation failure caused by medium-frequency time oscillation, we propose a gated gradient correlation correction strategy to enforce temporal coherence in collocation-based methods. 
As a result, the approximation of the remainder term is no longer affected by propagation failure.
Comparative experiments with existing collocation-based methods demonstrate the superior performance of our method for solving the NKGE with various regularities of initial data over the whole regime. 
\end{abstract}


\section{Introduction}
\label{sec:introduction}
The dimensionless nonlinear Klein--Gordon equation (NKGE) with real cubic nonlinearity in $d$-dimensions ($d=1,2,3$) \citep{bao2014mti,bao2012analy,grundland1992,messiah2014} is given by
\begin{equation}
\label{eqn:KG}
\left\{
  \begin{aligned}
    & \varepsilon^{2}\partial_{tt}u(\bm{x},t)-\Delta u(\bm{x},t)
      +\varepsilon^{-2}u(\bm{x},t)
      +\lambda u(\bm{x},t)^{3}=0,
      \bm{x}\in\mathbb{R}^{d},\ t>0,\\[2mm]
    & u(\bm{x},0)=\phi_{1}(\bm{x}),\qquad
      \partial_{t}u(\bm{x},0)=\varepsilon^{-2}\phi_{2}(\bm{x}),\qquad
      \bm{x}\in\mathbb{R}^{d},
  \end{aligned}
\right.
\end{equation}
Here $u=u(\bm{x},t)$ is real-valued, and $0<\varepsilon\le1$ is dimensionless.
$\lambda\in\mathbb{R}$ controls the nonlinearity, with defocusing sign $\lambda>0$ and focusing sign $\lambda<0$.
$\phi_{1}$ and $\phi_{2}$ denote prescribed real-valued, $\varepsilon$-independent initial data.
We write the continuous model on $\mathbb R^d$ for notational simplicity.
The experiments use periodic boxes $\Omega\subset\mathbb R^d$, interpreted as periodic truncations.

For $\varepsilon=1$, corresponding to the $O(1)$ wave-speed regime, many numerical methods exist \citep{bratsos2009numer,duncan1997symp,liu2018sym}. 
In particular, finite difference time domain (FDTD) methods \citep{duncan1997symp,jimenez1990,strauss1978} are highly efficient and accurate in this regime. 
For $0<\varepsilon\ll1$, efficient and accurate methods become harder to design and analyze.
The main difficulty is the highly oscillatory temporal behavior.
Uniformly accurate (UA) schemes address this nonrelativistic stiffness with $\varepsilon$-uniform error bounds \citep{bao2012analy}. 
They permit time steps independent of the fast temporal scale $O(\varepsilon^2)$. 
Examples include MTI \citep{bao2014mti,bao2017uniform,bao2019comp}, TSF \citep{chartier2015}, and NPI \citep{cai2021uni,cai2022uni,li2025uniform}. 
These methods achieve temporal super-resolution for highly oscillatory nonrelativistic solutions.

However, these schemes remain local in time and rely on refined temporal discretizations.
Long simulations therefore carry a substantial computational burden.
Each new final time or domain also requires recomputing the time-dependent solution.
Physics-driven collocation methods avoid precomputed data and time-step accumulation \citep{raissi2019phy,wang2023sci,cuomo2022sci}. 
They aim to recover the whole time interval in a single optimization.
For oscillatory NKGE, however, they can suffer spectral bias and propagation failure \citep{rahaman2019,wang2022resp}. 
Their accuracy may deteriorate as $\varepsilon\to0$.

Deep networks trained by gradient methods often learn low frequencies first \citep{xu2019freq,rahaman2019,xu2025view}. 
This behavior is often called the frequency principle. 
It is problematic when PDE solutions contain many active scales. 
High-frequency modes may appear late in training or remain unresolved \citep{luo2022upper}. 
Neural solvers may then miss sharp layers, oscillations, and multiscale patterns \citep{xu2019freq,rahaman2019}. 
The issue worsens in high-dimensional, stochastic, or geometrically complex settings \citep{xu2025view}. 
Thus multiscale PDE solvers need architectures that mitigate spectral bias \citep{xu2019freq,luo2022upper}.
Recent works have identified that PINNs struggle to converge to solutions when target functions exhibit high-frequency patterns \citep{krishnapriyan2021,moseley2023finite}.
The application of differential operators to neural networks complicates the loss landscape and makes optimization more difficult. 
Neural tangent kernel analyses connect slow convergence to high-frequency components \citep{wang2022why,farhani2022mom}.
Adaptive reweighting and Adam-type optimizers offer limited relief for highly oscillatory targets.
Fourier feature mappings can improve high-frequency representation in several PDE settings \citep{hertz2021sape,tancik2020fourier,li2023deep}.

Propagation failure is a second challenge. 
Collocation methods can break down even on relatively simple PDEs \citep{krishnapri2021}. 
They may produce incorrect solutions despite apparently successful residual optimization. 
Prior work links these failures to sampling, loss imbalance, optimizers, differentiation, and coefficient sensitivity \citep{wu2023comp,wang2023r3,yu2022grad,wu2024ropinn,wu2025propinn,rathore2024,shi2024stoch,wang2021under,wang2022why}. 
Interior residual enforcement can leave interior points weakly supervised by initial or boundary data \citep{wang2023r3}. 
Several methods adapt collocation distributions toward regions with large residuals \citep{lu2021deepxde,nabian2021}. 
Others modify losses, residual norms, or optimization procedures \citep{yu2022grad,wu2024ropinn,wang2022l2,rathore2024,shi2024stoch}. 
These remedies mainly target stability, rather than temporal information propagation. 
Many methods still treat collocation points as independent samples \citep{wong2022learn,rohrhofer2022,rohrhofer2022role}. 
Causality-aware training can improve information flow for time-dependent PDEs \citep{wang2022resp}.

We propose NeuralMD to address oscillation-induced spectral bias and propagation failure in the NKGE.
Our main contributions can be summarized as follows:
\begin{itemize}
    \item NeuralMD decomposes the NKGE solution into high-frequency phases and a low-frequency envelope.
    A fixed-time WKB expansion splits the problem into NLSW modulation and a small-initial-data remainder. 
    The network therefore learns a mildly oscillatory envelope in the nonrelativistic limit regime. 
    This reduces the effective temporal frequency seen by the neural network.
    \item NeuralMD is a two-stage pretraining framework. 
    Stage I approximates the modulated NLSW.  Stage II corrects amplitude errors with a remainder network. 
    As $\varepsilon \to 0$, we solve only the NLSW and recover the oscillatory NKGE solution via a WKB expansion without remainder terms.
    As $\varepsilon \to 1$, we also solve the remainder equation and use a WKB expansion with remainder terms. 
    For intermediate $\varepsilon$, residual selection determines how much learned dynamic remainder is retained.
    \item The remainder equation inherits oscillatory bands from the NKGE.
    We use gated gradient correlation to promote temporal coherence among neighboring collocation points.
    Random temporal perturbations form neighborhoods and correlate residuals across perturbed times.
    Low-residual regions are downweighted, while high-residual regions receive additional samples.
    \item We also explore a KAN-based NeuralMD variant for structured interpretability \citep{liu2024kan}. 
    We test NeuralMD across initial-data regularities, dimensions, and representative $\varepsilon$ values.
    The results show lower errors than tested baselines in the reported settings.
\end{itemize}

The remainder of the paper is organized as follows. 
Section \ref{sec:problem} introduces the phenomenon of spectral bias and propagation failure induced by the time oscillation of the NKGE.
Section \ref{sec:method} presents the NeuralMD method and discusses its implementation. 
Section \ref{sec:experiment} reports numerical results. 
Finally, Section \ref{sec:conclusion} concludes the paper and discusses future directions.

\section{Main Problem}
\label{sec:problem}

\subsection{Oscillation in time induces spectral bias}
For \eqref{eqn:KG}, as $\varepsilon \to 0$, the solution develops rapidly oscillatory temporal phases.
Their carrier frequencies are $\pm 1/\varepsilon^{2}$, giving temporal wavelength $O(\varepsilon^{2})$ but spatial scale $O(1)$. 
This makes the problem extremely stiff in time. 
In the same regime, the conserved energy satisfies $E(t) = O(\varepsilon^{-2})$ and becomes unbounded as $\varepsilon \to 0$. 
This ``energy inflation'' complicates asymptotic analysis and error control.
Discretization errors are amplified at this scale and interact with the fast phase.

We next study spectral bias caused by time oscillation in collocation-based solvers such as PINNs.
The model $u_{\theta}$ approximates $u$ by minimizing residual, initial, and periodic-boundary losses:
\begin{equation}
\label{eq:auto-001}
\left\{
\begin{aligned}
    \mathcal{L}_{\mathrm{Res}}(\theta)
    &= \mathbb{E}_{(\bm{x},t)\in P_f}
    \left\| \varepsilon^2 \partial_{tt} u_\theta
        - \Delta u_\theta
        + \varepsilon^{-2} u_\theta
        + \lambda u_\theta^3
    \right\|_{L^2}^2, \\[4pt]
    \mathcal{L}_{\mathrm{Ic}}(\theta)
    &= \mathbb{E}_{(\bm{x},t)\in P_0}
    \left(
        \left\| u_\theta(x_0,0) - \phi_1 \right\|_{L^2}^2
        + \left\| \partial_t u_\theta(x_0,0) - \varepsilon^{-2} \phi_2 \right\|_{L^2}^2
    \right), \\[4pt]
    \mathcal{L}_{\mathrm{Bd}}(\theta)
    &= \mathbb{E}_{(\bm{x},t)\in P_b}
    \left(
        \left\| u_\theta(a,t_a) - u_\theta(b,t_b) \right\|_{L^2}^2
        + \left\| \partial_x u_\theta(a,t_a) - \partial_x u_\theta(b,t_b) \right\|_{L^2}^2
    \right).
\end{aligned}
\right.
\end{equation}
Here $P_f$, $P_0$, and $P_b$ are collocation sets.
They correspond to the residual, initial condition, and boundary condition, respectively.
$\mathbb{E}_{(\bm{x},t)\in P}$ denotes the empirical average over the set $P$. 
$[a,b]$ is the periodic spatial domain.

Here, we introduce the following composition \citep{dong2021period}
\begin{equation}
\label{eq:auto-002}
\begin{bmatrix} t & 1 & \cos\left(\frac{2\pi}{P}x\right) & \sin\left(\frac{2\pi}{P}x\right) & \cos\left(\frac{2\pi}{P}2x\right) & \cdots & \cos\left(\frac{2\pi}{P}mx\right) & \sin\left(\frac{2\pi}{P}mx\right) \end{bmatrix}^\top,
\end{equation}
where $m > 0$ is a positive integer.
A similar transformation can be applied to a higher-dimensional spatial domain $\bm{x}$. 
We then apply the following transformation to satisfy the initial condition automatically:
\begin{equation}
\label{eq:auto-003}
\tilde{u}_{\theta} = (1+t)e^{-t} u_0 + t e^{-t} \partial_t u_0 + \left(1 - e^{-t} - t e^{-t}\right) u_{\theta},
\end{equation}

The total loss used for training via automatic differentiation and backpropagation is the following
\begin{equation}
\label{eq:auto-004}
    \mathcal{L}(\theta)
    = \lambda_{\mathrm{Res}} \mathcal{L}_{\mathrm{Res}}(\theta)
    + \lambda_{\mathrm{Ic}} \mathcal{L}_{\mathrm{Ic}}(\theta)
    + \lambda_{\mathrm{Bd}} \mathcal{L}_{\mathrm{Bd}}(\theta).
\end{equation}

We analyze the impact of the time oscillation of NKGE on the training dynamics of PINNs. 
As an example, we consider the commonly used hyperbolic tangent activation function
\begin{equation}
\label{eq:auto-005}
    \sigma(t) = \tanh(t) = \frac{e^{t} - e^{-t}}{e^{t} + e^{-t}}, \quad t \in \mathbb{R}.
\end{equation}
and PINNs with one hidden layer, having $m$ neurons, a 1-dimensional input $t$, and a 1-dimensional output, i.e.,
\begin{equation}
\label{eq:auto-006}
h(t)=\sum_{j=1}^{m}a_{j}\sigma(w_{j}t+b_{j}),\quad a_{j},w_{j},b_{j}\in{\rm \mathbb{R}},
\end{equation}
where $\theta_{j} \triangleq \{w_{j}, b_{j}, a_{j}\}$, $w_{j}$, $a_{j}$, and $b_{j}$ are training parameters. 

Then, the Fourier transform of $h(x)$ can be computed as follows
\begin{equation}
\label{eqn:FTW}
    \hat{h}(k)=\sum_{j=1}^{m}\frac{2\pi a_{j}}{|w_{j}|}\exp \Big(\frac{b_{j}k}{w_{j}}\Big)\frac{1}{\exp(-\frac{\pi k}{2w_{j}})-\exp(\frac{\pi k}{2w_{j}})},
\end{equation}
where $k$ denotes the frequency.

By Parseval's theorem \citep{stein2011fourier}, this loss is identical to the standard mean-squared error in the Fourier domain, that is
\begin{equation}
\label{eq:auto-007}
    \mathcal{L}(\theta)
    = \int_{-\infty}^{+\infty}[\lambda_{\mathrm{Res}} \mathcal{L}_{\mathrm{Res}}(\theta,k)
    + \lambda_{\mathrm{Ic}} \mathcal{L}_{\mathrm{Ic}}(\theta,k)
    + \lambda_{\mathrm{Bd}} \mathcal{L}_{\mathrm{Bd}}(\theta,k))]d{k}.
\end{equation}

We study frequency-dependent error attenuation by expressing the loss in Fourier space. 
In the early stage of training, the weights typically satisfy $|w_j|\ll 1$, so the last term in \eqref{eqn:FTW} can be approximated by
\begin{equation}
\label{eq:auto-008}
\frac{1}{\exp(-\pi k/w_j) - \exp(\pi k/w_j)} \approx -\operatorname{sgn}(w_j) \exp\left(-\pi\left|\frac{k}{w_j}\right|\right)= \exp\left(-\left|\frac{\pi k}{2w_j}\right|\right),
\end{equation}

Hence, the magnitude of the contribution from frequency $k$ to the gradient with respect to $\theta_{lj}$ is
\begin{equation}
\label{eq:auto-009}
    \frac{\partial \mathcal{L}(\theta_j,k)}{\partial\theta_{j}}
    \approx |\frac{\partial \mathcal{L}(\theta_j,k)}{\partial w_{j}}|
    \exp\!\left(-\left|\frac{\pi k}{2w_{j}}\right|\right)
    F_{j}.
\end{equation}
where $F_{j}$ is an $O(1)$ function depending on $\theta_{j}$ and $k$.
$\exp\!\left(-\lvert \pi k/(2w_j)\rvert\right)$ indicates that low-frequency components dominate the initialized weights, while high-frequency components are exponentially suppressed.

In the NKGE, the temporal wavelength is $O(\varepsilon^{2})$; as $\varepsilon$ decreases, the carrier frequencies move toward $|k|=O(\varepsilon^{-2})$, where the factor $\exp\!\left(-\lvert \pi k/(2w_j)\rvert\right)$ becomes exponentially small for networks initialized with small weights.
When training with a residual-based physics loss, the linear part of the residual in the temporal frequency domain carries the coefficient
\begin{equation}
\label{eq:auto-010}
    \left| -\varepsilon^{2} k^{2} + \varepsilon^{-2} \right|,
\end{equation}

Near the Klein--Gordon characteristic band $k\approx\pm\varepsilon^{-2}$, temporal and mass terms may cancel.
Away from this set, the coefficient can reach $O(\varepsilon^{-2})$ on high temporal frequencies.
Thus the residual loss combines two difficult effects.
Network sensitivity is exponentially small at the carrier scale.
Residual weights are stiff for off-characteristic artifacts.
A schematic contribution to the gradient at temporal frequency $k$ is
\begin{equation}
\label{eq:auto-011}
    \frac{\partial L(k)}{\partial \theta_{j}}
    \approx \bigl| -\varepsilon^{2}k^{2} + \varepsilon^{-2} \bigr|
    \,\bigl\lVert u_{\theta}(k) - u(k) \bigr\rVert
    \exp\!\left(-\frac{\pi |k|}{2|w_{j}|}\right)\, F_{j}.
\end{equation}
Here $F_{lj}$ is an $O(1)$ factor.
It incorporates smoothing from the nonlinearity and spatial derivatives on this band. 

\begin{figure}[!htb]
    \centering
    \includegraphics[width=0.85\textwidth]{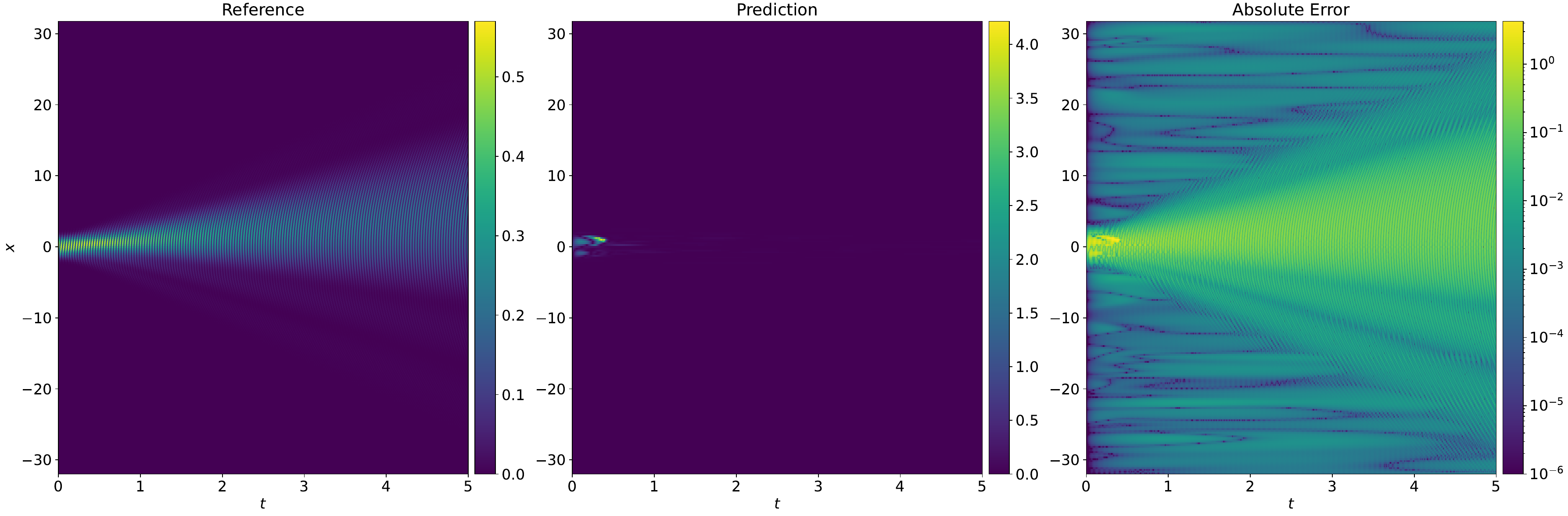}
    \vspace{-8pt}
    \caption{The prediction solution of PINNs for $\varepsilon=0.1$.}
    \label{fig:sp_bias_comp}
\end{figure}

Consequently, training is dominated by a mismatch in the nonrelativistic limit regime.
The solution contains high temporal frequencies, while the initialized network prefers low frequencies.
The residual operator also penalizes off-characteristic high-frequency artifacts.
This distinction is used later in the remainder-selection analysis.
Figure~\ref{fig:sp_bias_comp} shows that, in the pronounced nonrelativistic limit regime ($\varepsilon = 0.1$), PINNs exhibit optimization stagnation and high-frequency drift in the presence of strong time oscillation.

For reference, the key displayed relations in this subsection are \eqref{eq:auto-001}, \eqref{eq:auto-002}, \eqref{eq:auto-003}, \eqref{eq:auto-004}, \eqref{eq:auto-005}, \eqref{eq:auto-006}.
The continuation of this numbered set is \eqref{eq:auto-007}, \eqref{eq:auto-008}, \eqref{eq:auto-009}, \eqref{eq:auto-010}, \eqref{eq:auto-011}.

\subsection{Oscillation in time induces propagation failure}
We next ask whether a mild decrease in $\varepsilon$ already affects propagation in collocation-based solvers. 
The NKGE contains a potential term of size $\varepsilon^{-2}$.
This term induces high-frequency time oscillation.
More importantly, the initial time derivative also scales like $\varepsilon^{-2}$. 
Thus, as $\varepsilon \to 0$, high-frequency difficulty is present from the very beginning of the evolution.
When the initial gradient is high-frequency, PINNs struggle to propagate these dynamics into the domain. 
The learned solution typically stalls after a short time and enters a propagation-failure regime. 
This failure is not merely a ``frozen'' solution.
It reflects a mismatch between global PINN training and local-in-time PDE causality.

Figure~\ref{fig:sp_bias_comp} shows a failure mode that is distinct from the spectral bias induced by spatial oscillation. 
Starting from well-prepared initial data, the solution evolves only up to about $t \approx 0.5$ and then ceases to propagate.
This behavior is consistent with propagation failure. 
In this section, we explore how gradually increasing the oscillation frequency leads to propagation failure.

\begin{figure}[!htb]
    \centering
    \includegraphics[width=0.95\textwidth]{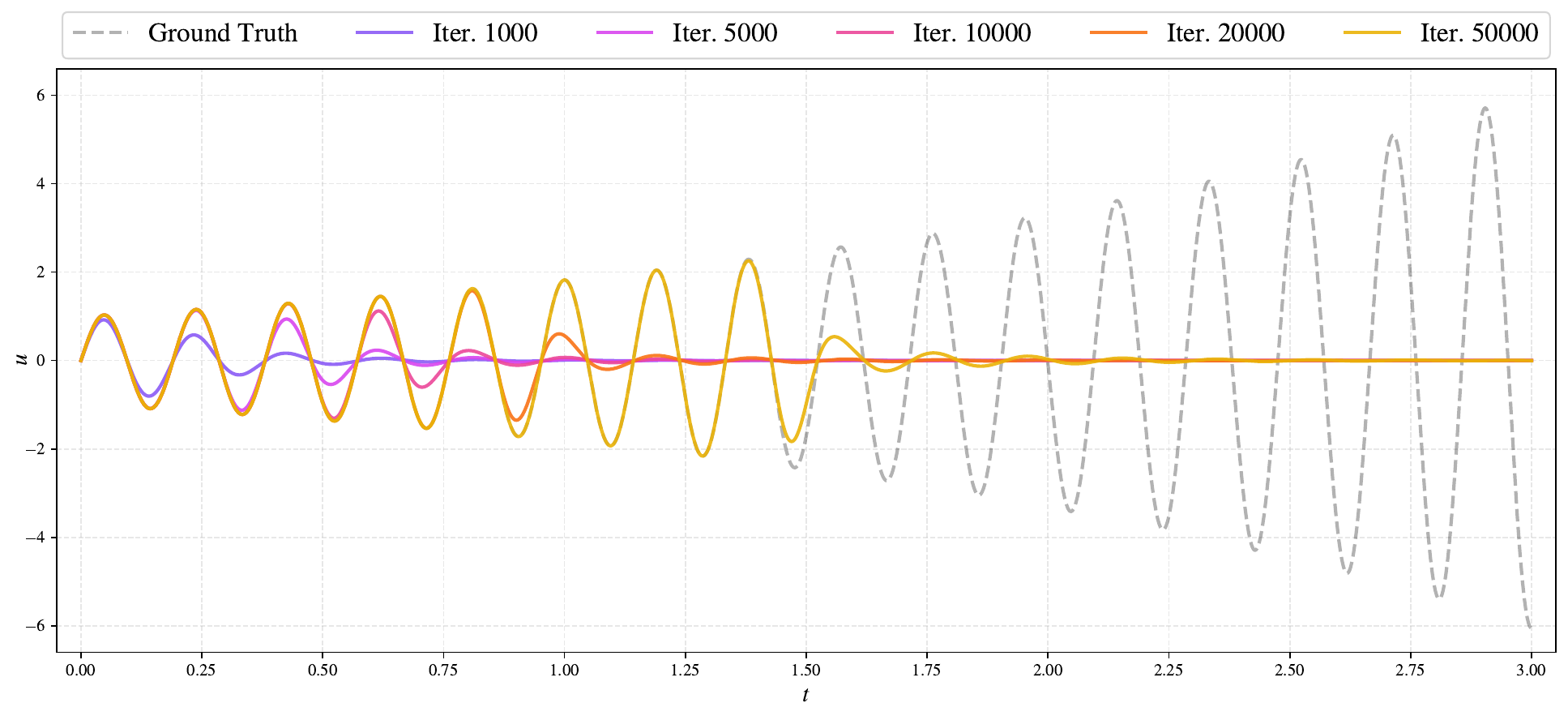}
    \vspace{-8pt}
    \caption{PINNs solutions for a temporally oscillatory ODE: $u_{tt} - 2\alpha u_t + (\alpha^2 + k^2)u = 0$ with $\alpha = 0.8$ and $k = 30$. The analytical solution is $u(t) = e^{\alpha t} \sin(kt)$, exhibiting exponentially growing oscillations in time.}
    \label{fig:pro_fail_intro}
\end{figure}

For the NKGE in the transition regime, the solution exhibits a temporal wavelength of $O(\varepsilon^{2})$, with mid-frequency time oscillation as $\varepsilon$ decreases. Both the potential term and initial time derivative scale like $\varepsilon^{-2}$, introducing substantial mid-frequency content from the start. 
In the PINN setting, the correct solution must propagate from initial/boundary points to interior collocation points. 
However, mid-frequency time oscillation complicates this propagation, making it difficult for the network to transmit the correct dynamics. 
As a result, some interior points may converge to trivial or low-frequency solutions before the correct solution can reach them. 
These trivial solutions then spread to nearby points, leading to large regions with incorrect solutions. 
Figure~\ref{fig:pro_fail_intro} illustrates this PINN propagation failure on a simple oscillatory problem.
As training iterations increase, optimization stagnates and the solution stops propagating correctly. 
Thus, mid-frequency oscillation and global PINN training together cause propagation failure.

For wave-type equations such as the NKGE, the time evolution is subject to causality and a finite propagation speed of information. 
The solution at a collocation point $(x_i,t_i)$ depends only on data in a certain region of the past. 
Traditional time-marching numerical schemes naturally respect this causal structure and effectively incorporate historical information. 
In contrast, PINNs optimize over the entire time domain at once.
This makes causal relations among temporal collocation points difficult to enforce.

Early in training, the network parameters $\theta$ are random.
Yet PINNs require $u_\theta$ to satisfy initial data at $t=0$.
They also impose the residual at $t=t_i$, often near the terminal time $T$. 
The network has not yet learned propagation from $t=0$ to $t=t_i$.
Rigid residual enforcement can then cause premature optimization stagnation.
This stagnation is not due to the failure of physical information propagation, but rather the failure of gradient information propagation during training. 
As a diagnostic, consider a local parameter update induced near time $t$.
We inspect its effect at the neighboring time $t+\delta t$:
\begin{equation}
\label{eq:auto-012}
D_{\mathrm{PINNs}}\left(t,t+\delta t\right)
=\lim_{\alpha\rightarrow0}
  \frac{\|u_\theta(\bm{x},t+\delta t)
  -u_{\theta-\alpha d_t}(\bm{x},t+\delta t)\|}{\alpha},
\end{equation}
where $d_t$ denotes a local training direction computed from samples near
time $t$.  This quantity is not an output-Jacobian inner product.
It only motivates temporal gradient alignment as a propagation diagnostic.
The formal method in Section~\ref{sec:gated_gradient} uses the
signed, normalized correlation of residual-loss gradients rather than the
absolute output-Jacobian correlation.  For intuition, define
\begin{equation}
\label{eq:auto-013}
G_{u_\theta}\left(t,t+\delta t\right)
=
\frac{\left\langle
  \left.\partial_\theta u_\theta\right|_t,
  \left.\partial_\theta u_\theta\right|_{t+\delta t}
\right\rangle}
{\|\left.\partial_\theta u_\theta\right|_t\|
 \|\left.\partial_\theta u_\theta\right|_{t+\delta t}\|+\eta_0}.
\end{equation}
Small or negative values indicate weak or conflicting influence between neighboring time slices.
As $\varepsilon \to 0$, parameter gradients at nearby times often contain high-frequency components.
At the same spatial location $x$, this motivates
\begin{equation}
\label{eq:auto-014}
G_{u_\theta}\left(\left(x,t\right),\left(x,t+\delta t\right)\right)\;\approx\;\int_{0}^{\infty}{S_x\left(k,\varepsilon\right)}\cos{\left(k\delta t\right)} dk,
\end{equation}
where $S_x\left(k,\varepsilon\right) \geq 0$ is the temporal angular frequency spectrum density of $\partial_\theta u_\theta$.  In the NKGE, the relevant high-frequency bands are centered near the carrier harmonics $m/\varepsilon^2$ rather than near $1/\varepsilon$.  A diagnostic model is
\begin{equation}
\label{eq:auto-015}
S_x(k,\varepsilon)
\approx S_{x,0}(k)+
\sum_{m\ne0}a_m(\varepsilon)
f_m\!\left(k-\frac{m}{\varepsilon^2}\right),
\end{equation}
where each $f_m$ has $O(1)$ width.  When at least one nonzero harmonic has
non-negligible weight, the second spectral moment is of order
$\Theta(\varepsilon^{-4})$.

In the nonrelativistic limit, high-frequency weights can dominate.
For a fixed $\delta t$, they cancel more easily across modes.
This makes $G_{u_\theta}$ decay faster and possibly become negative.
Adjacent-time parameter updates then become more orthogonal or conflicting.
Correct gradient information becomes harder to propagate along time.
This mechanism can trigger the evolution-failure mode.

In the nonrelativistic limit, severe oscillation makes this failure easy to attribute to spectral bias. 
To separate the effects, we choose $\varepsilon$ in a mild-oscillation regime.
We gradually increase the final time and test long-time PINN propagation.
Although oscillation is mild, the distribution of $S_x(k)$ still depends on $\delta t$.
Its effect on gradient directions accumulates over time. 
In long-time prediction, this accumulation gradually lowers adjacent-time gradient correlation.
It eventually triggers the evolution-failure mode. 
As $\varepsilon$ decreases, the accumulation effect accelerates, and the evolution failure mode is triggered more quickly. 
The gradient correlation decay satisfies
\begin{equation}
\label{eq:auto-016}
G_{u_\theta}\left(x,\delta t\right)\approx A\left(x\right)-\frac{1}{2}A\left(x\right)\left\langle k^2\right\rangle_{S,x}\delta t^2+O\left(\delta t^4\right),
\end{equation}
where $A\left(x\right)$ represents the norm of the parameter gradient, quantifying the gradient autocorrelation
\begin{equation}
\label{eq:auto-017}
A\left(x\right)=\int_{0}^{\infty}{S_x\left(k,\varepsilon\right)}dk=G_{u_\theta}(x,0)=\|\partial_\theta u_\theta\left(x,t\right)\|,
\end{equation}

$\left\langle k^2\right\rangle_{S,x}$ is the weighted second moment of high-frequency components.
It quantifies how high-frequency time oscillation drives gradient decorrelation. 
It is expressed as
\begin{equation}
\label{eq:auto-018}
\left\langle k^2\right\rangle_{S,x}=\frac{\int_{0}^{\infty}{k^2S_x\left(k,\varepsilon\right)}dk}{\int_{0}^{\infty}{S_x\left(k,\varepsilon\right)}dk}
=\Theta(\varepsilon^{-4})
\end{equation}

This leads to the following expression for $G_{u_\theta}$, as
\begin{equation}
\label{eq:auto-019}
G_{u_\theta}\left(x,\delta t\right)
\approx A\left(x\right)-\frac{1}{2}A\left(x\right)
\left\langle k^2\right\rangle_{S,x}\delta t^2+O\left(\delta t^4\right).
\end{equation}

As $\varepsilon$ decreases, $\left\langle k^2\right\rangle_{S,x}$ increases, causing the gradient correlation between adjacent time steps to decay faster. 
For $\varepsilon$ below 1, gradient correlation can decay enough to trigger evolution failure.
This failure is more pronounced in long-time predictions.

\begin{figure}[!htb]
    \centering
    \includegraphics[width=0.85\textwidth]{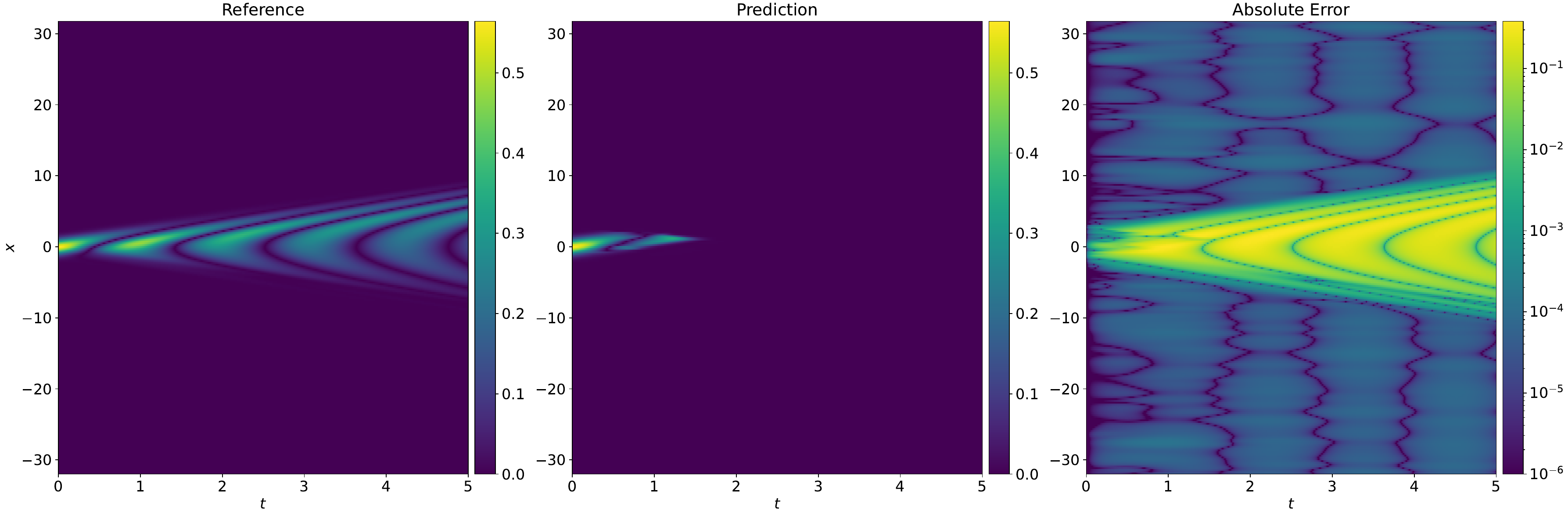}
    \vspace{-8pt}
    \caption{The prediction solution of PINNs for $\varepsilon=0.6,t=5.0$.}
    \label{fig:pro_fail_comp}
\end{figure}

To assess this effect, we set $T = 5$ and reduce $\varepsilon$ to $0.6$.
The network architecture and training hyperparameters are fixed.
Figure~\ref{fig:pro_fail_comp} shows that isolated temporal collocation does not enforce causality.
It also does not enforce recursive consistency.
As a result, phase errors from different sampling times accumulate and decorrelate.
This reduces the model's ability to track the true solution. 
When $\varepsilon = 0.6$, optimization stagnation appears near $t = 2$, where both phase and amplitude errors grow and fail to recover to the correct trajectory.

\begin{figure}[!htb]
    \centering
    \includegraphics[width=0.85\textwidth]{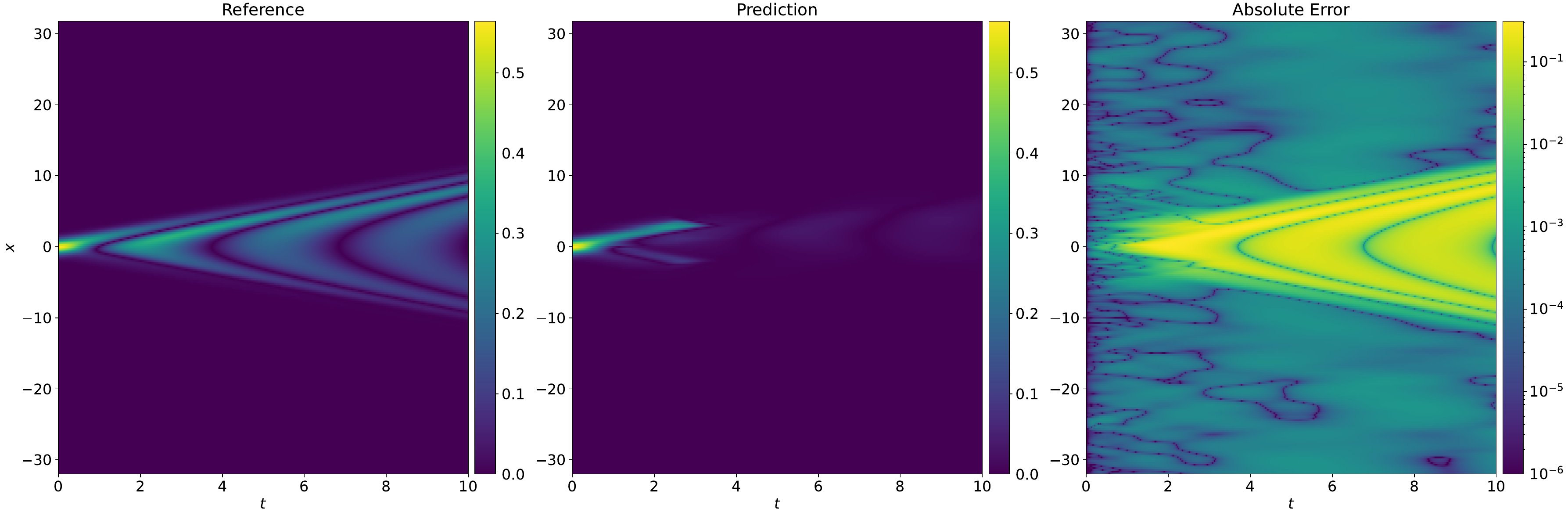}
    \vspace{-8pt}
    \caption{The prediction solution of PINNs for $\varepsilon=1.0,T=10.0$.}
    \label{fig:pro_fail_comp_long}
\end{figure}

\begin{figure}[!htbp]
     \centering  
     \subfigure[$\varepsilon=1.0,T=5.0$]{\includegraphics[width=0.80\textwidth]{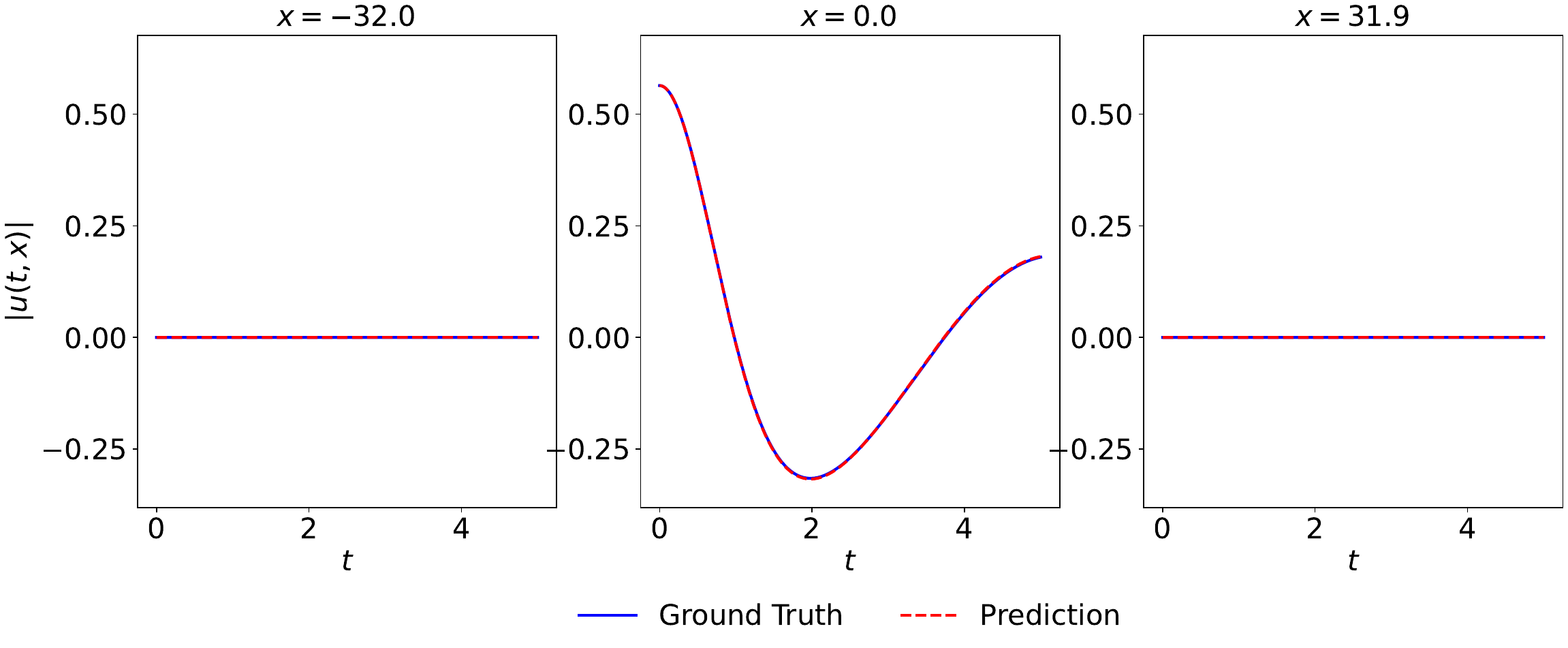}}\\
     \vspace{-4pt}
     \subfigure[$\varepsilon=1.0,T=10.0$]{\includegraphics[width=0.80\textwidth]{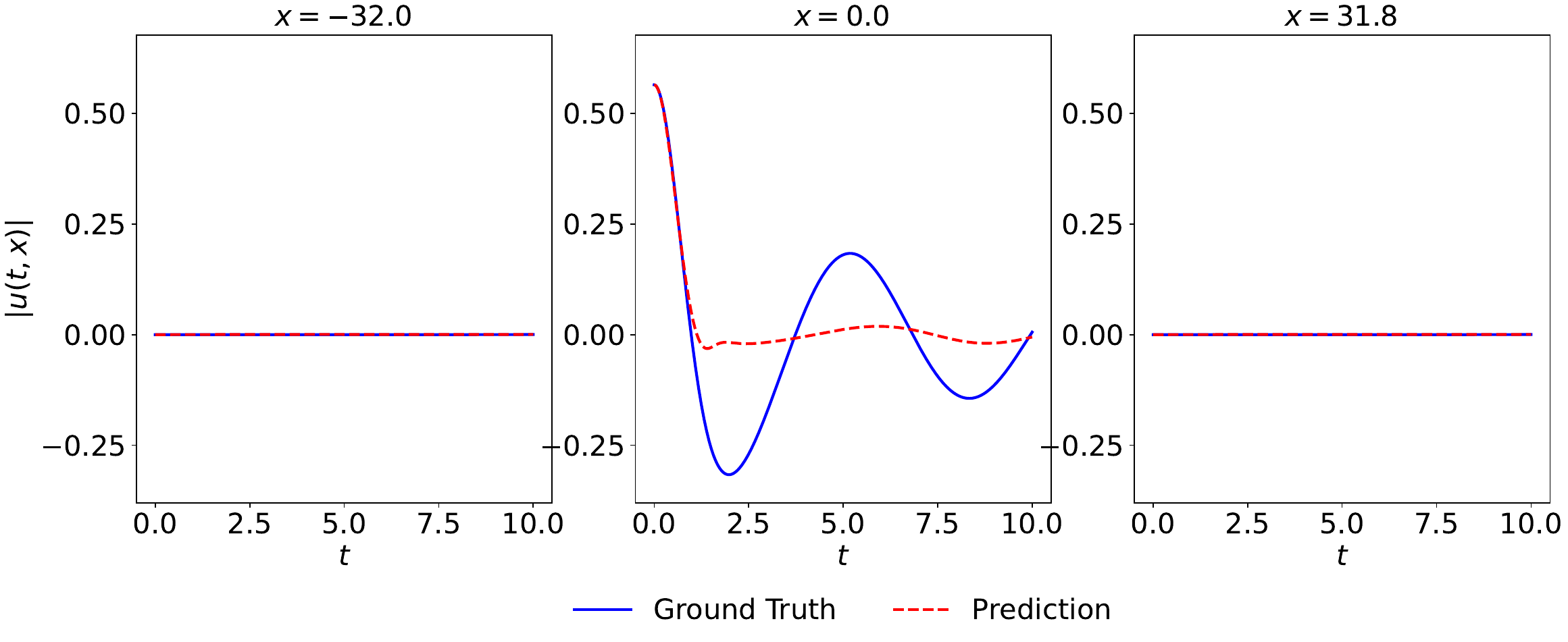}}\vspace{-8pt}
     \caption{The prediction solution of PINNs for $\varepsilon=1.0$ in different T.}  
     \label{fig:pro_fail_T}
\end{figure}

For $T = 10$, Figure~\ref{fig:pro_fail_comp_long} shows failure near $t = 3$.
This occurs even in the non-oscillatory case $\varepsilon = 1.0$. 
Figure~\ref{fig:pro_fail_T} shows progressive stagnation when $T$ increases from 5 to 10. 
These results identify missing causal and recursive constraints as the primary cause.
Time oscillation mainly accelerates the onset of propagation failure.

For reference, the key displayed relations in this subsection are \eqref{eq:auto-012}, \eqref{eq:auto-013}, \eqref{eq:auto-014}, \eqref{eq:auto-015}, \eqref{eq:auto-016}, \eqref{eq:auto-017}.
The continuation of this numbered set is \eqref{eq:auto-018}, \eqref{eq:auto-019}.

\section{Our proposed method}
\label{sec:method}

NeuralMD combines a WKB-inspired multiscale decomposition with physics-informed
learning.  The method has three components.  First, the rapidly oscillatory
carrier phases are represented analytically and a network learns the slowly
varying modulation.  Second, a separate network learns the dynamic remainder
subject to hard initial-data compatibility.  Third, a scalar, data-free
blending variable determines how much of the learned dynamic remainder should
be retained.  The blending variable is selected by the full NKGE residual and
is accompanied by an explicit near-oracle guarantee.  To reduce temporal
propagation failure during the first two stages, we use a stability-aware causal
gate together with multiscale time perturbations.  The gate is advanced only
when local loss gradients are mutually aligned.
It also requires a sufficiently small residual at the current temporal front.

\subsection{Multiscale decomposition by frequency}
\label{sec:MTI}

The WKB expansion, also called a modulated Fourier expansion, is a classical
tool for oscillatory differential equations
\citep{cohen2003modul,faou2013sob,hairer2006}.  For the NKGE, it separates the
rapid phase from slowly varying amplitudes through
\begin{equation}
\label{eq:auto-020}
  u(\bm{x},t)
  = \sum_{m\in\mathbb Z} e^{imt/\varepsilon^2}u_m(\bm{x},t).
\end{equation}
For sufficiently regular solutions, the derivatives of the modulated
coefficients remain bounded on fixed time intervals.  Retaining the leading
pair of modes gives
\begin{equation}
\label{eq:auto-021}
  u(\bm{x},t)
  = e^{it/\varepsilon^2}z(\bm{x},t)
  + e^{-it/\varepsilon^2}\overline{z(\bm{x},t)}
  + O(\varepsilon^2),
  \qquad \varepsilon\to0,
\end{equation}
where $z$ is complex-valued.  Under well-prepared initial data, the modulation
satisfies the nonlinear Schr\"odinger equation with wave operator (NLSW)
\citep{bao2012uni,bao2014uni},
\begin{equation}
\label{eq:auto-022}
\left\{
\begin{aligned}
  &2i\partial_tz+\varepsilon^2\partial_{tt}z-\Delta z
    +3\lambda|z|^2z=0,
    &&\bm{x}\in\mathbb R^d,\ t>0,\\
  &z(\bm{x},0)=z_0(\bm{x})
    :=\frac12\bigl(\phi_1(\bm{x})-i\phi_2(\bm{x})\bigr),\\
  &\partial_tz(\bm{x},0)=z_1(\bm{x})
    :=\frac{i}{2}\bigl(-\Delta z_0+3\lambda|z_0|^2z_0\bigr).
\end{aligned}
\right.
\end{equation}
Dropping the $\varepsilon^2\partial_{tt}z$ term yields the limiting NLSE
\citep{machihara2002,masmoudi2002},
\begin{equation}
\label{eq:auto-023}
\left\{
\begin{aligned}
  &2i\partial_tz-\Delta z+3\lambda|z|^2z=0,
    &&\bm{x}\in\mathbb R^d,\ t>0,\\
  &z(\bm{x},0)=z_0(\bm{x}).
\end{aligned}
\right.
\end{equation}
The corresponding leading-order reconstructions are
\begin{equation}
\label{eq:auto-024}
\begin{aligned}
  u_{\rm nlsw}(\bm{x},t)
  &=e^{it/\varepsilon^2}z_{\rm nlsw}(\bm{x},t)
    +e^{-it/\varepsilon^2}\overline{z_{\rm nlsw}(\bm{x},t)},\\
  u_{\rm nlse}(\bm{x},t)
  &=e^{it/\varepsilon^2}z_{\rm nlse}(\bm{x},t)
    +e^{-it/\varepsilon^2}\overline{z_{\rm nlse}(\bm{x},t)}.
\end{aligned}
\end{equation}
On a fixed time interval $[0,T]$, we use standard nonrelativistic-limit assumptions.
These include well-preparedness, regularity, and existence for the NKGE.
The estimates used in this work are summarized as
\begin{subequations}\label{eqn:modulation_error}
\begin{align}
  \|u(\cdot,t)-u_{\rm nlsw}(\cdot,t)\|_{H^1}
  &\le C_0\varepsilon^2,
  &&0\le t\le T,
  &&\phi_1,\phi_2\in H^2(\Omega),
  \\
  \|u(\cdot,t)-u_{\rm nlse}(\cdot,t)\|_{H^1}
  &\le (C_1+C_2T)\varepsilon^2,
  &&0\le t\le T,
  &&\phi_1,\phi_2\in H^3(\Omega),
\end{align}
\end{subequations}
where the constants may depend on $T$, $\Omega$, $\lambda$, and the dimension.
They may also depend on Sobolev norms of the limiting solution.
They are independent of $\varepsilon$ on the stated interval.
We do not use these estimates as global-in-time bounds.

\begin{figure}[!htb]
  \centering
  \subfigure[]{\label{subfig:Mfe}
    \includegraphics[width=0.42\textwidth]{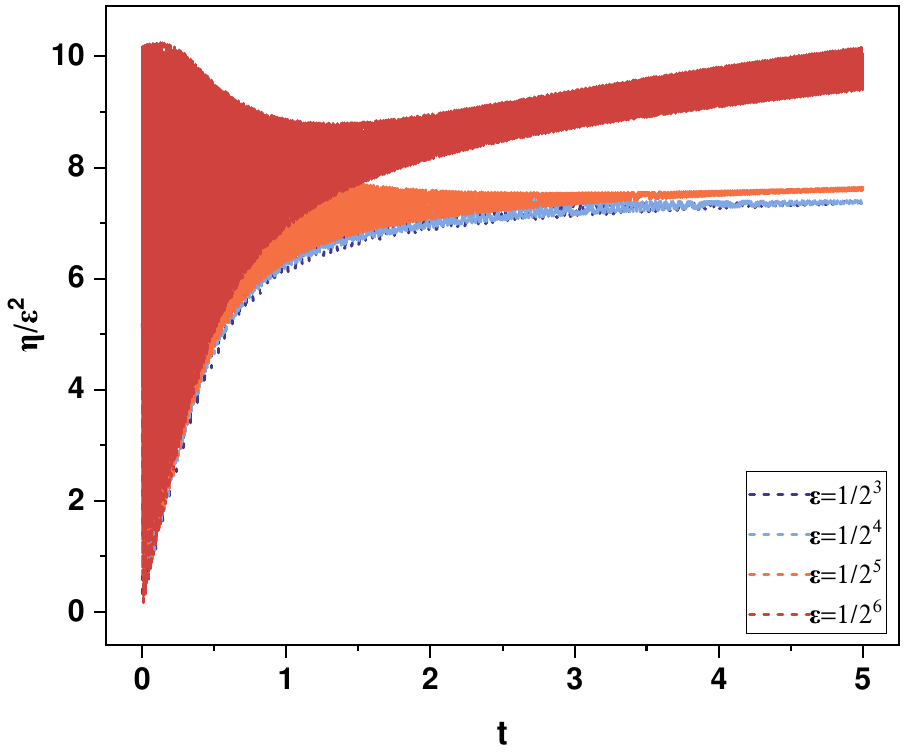}}
  \qquad
  \subfigure[]{\label{subfig:Mfo}
    \includegraphics[width=0.42\textwidth]{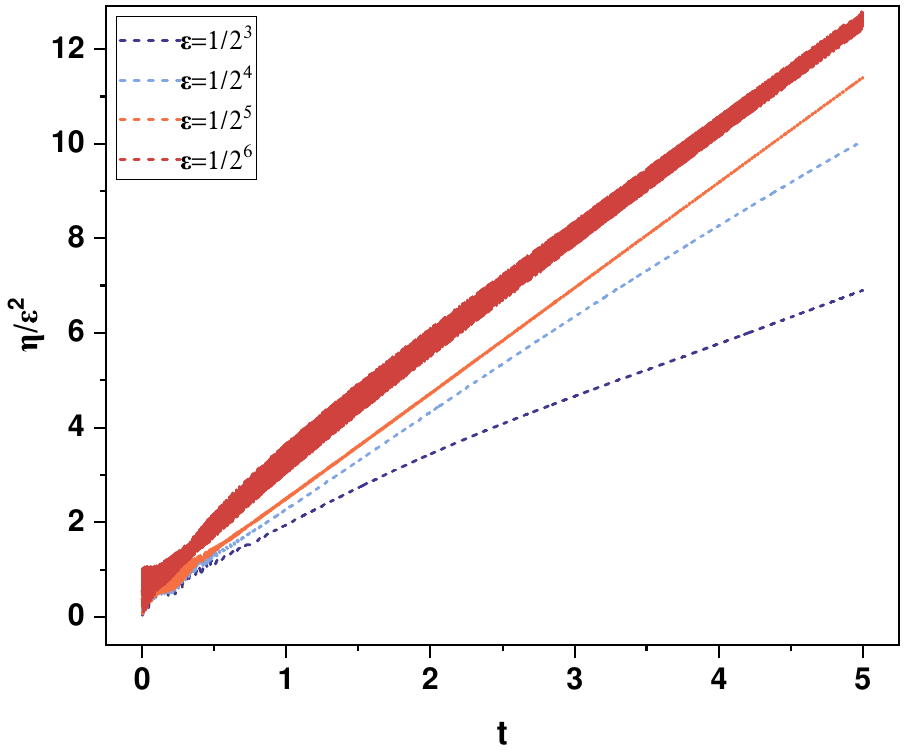}}
  \vspace{-6pt}
  \caption{Error convergence curves under different $\varepsilon$.
  (a) Convergence of the NKGE reconstruction to the NLSW reconstruction.
  (b) Convergence to the NLSE reconstruction with the expected time-dependent
  prefactor.}
  \label{fig:convergence}
\end{figure}

The multiscale time integrator (MTI) applies the same idea locally on each
interval $[t_n,t_{n+1}]$.  Writing $s=t-t_n$ and
\begin{equation}
\label{eq:auto-025}
  u(\bm{x},t_n+s)
  =e^{is/\varepsilon^2}z^n(\bm{x},s)
   +e^{-is/\varepsilon^2}\overline{z^n(\bm{x},s)}
   +r^n(\bm{x},s),
\end{equation}
leads to

\begin{align*}
  &2i\partial_sz^n+\varepsilon^2\partial_{ss}z^n-\Delta z^n
    +3\lambda|z^n|^2z^n=0,
  \\
  &\varepsilon^2\partial_{ss}r^n-\Delta r^n+\varepsilon^{-2}r^n
    +f_r(z^n,r^n;s)=0,
\end{align*}

with
\begin{equation}
\label{eq:auto-026}
\left\{
\begin{aligned}
  z^n(\bm{x},0)&=\frac12\bigl(\phi_1^n-i\phi_2^n\bigr),\\
  \partial_sz^n(\bm{x},0)
  &=\frac{i}{2}\bigl(-\Delta z^n(\bm{x},0)
    +3\lambda|z^n(\bm{x},0)|^2z^n(\bm{x},0)\bigr),\\
  r^n(\bm{x},0)&=0,\\
  \partial_sr^n(\bm{x},0)
  &=-\partial_sz^n(\bm{x},0)
    -\partial_s\overline{z^n(\bm{x},0)}.
\end{aligned}
\right.
\end{equation}
The nonlinear forcing is
\begin{equation}
\label{eq:fr}
\begin{aligned}
  f_r(z,r;s)
  ={}&\lambda e^{3is/\varepsilon^2}z^3
      +\lambda e^{-3is/\varepsilon^2}\overline z^{\,3}\\
   &+3\lambda\bigl(e^{2is/\varepsilon^2}z^2
      +e^{-2is/\varepsilon^2}\overline z^{\,2}\bigr)r\\
   &+3\lambda\bigl(e^{is/\varepsilon^2}z
      +e^{-is/\varepsilon^2}\overline z\bigr)r^2
      +6\lambda|z|^2r+\lambda r^3.
\end{aligned}
\end{equation}
For NeuralMD we use the same algebraic decomposition on
$Q_T=\Omega\times(0,T)$.  The identity is exact once the remainder is defined
by subtracting the modulated reconstruction from the NKGE solution.
Only the smallness claim is asymptotic.
It requires the regularity assumptions behind~\eqref{eqn:modulation_error}.

For reference, the key displayed relations in this subsection are \eqref{eq:auto-020}, \eqref{eq:auto-021}, \eqref{eq:auto-022}, \eqref{eq:auto-023}, \eqref{eq:auto-024}, \eqref{eq:auto-025}.
The continuation of this numbered set is \eqref{eq:auto-026}.

\subsection{Neural multiscale decomposition (NeuralMD)}
\label{sec:NeuralMD}

Define the NKGE operator
\begin{equation}
\label{eq:auto-027}
  \mathcal N_\varepsilon[v]
  :=\varepsilon^2\partial_{tt}v-\Delta v
    +\varepsilon^{-2}v+\lambda v^3.
\end{equation}
NeuralMD uses two pretrained networks and one scalar blending variable.  We use
hard initial-data lifts so that every member of the final blending family
satisfies the original initial conditions.

\begin{figure}[!htb]
    \centering
    \includegraphics[width=1.0\textwidth]{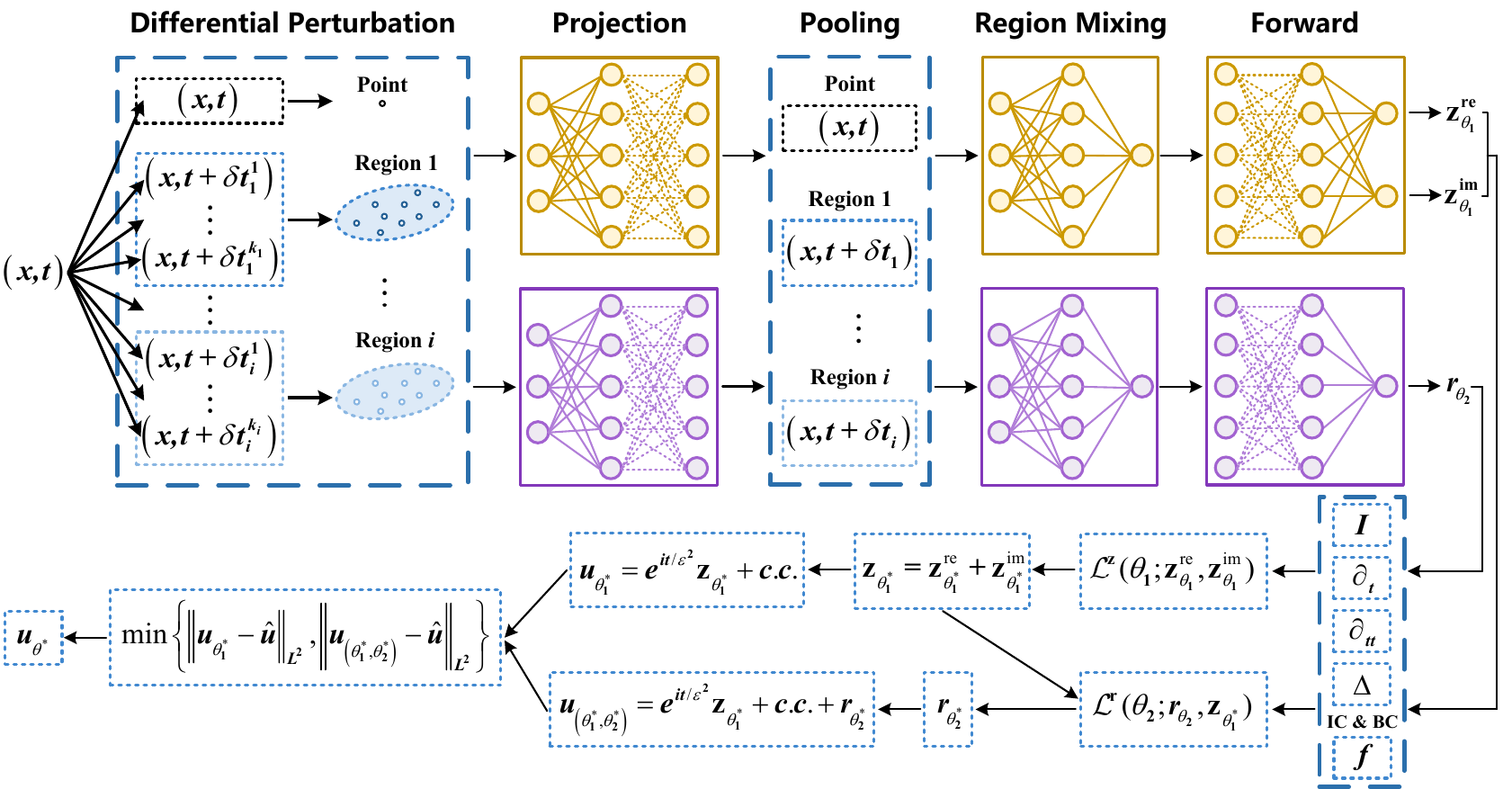}
    \vspace{-8pt}
    \caption{Overall architecture of NeuralMD. The single input point is perturbed in time to form multiscale temporal neighborhoods.
    The model gradients are propagated backward across multiple regions.}
    \label{fig:NeuralMD}
\end{figure}

For reference, the key displayed relations in this subsection are \eqref{eq:auto-027}.

\paragraph{Initial-data lift.}
Let
\begin{equation}
\label{eq:auto-028}
  \chi_0(t)=(1+t)e^{-t},\qquad
  \chi_1(t)=te^{-t},\qquad
  \chi_2(t)=1-e^{-t}-te^{-t}.
\end{equation}
Then
$\chi_0(0)=1$, $\chi_0'(0)=0$,
$\chi_1(0)=0$, $\chi_1'(0)=1$, and
$\chi_2(0)=\chi_2'(0)=0$.
For any unconstrained complex-valued network
$\widetilde z_{\theta_1}$, define
\begin{equation}
\label{eq:z_hard_lift}
  z_{\theta_1}(\bm{x},t)
  =\chi_0(t)z_0(\bm{x})
   +\chi_1(t)z_1(\bm{x})
   +\chi_2(t)\widetilde z_{\theta_1}(\bm{x},t).
\end{equation}
Thus $z_{\theta_1}(\bm{x},0)=z_0(\bm{x})$ and
$\partial_t z_{\theta_1}(\bm{x},0)=z_1(\bm{x})$ exactly.
Periodic boundary conditions are imposed through the periodic coordinate map
used by the network.  If a soft boundary treatment is used instead, the
boundary losses below are retained explicitly.

For reference, the key displayed relations in this subsection are \eqref{eq:auto-028}.

\paragraph{Pre-training stage I: modulation network.}
The NLSW residual is
\begin{equation}
\label{eq:Rz}
  \mathcal R_z[ z_{\theta_1}]
  =2i\partial_t z_{\theta_1}
   +\varepsilon^2\partial_{tt}z_{\theta_1}
   -\Delta z_{\theta_1}
   +3\lambda|z_{\theta_1}|^2z_{\theta_1}.
\end{equation}
The stage-I objective is
\begin{equation}
\label{eq:Lz}
  \mathcal L^z(\theta_1)
  =\frac1{N_f}\sum_{j=1}^{N_f}
    \bigl|\mathcal R_z[z_{\theta_1}](\bm{x}_j,t_j)\bigr|^2
   +\lambda_{\rm Bd}\mathcal L^z_{\rm Bd}(\theta_1).
\end{equation}
The complex residual is implemented by two real output channels.  The hard lift
in~\eqref{eq:z_hard_lift} removes the need to balance an additional initial-loss
weight.

The temporal symbol of the linear part of~\eqref{eq:Rz} is
\begin{equation}
\label{eq:auto-029}
  \sigma_z(k,\xi)=-2k-\varepsilon^2k^2+|\xi|^2.
\end{equation}
The physically relevant modulation band has $|k|=O(1)$; hence
$|\sigma_z(k,\xi)|=O(1)$ on bounded spatial bands, uniformly in
$\varepsilon$.  The network therefore learns a slow envelope rather than the
carrier frequency $O(\varepsilon^{-2})$.

For reference, the key displayed relations in this subsection are \eqref{eq:auto-029}.

\paragraph{Pre-training stage II: dynamic remainder network.}
Freeze $z_{\theta_1^\ast}$ and define the exact remainder initial data
\begin{equation}
\label{eq:auto-030}
  r_0(\bm{x})=0,
  \qquad
  r_1(\bm{x})=-z_1(\bm{x})-\overline{z_1(\bm{x})}.
\end{equation}
We split the learned remainder into a fixed initial-data correction and a
dynamic correction,
\begin{equation}
\label{eq:r_split}
\begin{aligned}
  r_{\rm ic}(\bm{x},t)
  &=r_0(\bm{x})\cos(t/\varepsilon^2)
    +\varepsilon^2r_1(\bm{x})\sin(t/\varepsilon^2),\\
  q_{\theta_2}(\bm{x},t)
  &=\chi_2(t)\widetilde q_{\theta_2}(\bm{x},t),\\
  r_{\theta_2}(\bm{x},t)
  &=r_{\rm ic}(\bm{x},t)+q_{\theta_2}(\bm{x},t).
\end{aligned}
\end{equation}
Consequently,
$q_{\theta_2}(\bm{x},0)=\partial_t q_{\theta_2}(\bm{x},0)=0$.
The oscillatory lift matches the homogeneous stiff Klein--Gordon oscillator.
It avoids an $O(\varepsilon^{-2})$ residual artifact from a slowly varying velocity lift.
The remainder residual is
\begin{equation}
\label{eq:auto-031}
  \mathcal R_r[q_{\theta_2}]
  =\varepsilon^2\partial_{tt}r_{\theta_2}
   -\Delta r_{\theta_2}
   +\varepsilon^{-2}r_{\theta_2}
   +f_r(z_{\theta_1^\ast},r_{\theta_2};t),
\end{equation}
and the stage-II objective is
\begin{equation}
\label{eq:Lr}
  \mathcal L^r(\theta_2)
  =\frac1{N_f}\sum_{j=1}^{N_f}
   \bigl|\mathcal R_r[q_{\theta_2}](\bm{x}_j,t_j)\bigr|^2
   +\lambda_{\rm Bd}\mathcal L^r_{\rm Bd}(\theta_2).
\end{equation}

The principal linear symbol of the remainder equation is
\begin{equation}
\label{eq:auto-032}
  \sigma_\varepsilon(k,\xi)
  =-\varepsilon^2k^2+|\xi|^2+\varepsilon^{-2}.
\end{equation}
The forcing~\eqref{eq:fr} generates bands centered near
$k=m\varepsilon^{-2}$, $m\in\{\pm1,\pm2,\pm3\}$.  These bands are difficult
for standard networks because of their large temporal frequencies.  However,
$|\sigma_\varepsilon|$ is not uniformly large on every high-frequency mode.
Near the characteristic set, it may be small because temporal and mass terms cancel.
This distinction is essential for the remainder
selection theory below; the PDE residual strongly penalizes off-characteristic
artifacts, not all high-frequency physical components.

For reference, the key displayed relations in this subsection are \eqref{eq:auto-030}, \eqref{eq:auto-031}, \eqref{eq:auto-032}.

\paragraph{Error criterion.}
\label{sec:soft_blending}

After stages I and II, NeuralMD uses a physics-aware error criterion to decide
whether the learned remainder should be retained.  The criterion keeps the
original initial data exactly by blending only the dynamic part of the
remainder.  Set
\begin{equation}
\label{eq:auto-033}
  u_{\rm base}(\bm{x},t)
  =e^{it/\varepsilon^2}z_{\theta_1^\ast}(\bm{x},t)
   +e^{-it/\varepsilon^2}\overline{z_{\theta_1^\ast}(\bm{x},t)}
   +r_{\rm ic}(\bm{x},t).
\end{equation}
The final one-parameter family is
\begin{equation}
\label{eq:blend_family}
  u_\omega(\bm{x},t)
  =u_{\rm base}(\bm{x},t)+\omega q_{\theta_2^\ast}(\bm{x},t),
  \qquad 0\le\omega\le1.
\end{equation}
Because the dynamic correction and its first time derivative vanish at
$t=0$, every $u_\omega$ satisfies
\begin{equation}
\label{eq:auto-034}
  u_\omega(\bm{x},0)=\phi_1(\bm{x}),
  \qquad
  \partial_tu_\omega(\bm{x},0)=\varepsilon^{-2}\phi_2(\bm{x}).
\end{equation}
If periodicity is hard encoded, every member of the family also satisfies the
boundary condition.  This avoids the incompatibility caused by multiplying the
entire remainder, including its nonzero initial velocity, by a free scalar.

\begin{proposition}[Initial-data invariance of dynamic blending]
\label{prop:blend_ic_invariance}
Assume $u_{\rm base}$ satisfies the prescribed initial data and boundary
condition, and assume
$q_{\theta_2^\ast}(\bm{x},0)=\partial_tq_{\theta_2^\ast}(\bm{x},0)=0$.
Then every $u_\omega$ in~\eqref{eq:blend_family} satisfies the same initial
data.  If a shared hard lift imposes periodic, Dirichlet, or Neumann data,
then every $u_\omega$ satisfies the same boundary condition.
\end{proposition}

\begin{proof}
The initial identities follow by substituting $t=0$ in
\eqref{eq:blend_family}.  The boundary statement follows from linearity of the
boundary trace for the homogeneous dynamic correction.  Thus the scalar
selection step changes only the dynamic remainder amplitude, not the
constraint class.
\end{proof}

The population physics objective is
\begin{equation}
\label{eq:L_blend}
  \mathcal J_{\rm phys}(\omega)
  =\|\mathcal N_\varepsilon[u_\omega]\|_{\mathcal Y}^2,
  \qquad \mathcal Y=L^2(Q_T,\mu),
\end{equation}
where $\mu$ is the normalized collocation measure used for the population
residual norm.
and the empirical selector uses an independent collocation set
$\widehat P_f=\{(\widehat{\bm{x}}_j,\widehat t_j)\}_{j=1}^{N_s}$,
\begin{equation}
\label{eq:L_blend_empirical}
  \widehat{\mathcal J}_{N_s}(\omega)
  =\frac1{N_s}\sum_{j=1}^{N_s}
   \left|\mathcal N_\varepsilon[u_\omega]
   (\widehat{\bm{x}}_j,\widehat t_j)\right|^2.
\end{equation}
Using an independent set reduces selection overfitting.  If the initial or
boundary conditions are imposed softly, their losses must also be included in
both~\eqref{eq:L_blend} and~\eqref{eq:L_blend_empirical}.

For the cubic NKGE nonlinearity, the scalar objective is explicit.  With
$u_0=u_{\rm base}$ and $q=q_{\theta_2^\ast}$, define
\begin{align*}
  \mathcal L_{u_0}q
  &:=\varepsilon^2\partial_{tt}q-\Delta q+\varepsilon^{-2}q
    +3\lambda u_0^2q,\\
  B_2(u_0,q)&:=3u_0q^2,\qquad
  B_3(q):=q^3 .
\end{align*}
Then
\begin{equation}
\label{eq:auto-035}
  \mathcal N_\varepsilon[u_0+\omega q]
  =\mathcal N_\varepsilon[u_0]
   +\omega\mathcal L_{u_0}q
   +\lambda\omega^2B_2(u_0,q)
   +\lambda\omega^3B_3(q).
\end{equation}
Consequently $\mathcal J_{\rm phys}$ is a degree-at-most-six polynomial in
$\omega$.  It is continuous on $[0,1]$.
Therefore a global scalar minimizer exists, although convexity is not assumed.
Its derivatives are
\begin{align*}
  \mathcal J_{\rm phys}'(\omega)
  &=2
    \langle R_\omega,R_\omega'\rangle_{\mathcal Y},\\
  \mathcal J_{\rm phys}''(\omega)
  &=2\|R_\omega'\|_{\mathcal Y}^2
    +2
    \langle R_\omega,R_\omega''\rangle_{\mathcal Y},
\end{align*}
where $R_\omega=\mathcal N_\varepsilon[u_0+\omega q]$,
$R_\omega'=\mathcal L_{u_0}q+2\lambda\omega B_2+3\lambda\omega^2B_3$, and
$R_\omega''=2\lambda B_2+6\lambda\omega B_3$.

For reference, the key displayed relations in this subsection are \eqref{eq:auto-033}, \eqref{eq:auto-034}, \eqref{eq:auto-035}.

\paragraph{Exact oracle weight.}
Let $u^\dagger$ be the exact NKGE solution and let $\mathcal H_{\rm sol}$ be a
real solution-error Hilbert space used for evaluation.  For an $L^2$ error
metric one may take $\mathcal H_{\rm sol}=L^2(Q_T,\mu)$.  For an energy metric,
one may use the fixed-time norm associated with
$H^1(\Omega)\times L^2(\Omega)$.  This requires finite solution and remainder norms.
We keep the notation abstract because the oracle
argument only uses Hilbert-space geometry.
Define
\begin{equation}
\label{eq:auto-036}
  \rho=u^\dagger-u_{\rm base},
  \qquad q=q_{\theta_2^\ast}.
\end{equation}
The oracle error is
\begin{equation}
\label{eq:E_omega}
  \mathcal E_{\rm or}(\omega)
  =\|u^\dagger-u_\omega\|_{\mathcal H_{\rm sol}}^2
  =\|\rho-\omega q\|_{\mathcal H_{\rm sol}}^2.
\end{equation}

\begin{proposition}[Exact oracle blending weight]\label{prop:oracle_weight}
Assume $q\ne0$.  The unique minimizer of~\eqref{eq:E_omega} over $[0,1]$ is
\begin{equation}
\label{eq:oracle_exact}
  \omega_{\rm or}
  =\Pi_{[0,1]}\!\left(
    \frac{\langle\rho,q\rangle_{\mathcal H_{\rm sol}}}
         {\|q\|_{\mathcal H_{\rm sol}}^2}
  \right),
\end{equation}
where $\Pi_{[0,1]}$ denotes Euclidean projection.  Moreover,
$\omega_{\rm or}\in(0,1)$ if and only if
\begin{equation}
\label{eq:oracle_interior_condition}
  0<\langle\rho,q\rangle_{\mathcal H_{\rm sol}}<\|q\|_{\mathcal H_{\rm sol}}^2.
\end{equation}
Writing $a_q=\langle\rho,q\rangle_{\mathcal H_{\rm sol}}$ and
$d_q=\|q\|_{\mathcal H_{\rm sol}}^2$, the KKT alternatives are
\begin{equation}
\label{eq:oracle_kkt_cases}
  \omega_{\rm or}=0\ \Longleftrightarrow\ a_q\le0,\qquad
  \omega_{\rm or}=a_q/d_q\ \Longleftrightarrow\ 0<a_q<d_q,\qquad
  \omega_{\rm or}=1\ \Longleftrightarrow\ a_q\ge d_q .
\end{equation}
The minimum oracle error is
\begin{equation}
\label{eq:oracle_min_error}
  \mathcal E_{\rm or}(\omega_{\rm or})
  =
  \begin{cases}
    \|\rho\|_{\mathcal H_{\rm sol}}^2, & a_q\le0,\\
    \|\rho\|_{\mathcal H_{\rm sol}}^2-a_q^2/d_q, & 0<a_q<d_q,\\
    \|\rho-q\|_{\mathcal H_{\rm sol}}^2, & a_q\ge d_q .
  \end{cases}
\end{equation}
If the minimizer is interior, namely $0<a_q<d_q$, then
\begin{equation}
\label{eq:strict_improve}
  \mathcal E_{\rm or}(\omega_{\rm or})
  <\min\{\mathcal E_{\rm or}(0),\mathcal E_{\rm or}(1)\}.
\end{equation}
More precisely, for an interior minimizer,
\begin{equation}
\label{eq:oracle_improvement_amounts}
  \mathcal E_{\rm or}(0)-\mathcal E_{\rm or}(\omega_{\rm or})
  =\frac{a_q^2}{d_q},
  \qquad
  \mathcal E_{\rm or}(1)-\mathcal E_{\rm or}(\omega_{\rm or})
  =\frac{(d_q-a_q)^2}{d_q}.
\end{equation}
If the minimizer is on the boundary, continuous blending is equal to, rather
than strictly better than, the best binary choice.
\end{proposition}

\begin{proof}
Expanding~\eqref{eq:E_omega} gives
\begin{equation}
\label{eq:auto-037}
  \mathcal E_{\rm or}(\omega)
  =\|\rho\|_{\mathcal H_{\rm sol}}^2
   -2\omega\langle\rho,q\rangle_{\mathcal H_{\rm sol}}
   +\omega^2\|q\|_{\mathcal H_{\rm sol}}^2.
\end{equation}
The unconstrained minimizer is
$\langle\rho,q\rangle_{\mathcal H_{\rm sol}}/\|q\|_{\mathcal H_{\rm sol}}^2$; projection gives
\eqref{eq:oracle_exact}.  The one-dimensional KKT conditions give
\eqref{eq:oracle_kkt_cases}, and evaluating the quadratic at the projected
minimizer gives~\eqref{eq:oracle_min_error}.  Strict convexity and the
location of the unconstrained minimizer give
\eqref{eq:oracle_interior_condition}--\eqref{eq:oracle_improvement_amounts}.
\end{proof}

Writing $q=\rho+\delta$ makes the signal--error trade-off explicit:
\begin{equation}
\label{eq:auto-038}
  \omega_{\rm or}
  =\Pi_{[0,1]}\!\left(
  \frac{\|\rho\|_{\mathcal H_{\rm sol}}^2
        +\langle\rho,\delta\rangle_{\mathcal H_{\rm sol}}}
       {\|\rho+\delta\|_{\mathcal H_{\rm sol}}^2}
  \right).
\end{equation}
Only when the cross term is negligible does this reduce to the familiar
Wiener-type expression
$\|\rho\|^2/(\|\rho\|^2+\|\delta\|^2)$.
Indeed, if
$a=\|\rho\|_{\mathcal H_{\rm sol}}^2$,
$b=\langle\rho,\delta\rangle_{\mathcal H_{\rm sol}}$, and
$c=\|\delta\|_{\mathcal H_{\rm sol}}^2$, then the unconstrained weight is
\begin{equation}
\label{eq:auto-039}
  \omega_0=\frac{a+b}{a+2b+c}.
\end{equation}
Whenever $0<\omega_0<1$, the exact minimum error is
\begin{equation}
\label{eq:auto-040}
  \mathcal E_{\rm or}(\omega_0)
  =\frac{ac-b^2}{a+2b+c}.
\end{equation}
Thus the cross term is part of the optimal trade-off; the orthogonal
signal-noise formula is only the special case $b=0$.

\begin{remark}[Asymptotic smallness is not a removal criterion]
If the learned dynamic remainder is exact, $q=\rho$, then
$\omega_{\rm or}=1$ regardless of how small $\|\rho\|_{\mathcal H_{\rm sol}}$ is.
In particular, $\|\rho\|=O(\varepsilon^2)$ alone does not imply that the
accuracy-optimal weight tends to zero.  Complete removal requires relative
learning error, off-characteristic artifacts, or an explicit retention-cost
term.
\end{remark}

\begin{corollary}[A sufficient condition for a transition-regime interior point]
\label{cor:oracle_interior}
Suppose $\rho\ne0$, $\delta\ne0$, and for some $0\le\beta<1$,
\begin{equation}
\label{eq:cross_term_assumption}
  |\langle\rho,\delta\rangle_{\mathcal H_{\rm sol}}|
  \le\beta\min\{\|\rho\|_{\mathcal H_{\rm sol}}^2,
                 \|\delta\|_{\mathcal H_{\rm sol}}^2\}.
\end{equation}
Then $0<\omega_{\rm or}<1$ and the strict inequality
\eqref{eq:strict_improve} holds.
\end{corollary}

\begin{proof}
Let $a=\|\rho\|_{\mathcal H_{\rm sol}}^2$, $c=\|\delta\|_{\mathcal H_{\rm sol}}^2$, and
$b=\langle\rho,\delta\rangle_{\mathcal H_{\rm sol}}$.  Assumption
\eqref{eq:cross_term_assumption} gives $a+b>0$ and $c+b>0$.
Since $\langle\rho,q\rangle=a+b$ and
$\|q\|^2-\langle\rho,q\rangle=c+b$, condition
\eqref{eq:oracle_interior_condition} follows.
\end{proof}

For reference, the key displayed relations in this subsection are \eqref{eq:auto-036}, \eqref{eq:auto-037}, \eqref{eq:auto-038}, \eqref{eq:auto-039}, \eqref{eq:auto-040}.

\paragraph{Near-oracle optimality of the PDE-residual selector.}
Set
\begin{equation}
\label{eq:auto-041}
  V=\operatorname{span}\{\rho,q\},
  \qquad
  A_\varepsilon=\mathcal N_\varepsilon'[u^\dagger].
\end{equation}
The following assumptions are restricted to the two-dimensional blending
subspace.  They are therefore weaker than global coercivity on the full
neural-network hypothesis class.

The linearized residual objective on this subspace has an analytic minimizer.
For $v_\omega=\omega q-\rho$,
\begin{equation}
\label{eq:auto-042}
  \omega_{\rm lin}
  =
  \Pi_{[0,1]}\!\left(
  \frac{\langle A_\varepsilon\rho,A_\varepsilon q\rangle_{\mathcal Y}}
       {\|A_\varepsilon q\|_{\mathcal Y}^2}
  \right).
\end{equation}
Assumption~\ref{assum:metric} states when this residual metric is stable.
The comparison uses the subspace $V=\operatorname{span}\{\rho,q\}$.
In practice, audit the perturbation using two $2\times2$ Gram matrices:
\begin{equation}
\label{eq:auto-043}
[\langle v_i,v_j\rangle_{\mathcal H_{\rm sol}}],
\qquad
[\langle A_\varepsilon v_i,A_\varepsilon v_j\rangle_{\mathcal Y}].
\end{equation}
Here $\{v_1,v_2\}=\{\rho,q\}$.

\begin{assumption}[Restricted metric compatibility]\label{assum:metric}
There exist $\kappa>0$ and $0\le\eta_{\rm op}<1$ such that
\begin{equation}
\label{eq:metric_compatibility}
  \left|
  \langle A_\varepsilon v,A_\varepsilon w\rangle_{\mathcal Y}
  -\kappa\langle v,w\rangle_{\mathcal H_{\rm sol}}
  \right|
  \le\eta_{\rm op}\kappa
  \|v\|_{\mathcal H_{\rm sol}}\|w\|_{\mathcal H_{\rm sol}},
  \qquad v,w\in V.
\end{equation}
\end{assumption}

\begin{assumption}[Controlled nonlinear remainder]\label{assum:taylor}
Let $s>d/2$ and assume that $u^\dagger$ is bounded in
$L^\infty(0,T;H^s(\Omega))$.  Assume the same bound along the path
$\{u^\dagger+\omega q-\rho:0\le\omega\le1\}$.  The corresponding time derivatives
needed for the residual norm in $\mathcal Y$ are also assumed finite.
There exists $L_{\rm nl}>0$ such that,
for every $v=\omega q-\rho$, $0\le\omega\le1$,
\begin{equation}
\label{eq:auto-044}
  \|\mathcal N_\varepsilon[u^\dagger+v]-A_\varepsilon v\|_{\mathcal Y}
  \le L_{\rm nl}\|v\|_{\mathcal H_{\rm sol}}^2.
\end{equation}
\end{assumption}
For the cubic NKGE this assumption is the norm form of the algebraic identity
\begin{equation}
\label{eq:cubic_exact_remainder}
  \mathcal N_\varepsilon[u^\dagger+v]-A_\varepsilon v
  =\lambda\bigl(3u^\dagger v^2+v^3\bigr),
\end{equation}
so $L_{\rm nl}$ is controlled by Sobolev product estimates.
It also depends on the fixed time interval, $R_V$, and the stated bounds.  It may
depend on $\varepsilon$ through the residual norm and the regularity constants;
no uniform-in-$\varepsilon$ nonlinear stability is asserted here.  The
quadratic and cubic terms in \eqref{eq:cubic_exact_remainder} are the only
source of the nonlinear perturbation $\Delta_V$ below.

\begin{assumption}[Sampling and optimization accuracy]\label{assum:sampling_opt}
With probability at least $1-\delta_s$,
\begin{equation}
\label{eq:sampling_uniform}
  \sup_{\omega\in[0,1]}
  |\widehat{\mathcal J}_{N_s}(\omega)-\mathcal J_{\rm phys}(\omega)|
  \le\eta_{\rm samp}.
\end{equation}
The returned scalar $\widehat\omega$ satisfies
\begin{equation}
\label{eq:optimization_tolerance}
  \widehat{\mathcal J}_{N_s}(\widehat\omega)
  \le\inf_{\omega\in[0,1]}
      \widehat{\mathcal J}_{N_s}(\omega)+\eta_{\rm opt}.
\end{equation}
One concrete sufficient condition is as follows.  Let
$Z_\omega(\bm{x},t)=|\mathcal N_\varepsilon[u_\omega](\bm{x},t)|^2$ and assume
$0\le Z_\omega\le M_s$ and
$|Z_\omega-Z_{\omega'}|\le L_{\rm path}|\omega-\omega'|$ uniformly on $Q_T$.  For
frozen pretrained networks, this is an a posteriori boundedness assumption.
It applies only on the selector path.
It may be replaced by a sub-exponential concentration condition for unbounded residuals.  A
$\tau$-net of $[0,1]$ and Hoeffding's inequality give
\begin{equation}
\label{eq:auto-045}
  \eta_{\rm samp}
  \le
  L_{\rm path}\tau
  +M_s\sqrt{
    \frac{\log\!\left(2\lceil1/\tau\rceil/\delta_s\right)}
         {2N_s}} .
\end{equation}
If scalar search uses grid width $h_{\rm grid}$ and refinement tolerance
$\tau_{\rm ref}$, the following suffices.
Assume $\widehat{\mathcal J}_{N_s}$ is $L_{\rm emp}$-Lipschitz and take
\begin{equation}
\label{eq:auto-046}
  \eta_{\rm opt}\le L_{\rm emp}\left(\frac{h_{\rm grid}}2+\tau_{\rm ref}\right).
\end{equation}
\end{assumption}

Define the radius of the blending path by
\begin{equation}
\label{eq:auto-047}
  R_V=\sup_{\omega\in[0,1]}
      \|\omega q-\rho\|_{\mathcal H_{\rm sol}}
  \le\|q\|_{\mathcal H_{\rm sol}}+\|\rho\|_{\mathcal H_{\rm sol}},
\end{equation}
and set
\begin{equation}
\label{eq:auto-048}
\begin{aligned}
  \Delta_V={}&
   \eta_{\rm op}\kappa R_V^2
   +2L_{\rm nl}\sqrt{\kappa(1+\eta_{\rm op})}\,R_V^3
   +L_{\rm nl}^2R_V^4.
\end{aligned}
\end{equation}

\begin{theorem}[Finite-sample near-oracle remainder selection]
\label{thm:blend}
Under Assumptions~\ref{assum:metric}--\ref{assum:sampling_opt}, with
probability at least $1-\delta_s$,
\begin{equation}
\label{eq:near_opt}
  \mathcal E_{\rm or}(\widehat\omega)
  \le\mathcal E_{\rm or}(\omega_{\rm or})
  +\frac{2\Delta_V+2\eta_{\rm samp}+\eta_{\rm opt}}{\kappa}.
\end{equation}
If $q\ne0$, then
\begin{equation}
\label{eq:weight_near_opt}
  |\widehat\omega-\omega_{\rm or}|
  \le
  \left(
  \frac{2\Delta_V+2\eta_{\rm samp}+\eta_{\rm opt}}
       {\kappa\|q\|_{\mathcal H_{\rm sol}}^2}
  \right)^{1/2}.
\end{equation}
Thus the residual selector is near-oracle only under four controlled effects.
The linearized operator must approximately preserve the error metric on $V$.
The nonlinear Taylor remainder must be small.
The collocation estimate and scalar optimization must also be accurate.
\end{theorem}

\begin{proof}
Let $v_\omega=u_\omega-u^\dagger=\omega q-\rho$.  Assumption~\ref{assum:metric} gives
\begin{equation}
\label{eq:auto-049}
  \left|\|A_\varepsilon v_\omega\|_{\mathcal Y}^2
  -\kappa\|v_\omega\|_{\mathcal H_{\rm sol}}^2\right|
  \le\eta_{\rm op}\kappa R_V^2.
\end{equation}
Writing
$\mathcal N_\varepsilon[u^\dagger+v_\omega]
=A_\varepsilon v_\omega+\mathcal Q(v_\omega)$,
Assumption~\ref{assum:taylor} and
$\|A_\varepsilon v_\omega\|_{\mathcal Y}
\le\sqrt{\kappa(1+\eta_{\rm op})}R_V$ yield
\begin{equation}
\label{eq:uniform_proxy_gap}
  |\mathcal J_{\rm phys}(\omega)
    -\kappa\mathcal E_{\rm or}(\omega)|
  \le\Delta_V,
  \qquad 0\le\omega\le1.
\end{equation}
Using~\eqref{eq:sampling_uniform},
\eqref{eq:optimization_tolerance}, and~\eqref{eq:uniform_proxy_gap},
\begin{align*}
  \kappa\mathcal E_{\rm or}(\widehat\omega)
  &\le\mathcal J_{\rm phys}(\widehat\omega)+\Delta_V\\
  &\le\widehat{\mathcal J}_{N_s}(\widehat\omega)
      +\eta_{\rm samp}+\Delta_V\\
  &\le\widehat{\mathcal J}_{N_s}(\omega_{\rm or})
      +\eta_{\rm opt}+\eta_{\rm samp}+\Delta_V\\
  &\le\mathcal J_{\rm phys}(\omega_{\rm or})
      +2\eta_{\rm samp}+\eta_{\rm opt}+\Delta_V\\
  &\le\kappa\mathcal E_{\rm or}(\omega_{\rm or})
      +2\Delta_V+2\eta_{\rm samp}+\eta_{\rm opt}.
\end{align*}
This proves~\eqref{eq:near_opt}.  Since
$\mathcal E_{\rm or}$ is $2\|q\|_{\mathcal H_{\rm sol}}^2$-strongly convex and
$\omega_{\rm or}$ is its constrained minimizer,
\begin{equation}
\label{eq:auto-050}
  \mathcal E_{\rm or}(\omega)-\mathcal E_{\rm or}(\omega_{\rm or})
  \ge\|q\|_{\mathcal H_{\rm sol}}^2|\omega-\omega_{\rm or}|^2.
\end{equation}
Combining this inequality with~\eqref{eq:near_opt} proves
\eqref{eq:weight_near_opt}.
\end{proof}

Theorem~\ref{thm:blend} is deliberately conditional.  In particular, the PDE
residual is not asserted to be a perfect proxy for solution error.  The
restricted Gram matrices associated with~\eqref{eq:metric_compatibility}, the
uniform empirical gap in~\eqref{eq:sampling_uniform}, and the scalar objective
curve should be reported in the experiments.

For reference, the key displayed relations in this subsection are \eqref{eq:auto-041}, \eqref{eq:auto-042}, \eqref{eq:auto-043}, \eqref{eq:auto-044}, \eqref{eq:auto-045}, \eqref{eq:auto-046}.
The continuation of this numbered set is \eqref{eq:auto-047}, \eqref{eq:auto-048}, \eqref{eq:auto-049}, \eqref{eq:auto-050}.

\paragraph{Accuracy-aware versus cost-aware removal.}
The weight $\omega_{\rm or}$ minimizes solution error, whereas deciding whether
a remainder is worth retaining may also involve computational cost.  We use the
convex retention-cost surrogate
\begin{equation}
\label{eq:cost_aware_selector}
  \mathcal J_{\rm cost}(\omega)
  =\mathcal J_{\rm phys}(\omega)+\lambda_{\rm comp}\omega,
  \qquad \lambda_{\rm comp}\ge0.
\end{equation}
If $\mathcal J_{\rm phys}$ is convex on $[0,1]$, the KKT condition gives
\begin{equation}
\label{eq:auto-051}
  \omega_{\rm cost}=0
  \quad\Longleftrightarrow\quad
  \mathcal J_{\rm phys}'(0^+)+\lambda_{\rm comp}\ge0.
\end{equation}
For the oracle quadratic with a linear retention cost,
\begin{equation}
\label{eq:auto-052}
  \min_{0\le\omega\le1}
  \left\{\|\rho-\omega q\|_{\mathcal H_{\rm sol}}^2
        +\lambda_{\rm comp}\omega\right\},
\end{equation}
the exact cost-aware weight is
\begin{equation}
\label{eq:auto-053}
  \omega_{\rm cost,or}
  =
  \Pi_{[0,1]}\!\left(
  \frac{\langle\rho,q\rangle_{\mathcal H_{\rm sol}}-\lambda_{\rm comp}/2}
       {\|q\|_{\mathcal H_{\rm sol}}^2}
  \right).
\end{equation}
Hence complete removal is selected by this quadratic criterion precisely when
$\langle\rho,q\rangle_{\mathcal H_{\rm sol}}\le\lambda_{\rm comp}/2$.
For a literal on/off inference cost, one may instead use
$\lambda_{\rm comp}\mathbf 1_{\{\omega>0\}}$ and compare the global scalar
objective values directly.  Thus complete removal is optimal only for a specified
cost-accuracy objective.  Selecting
$\omega=0$ permits pruning the remainder network at inference.
It does not remove the Stage-II training cost.
Avoiding that cost requires a pilot or early-stopping policy.

For reference, the key displayed relations in this subsection are \eqref{eq:auto-051}, \eqref{eq:auto-052}, \eqref{eq:auto-053}.

\paragraph{Off-characteristic suppression in the nonrelativistic limit.}
Let
\begin{equation}
\label{eq:auto-054}
  \Gamma_{\varepsilon,\kappa}^{\rm off}
  =\{(k,\xi):|\sigma_\varepsilon(k,\xi)|
    \ge\kappa\varepsilon^{-2}\},
  \qquad \kappa>0,
\end{equation}
and let $P_{\rm off}$ denote the corresponding spectral projector.  The next
result concerns learned artifacts with non-vanishing mass away from the
Klein--Gordon characteristic set.

\begin{corollary}[Linearized off-characteristic suppression]
\label{cor:nonrel}
Let $F_\varepsilon=\mathcal N_\varepsilon[u_{\rm base}]$ and
$A_{\rm base}=\mathcal N_\varepsilon'[u_{\rm base}]$.  Suppose
\begin{equation}
\label{eq:off_characteristic_assumption}
  \|F_\varepsilon\|_{\mathcal Y}\le C_F,
  \qquad
  \|P_{\rm off}q\|_{\mathcal H_{\rm sol}}\ge c_0,
  \qquad
  \|A_{\rm base}q\|_{\mathcal Y}
  \ge c_A\varepsilon^{-2}\|P_{\rm off}q\|_{\mathcal H_{\rm sol}},
\end{equation}
with constants independent of $\varepsilon$.  The minimizer of the linearized
selector
\begin{equation}
\label{eq:auto-055}
  \mathcal J_{\rm lin}(\omega)
  =\|F_\varepsilon+\omega A_{\rm base}q\|_{\mathcal Y}^2
\end{equation}
satisfies
\begin{equation}
\label{eq:auto-056}
  0\le\omega_{\rm lin}
  \le\frac{C_F}{c_Ac_0}\varepsilon^2.
\end{equation}
Hence a non-vanishing off-characteristic artifact is suppressed as
$\varepsilon\to0$.  Moreover, suppose
\begin{equation}
\label{eq:off_characteristic_empirical_gap}
  \sup_{\omega\in[0,1]}
  |\widehat{\mathcal J}_{N_s}(\omega)-\mathcal J_{\rm lin}(\omega)|
  \le\eta_{\rm sel}
\end{equation}
and the scalar optimizer has tolerance $\eta_{\rm opt}$.  Then
\begin{equation}
\label{eq:off_characteristic_empirical_bound}
  \widehat\omega
  \le
  \frac{\varepsilon^2}{c_Ac_0}
  \left(C_F+\sqrt{2\eta_{\rm sel}+\eta_{\rm opt}}\right).
\end{equation}
This conclusion does not apply to a genuine component concentrated near the
characteristic set, where $|\sigma_\varepsilon|$ may be small.
\end{corollary}

\begin{proof}
The unconstrained minimizer is
\begin{equation}
\label{eq:auto-057}
  -\frac{\langle F_\varepsilon,A_{\rm base}q\rangle_{\mathcal Y}}
         {\|A_{\rm base}q\|_{\mathcal Y}^2}.
\end{equation}
After projection onto $[0,1]$, Cauchy--Schwarz and
\eqref{eq:off_characteristic_assumption} give
\begin{equation}
\label{eq:auto-058}
  \omega_{\rm lin}
  \le\frac{\|F_\varepsilon\|_{\mathcal Y}}
           {\|A_{\rm base}q\|_{\mathcal Y}}
  \le\frac{C_F}{c_Ac_0}\varepsilon^2.
\end{equation}
The linearized quadratic is
$2\|A_{\rm base}q\|_{\mathcal Y}^2$-strongly convex.  The uniform gap
in~\eqref{eq:off_characteristic_empirical_gap} and the optimization tolerance
therefore give
$|\widehat\omega-\omega_{\rm lin}|
\le\sqrt{2\eta_{\rm sel}+\eta_{\rm opt}}/
\|A_{\rm base}q\|_{\mathcal Y}$, which proves
\eqref{eq:off_characteristic_empirical_bound}.
\end{proof}

The exact nonlinear selector inherits this behavior only when the nonlinear,
sampling, and optimization perturbations are controlled on the
$O(\varepsilon^2)$ neighborhood.  Accordingly, we describe the mechanism as
off-characteristic artifact suppression rather than universal high-frequency
suppression.

In practice, both pretrained networks are frozen and the scalar is optimized
directly on the closed interval $[0,1]$.  We first evaluate
$\widehat{\mathcal J}_{N_s}$ on a fixed grid and then refine the best brackets
with a safeguarded one-dimensional search.  This procedure can attain the exact
boundary choices $\omega=0$ and $\omega=1$; its remaining numerical error is
recorded as $\eta_{\rm opt}$ in Theorem~\ref{thm:blend}.

For reference, the key displayed relations in this subsection are \eqref{eq:auto-054}, \eqref{eq:auto-055}, \eqref{eq:auto-056}, \eqref{eq:auto-057}, \eqref{eq:auto-058}.

\subsection{Gated gradient correlation correction}
\label{sec:gated_gradient}

The modulation removes the carrier from stage I, but the dynamic remainder in
stage II still contains harmonic bands near $m\varepsilon^{-2}$.  We therefore
use a causal temporal curriculum with two diagnostics.
Gradient compatibility measures whether a local update can help neighboring times.
The front residual measures whether the revealed time window is sufficiently trained.
The gate is not included in the
correlation diagnostic itself, thereby avoiding the circular conclusion that a
closed gate implies poor correlation.

\paragraph{Gated gradient flow.}

For either the modulation or dynamic-remainder network, let
\begin{equation}
\label{eq:auto-059}
  \ell_\theta(t)
  =\frac12\|\mathcal R_\theta(\cdot,t)\|_{L^2(\Omega)}^2,
  \qquad
  g_\theta(t)=\nabla_\theta\ell_\theta(t),
\end{equation}
where $\mathcal R_\theta$ is the corresponding PDE residual.  Equivalently,
if $J_{\mathcal R,\theta}(t)$ denotes the residual Jacobian with respect to
the trainable parameters, then
\begin{equation}
\label{eq:auto-060}
  g_\theta(t)
  =J_{\mathcal R,\theta}(t)^\ast
   \mathcal R_\theta(\cdot,t),
  \qquad
  J_{\mathcal R,\theta}(t)
  :=D_\theta\mathcal R_\theta(\cdot,t).
\end{equation}
Thus the diagnostic uses the actual local residual-loss gradient propagated by
backpropagation, not the output Jacobian of $u_\theta$ alone.  The signed,
normalized correlation is
\begin{equation}
\label{eq:grad_corr}
  C_\theta(t,s)
  =\frac{\langle g_\theta(t),g_\theta(s)\rangle}
  {\|g_\theta(t)\|\,\|g_\theta(s)\|+\eta_0},
  \qquad -1\le C_\theta(t,s)\le1,
\end{equation}
with a small $\eta_0>0$.  Unlike an absolute inner product,
\eqref{eq:grad_corr} distinguishes helpful alignment from gradient conflict.
It is also less sensitive to changes in network width and loss scale than an
unnormalized inner product.

A useful diagnostic model follows from the temporal spectrum of the local-loss
gradient.  If $S_g(k;\varepsilon)$ denotes its nonnegative spectral density,
then the normalized autocorrelation has the form
\begin{equation}
\label{eq:spectral_correlation}
  \mathcal C_g(\delta t)
  \approx
  \frac{\int_{\mathbb R}S_g(k;\varepsilon)\cos(k\delta t)\,dk}
       {\int_{\mathbb R}S_g(k;\varepsilon)\,dk}.
\end{equation}
Whenever the fourth moment is finite,
\begin{equation}
\label{eq:small_lag_corr}
  \mathcal C_g(\delta t)
  =1-\frac12M_{2,g}\delta t^2
   +O(M_{4,g}\delta t^4),
  \qquad
  M_{2,g}=\frac{\int k^2S_g(k;\varepsilon)\,dk}
                 {\int S_g(k;\varepsilon)\,dk}.
\end{equation}
More explicitly, with
$M_{4,g}=\int k^4S_g(k;\varepsilon)\,dk/\int S_g(k;\varepsilon)\,dk$,
\begin{equation}
\label{eq:auto-061}
  \left|
  \mathcal C_g(\delta t)-1+\frac12M_{2,g}\delta t^2
  \right|
  \le\frac{M_{4,g}}{24}|\delta t|^4 .
\end{equation}
For the modulation network, the active band remains $|k|=O(1)$, so
$M_{2,g}^{z}=O(1)$ uniformly in $\varepsilon$.  For the dynamic remainder, a
more faithful model than a simple rescaling around the origin is the
harmonic-centered density
\begin{equation}
\label{eq:auto-062}
  S_g^{q}(k;\varepsilon)
  =S_{g,0}(k;\varepsilon)
   +\sum_{m\in\{\pm1,\pm2,\pm3\}}
    a_m(\varepsilon)
    f_m\!\left(k-\frac{m}{\varepsilon^2}\right),
\end{equation}
where each $f_m$ has an $O(1)$ width and finite moments.  If the total weight of
at least one nonzero harmonic is bounded below, then
\begin{equation}
\label{eq:remainder_second_moment}
  M_{2,g}^{q}=\Theta(\varepsilon^{-4}).
\end{equation}
Thus remainder-gradient alignment can change on time separations
$\delta t=O(\varepsilon^2)$, whereas the modulation network does not acquire an
$\varepsilon$-dependent decorrelation scale.  Equations
\eqref{eq:spectral_correlation}--\eqref{eq:remainder_second_moment} are a
spectral diagnostic, not a claim that every optimization trajectory is
stationary.

The default monotone gate is
\begin{equation}
\label{eq:auto-063}
  h_\gamma(t)
  =\frac{1-\tanh\bigl(\alpha(t/T-\gamma)\bigr)}{2},
  \qquad 0\le h_\gamma(t)\le1,
\end{equation}
where $\gamma$ is the normalized temporal front and $\alpha>0$ controls the
transition width.  The gated objective is
\begin{equation}
\label{eq:auto-064}
  J(\theta,\gamma)
  =\int_0^T h_\gamma(t)\ell_\theta(t)\,dt
   +\lambda_{\rm Bd}\mathcal L_{\rm Bd}(\theta).
\end{equation}

Let $I_{\rm front}(\gamma)$ be a short interval centered at the transition of
$h_\gamma$.  We estimate
\begin{align*}
  C_{\rm front}(\gamma)
  &=\mathbb E_{t\in I_{\rm front},\,\delta t\sim\nu_\varepsilon}
    [C_\theta(t,t+\delta t)],
  \\
  \overline R_{\rm front}(\gamma)
  &=\frac{\mathbb E_{t\in I_{\rm front}}[\ell_\theta(t)]}
  {R_{\rm ref}+\eta_0},
\end{align*}
where $R_{\rm ref}$ is a fixed reference or an exponential moving average.
The perturbation law $\nu_\varepsilon$ uses $O(1)$ radii for stage I and
$O(\varepsilon^2)$ radii for stage II, in agreement with
\eqref{eq:remainder_second_moment}.

Define the readiness score
\begin{equation}
\label{eq:readiness_score}
  S_i
  =\operatorname{sigmoid}
    \bigl(a(C_{\rm front}(\gamma^i)-C_0)\bigr)
    \exp\bigl(-b\overline R_{\rm front}(\gamma^i)\bigr),
\end{equation}
where $C_0\in(-1,1)$ and $a,b>0$.  Rather than assigning a heuristic speed, we
choose the next gate velocity as the solution of a one-step regularized
progress problem,
\begin{equation}
\label{eq:gate_speed_problem}
  v_i^\ast
  =\arg\max_{0\le v\le\Delta_{\max}}
    \left\{S_i v-\frac{\lambda_\gamma}{2}v^2\right\}
  =\min\left\{\Delta_{\max},\frac{S_i}{\lambda_\gamma}\right\}.
\end{equation}
When a strict nonincrease certificate for a local merit upper bound is needed,
we use the conservative variant
\begin{equation}
\label{eq:auto-065}
  v_{i,{\rm safe}}^\ast
  =
  \min\left\{
  \Delta_{\max},
  \frac{(S_i-B_i)_+}{L_\gamma}
  \right\},
  \qquad
  B_i=|\partial_\gamma J(\theta^i,\gamma^i)|,
\end{equation}
where $L_\gamma$ is a local Lipschitz constant of $\partial_\gamma J$.
Indeed, the following upper model is used:
\begin{equation}
\label{eq:auto-066}
J(\theta^i,\gamma^i+v)\le J(\theta^i,\gamma^i)+B_iv+L_\gamma v^2/2 .
\end{equation}
It is dominated by the progress credit $S_iv$ when
$0\le v\le (S_i-B_i)_+/L_\gamma$.
This safe variant is more conservative than
\eqref{eq:gate_speed_problem} and is used as a certificate, not as an
assumption of global optimization.
The front is then updated by
\begin{equation}
\label{eq:gamma_update}
  \gamma^{i+1}
  =\Pi_{[0,\gamma_{\max}]}
    \bigl(\gamma^i+\eta_\gamma v_i^\ast\bigr).
\end{equation}
This is a local optimality statement for the explicit surrogate
in~\eqref{eq:gate_speed_problem}.  It is not a global optimality claim for
the nonconvex neural training trajectory.

\begin{proposition}[Causal monotonicity and readiness response]
\label{prop:gate_schedule}
Let $\lambda_\gamma>0$, $\eta_\gamma>0$, and
$0\le\gamma^0\le\gamma_{\max}$.  Then:
\begin{enumerate}
  \item $v_i^\ast$ is the unique solution of
  \eqref{eq:gate_speed_problem};
  \item $\gamma^{i+1}\ge\gamma^i$ and
  $\gamma^i\le\gamma_{\max}$ for all $i$;
  \item before clipping, the gate speed is increasing in
  $C_{\rm front}$ and decreasing in $\overline R_{\rm front}$.
\end{enumerate}
Consequently, aligned local gradients and a small front residual accelerate the gate.
Gradient conflict or a poorly trained front slows it down.
\end{proposition}

\begin{proof}
The objective in~\eqref{eq:gate_speed_problem} is strictly concave in $v$, so
its projected stationary point is the unique maximizer shown in the equation.
Since $S_i\ge0$ and $v_i^\ast\ge0$, each update is nonnegative before clipping.
Projection onto $[0,\gamma_{\max}]$ makes the sequence bounded.
Finally, the sigmoid factor in~\eqref{eq:readiness_score} is increasing in
$C_{\rm front}$ and the exponential factor is decreasing in
$\overline R_{\rm front}$.
\end{proof}

\begin{corollary}[Correlation-residual response]
\label{cor:gate_response}
If
$C_{\rm front}(\gamma^i)\ge C_0+\kappa_c$ and
$\overline R_{\rm front}(\gamma^i)\le r_0$, then the unclipped progress rule
\eqref{eq:gate_speed_problem} satisfies
\begin{equation}
\label{eq:auto-067}
  v_i^\ast
  \ge
  \min\left\{
    \Delta_{\max},
    \frac{\operatorname{sigmoid}(a\kappa_c)\exp(-br_0)}
         {\lambda_\gamma}
  \right\}.
\end{equation}
Thus a positive signed-correlation margin and a controlled front residual give
a direct lower bound on the causal-front advance.
\end{corollary}

\begin{corollary}[Finite gate completion under persistent readiness]
\label{cor:finite_gate_completion}
Assume that before reaching $\gamma_{\max}$ the readiness score satisfies
$S_i\ge S_{\min}>0$.  Then the update~\eqref{eq:gamma_update} reaches
$\gamma_{\max}$ after at most
\begin{equation}
\label{eq:auto-068}
  \left\lceil
  \frac{\gamma_{\max}-\gamma^0}
       {\eta_\gamma\min\{\Delta_{\max},S_{\min}/\lambda_\gamma\}}
  \right\rceil
\end{equation}
outer gate updates.
\end{corollary}

A learnable gate may be used, but it must preserve causality.  One admissible
choice is a convex mixture of monotone sigmoids,
\begin{equation}
\label{eq:auto-069}
  h_{\psi,\gamma}(t)
  =\sum_{j=1}^{M}\pi_j
    \operatorname{sigmoid}
    \bigl(-\alpha_j(t/T-\gamma-\beta_j)\bigr),
  \quad
  \pi_j\ge0,\quad\sum_j\pi_j=1,\quad\alpha_j>0.
\end{equation}
Every component is decreasing in $t$ and increasing in $\gamma$, so the mixture
cannot reveal a future interval while hiding an earlier one.  Unconstrained
non-monotone Gaussian mixtures are not used as causal gates.

\begin{figure}[!htb]
  \centering
  \includegraphics[width=0.95\textwidth]{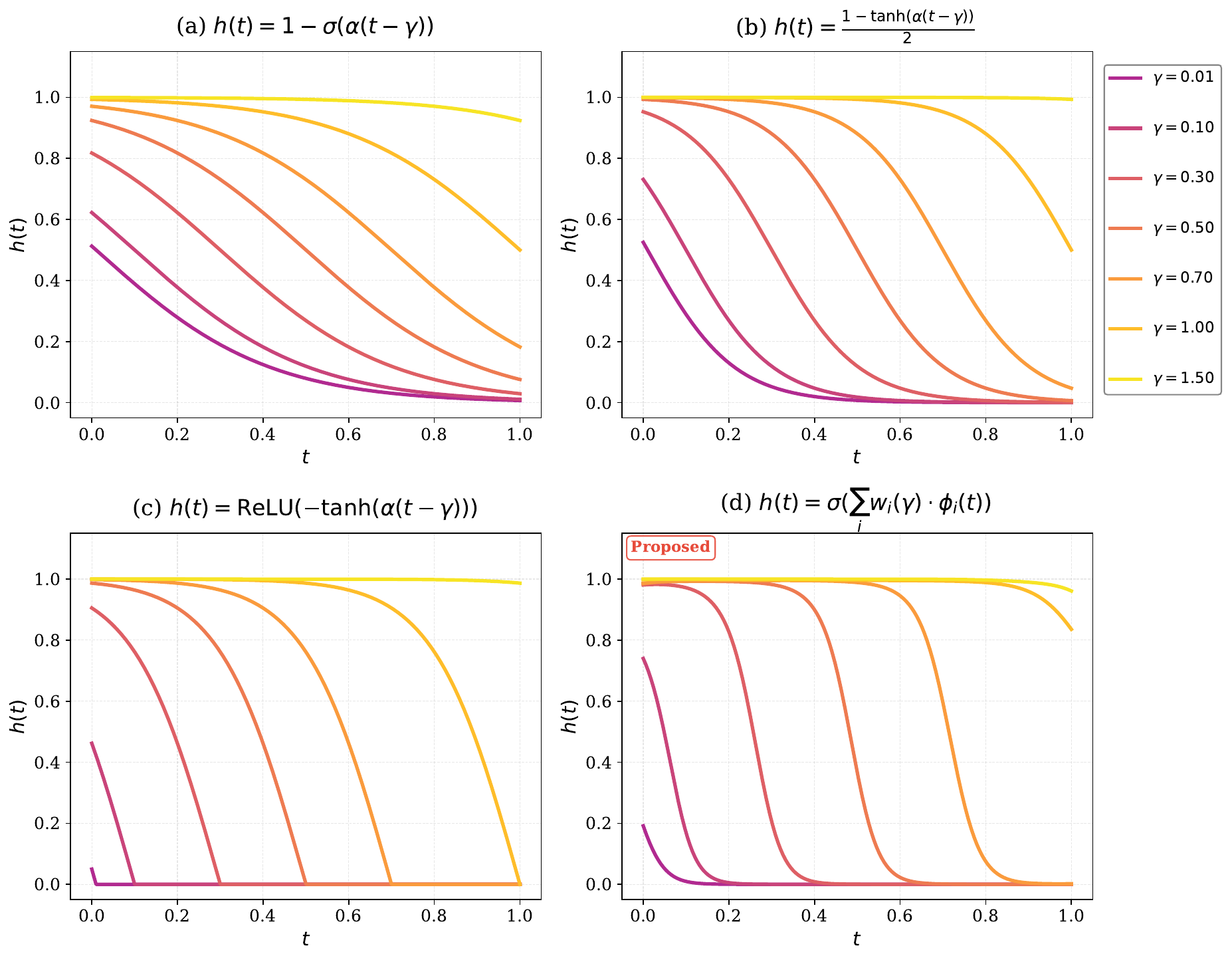}
  \vspace{-8pt}
  \caption{Comparison of different gate functions $h(t)$ under varying values
  of the shift parameter $\gamma$.  The hand-crafted sigmoid, tanh, and
  ReLU-tanh gates provide monotone causal curricula, while the learnable
  monotone mixture offers additional shape flexibility without revealing future
  times before earlier ones.}
  \label{fig:gate_func}
\end{figure}

For reference, the key displayed relations in this subsection are \eqref{eq:auto-059}, \eqref{eq:auto-060}, \eqref{eq:auto-061}, \eqref{eq:auto-062}, \eqref{eq:auto-063}, \eqref{eq:auto-064}.
The continuation of this numbered set is \eqref{eq:auto-065}, \eqref{eq:auto-066}, \eqref{eq:auto-067}, \eqref{eq:auto-068}, \eqref{eq:auto-069}.

\paragraph{Random time perturbation.}

For a base point $(\bm{x},t)$, choose $L_s$ temporal scales and sample
$\delta t_l^i\in[-R_l,R_l]$, $i=1,\ldots,k_l$.  We use
\begin{equation}
\label{eq:perturbation_radii}
  R_l^{z}=c_lR_0,
  \qquad
  R_l^{q}=c_l\varepsilon^2R_0,
\end{equation}
with fixed dimensionless $c_l$.  To avoid assuming a uniform correlation
radius over the entire time interval, these nominal radii are capped by the
empirical local radius
\begin{equation}
\label{eq:empirical_corr_radius}
  R_{\rm corr}(t)
  =\sup\left\{r\ge0:
    \min_{|\delta|\le r}
    C_\theta(t,t+\delta)\ge C_{\min}\right\},
  \qquad
  R_l(t)\leftarrow\min\{R_l,R_{\rm corr}(t)\},
\end{equation}
where the minimum is estimated from a small candidate set and treated as a
stop-gradient sampling statistic.  A deterministic local lower bound is
\begin{equation}
\label{eq:auto-070}
  R_{\rm corr}(t)
  \ \gtrsim\
  \frac{m_g(t)}{L_g(t)}
  \sqrt{\frac{1-C_{\min}}{2}},
  \qquad
  m_g(t)\le\|g_\theta(t)\|,\quad
  \|g_\theta(t+\delta)-g_\theta(t)\|\le L_g(t)|\delta|,
\end{equation}
whenever the displayed lower norm and local Lipschitz bounds hold and
$\eta_0$ is negligible compared with $m_g(t)^2$.  This deterministic estimate
is used only to motivate the empirical cap; the implementation still estimates
$R_{\rm corr}(t)$ directly from signed correlations.

The augmented set and its coordinate representations are
\begin{equation}
\label{eq:auto-071}
\begin{aligned}
  \operatorname{Diff\text{-}Aug}(\bm{x},t)
  &=\left\{(\bm{x},t),
  \{(\bm{x},t+\delta t_l^i)
  :i=1,\ldots,k_l\}_{l=1}^{L_s}\right\},\\
  \bm{x}_{\rm point}&=\mathcal P(\bm{x},t),\\
  \bm{x}_{\rm region}^{i,l}
  &=\mathcal P(\bm{x},t+\delta t_l^i),
\end{aligned}
\end{equation}
where $\mathcal P:\mathbb R^{d+1}\to\mathbb R^{d_{\rm model}}$ is a
coordinate encoder.

For the optimization analysis, define nonnegative normalized weights
\begin{equation}
\label{eq:auto-072}
  w_l^i(t)
  =\frac{[h_\gamma(t+\delta t_l^i)+\eta_0]
           [K_l(\delta t_l^i)+\eta_0]}
  {\sum_{j=1}^{k_l}[h_\gamma(t+\delta t_l^j)+\eta_0]
           [K_l(\delta t_l^j)+\eta_0]},
  \qquad \sum_{i=1}^{k_l}w_l^i(t)=1,
\end{equation}
where $K_l$ is a nonnegative local kernel.  The augmented local loss is
\begin{equation}
\label{eq:auto-073}
  \overline\ell_{\theta,l}(t)
  =\sum_{i=1}^{k_l}w_l^i(t)
    \ell_\theta(t+\delta t_l^i).
\end{equation}
The weights are independent of $\theta$ within an inner update, so its gradient
is the vector
\begin{equation}
\label{eq:correct_augmented_gradient}
  \nabla_\theta\overline\ell_{\theta,l}(t)
  =\sum_{i=1}^{k_l}w_l^i(t)
    g_\theta(t+\delta t_l^i).
\end{equation}
Equation~\eqref{eq:correct_augmented_gradient}, rather than an inner product of
Jacobians, is the gradient actually propagated by backpropagation.
It also gives the standard variance reduction mechanism.  Conditional on the
base time and fixed weights, if the active neighbor gradients have covariance
operator trace at most $\sigma_g^2$, then
\begin{equation}
\label{eq:auto-074}
  \mathbb E\left\|
  \sum_iw_l^i(t)g_\theta(t+\delta t_l^i)
  -\mathbb E[g_\theta(t+\delta t)]
  \right\|^2
  \le
  \sigma_g^2\sum_i\bigl(w_l^i(t)\bigr)^2 .
\end{equation}
Thus broad, positively correlated neighborhoods reduce stochastic gradient
variance, whereas the cap~\eqref{eq:empirical_corr_radius} prevents averaging
across sign-conflicting temporal phases.

\begin{proposition}[One-step cross-time descent]
\label{prop:cross_time_descent}
Assume each perturbed loss $\ell_i(\theta)
=\ell_\theta(t+\delta t_l^i)$ has an $L_\theta$-Lipschitz gradient.  Let
$\overline g=\sum_iw_l^ig_i$ and update
$\theta^+=\theta-\eta\overline g$.  Then
\begin{equation}
\label{eq:neighbor_descent_bound}
  \ell_j(\theta^+)
  \le\ell_j(\theta)
  -\eta\langle g_j,\overline g\rangle
  +\frac{L_\theta\eta^2}{2}\|\overline g\|^2.
\end{equation}
If $\langle g_j,\overline g\rangle>0$ for all active neighbors and
\begin{equation}
\label{eq:descent_stepsize}
  0<\eta<
  \frac{2\min_j\langle g_j,\overline g\rangle}
       {L_\theta\|\overline g\|^2},
\end{equation}
then every active neighboring loss decreases in one step.  In particular,
positive signed correlations in~\eqref{eq:grad_corr} provide a directly
verifiable sufficient condition for temporal coherence.
\end{proposition}

\begin{proof}
Apply the descent lemma to each $\ell_j$ with displacement
$-\eta\overline g$ to obtain~\eqref{eq:neighbor_descent_bound}.  Condition
\eqref{eq:descent_stepsize} makes the right-hand decrease term strictly
negative for every $j$.
\end{proof}

For reference, the key displayed relations in this subsection are \eqref{eq:auto-070}, \eqref{eq:auto-071}, \eqref{eq:auto-072}, \eqref{eq:auto-073}, \eqref{eq:auto-074}.

\paragraph{Multiscale time mixing.}

At each scale, the perturbed representations are pooled with the same causal
weights,
\begin{equation}
\label{eqn:multiscale}
\begin{aligned}
  \bm{x}_{\rm region}^{l}
  &=\operatorname{Pooling}
    \left(\{w_l^i(t)\bm{x}_{\rm region}^{i,l}\}_{i=1}^{k_l}\right),
    \qquad l=1,\ldots,L_s,\\
  \bm{x}_{\rm mix}
  &=\mathcal M\left(
    \bm{x}_{\rm point},\bm{x}_{\rm region}^{1},\ldots,
    \bm{x}_{\rm region}^{L_s}\right),
\end{aligned}
\end{equation}
where $\mathcal M$ is a lightweight mixing network.  A decoder
$\mathcal D$ maps $\bm{x}_{\rm mix}$ to either
$\widetilde z_{\theta_1}$ or $\widetilde q_{\theta_2}$.  Stage I uses
$O(1)$ temporal neighborhoods to learn the modulation.
Stage II uses $O(\varepsilon^2)$ neighborhoods to avoid averaging unrelated
remainder phases.  This architecture enlarges the temporal receptive field while
preserving the actual vector-valued gradient in
\eqref{eq:correct_augmented_gradient}.

\subsection{Justification of gradient correlation}
\label{sec:grad_corr_justify}

We keep the original gradient-correlation justification explicit.
It explains why temporal perturbation is useful beyond data augmentation.
For a base time $t$, consider an active perturbed neighborhood
$\mathcal T_l(t)=\{t+\delta t_l^i\}_{i=1}^{k_l}$.
Define $g_i=g_\theta(t+\delta t_l^i)$.
Also define the weighted augmented gradient
$\overline g_l(t)=\sum_i w_l^i(t)g_i$, as in
\eqref{eq:correct_augmented_gradient}.  The signed correlation
\eqref{eq:grad_corr} implies
\begin{equation}
\label{eq:positive_neighbor_alignment}
  \langle g_j,\overline g_l(t)\rangle
  =\sum_i w_l^i(t)\langle g_j,g_i\rangle>0
  \qquad \text{for every active neighbor } j
\end{equation}
whenever the active neighborhood lies inside a local positive-correlation radius.
Substituting~\eqref{eq:positive_neighbor_alignment} into
Proposition~\ref{prop:cross_time_descent} gives a direct descent condition.
With a small enough step size, all active neighboring losses decrease.
Thus multiscale perturbation improves temporal propagation by changing the update direction.
It does not merely add more collocation points.

The scale choice is also justified by the spectral estimates
\eqref{eq:small_lag_corr}--\eqref{eq:remainder_second_moment}.  The modulation
network has an $O(1)$ temporal correlation length, so stage I can use
$O(1)$ neighborhoods.  The dynamic remainder has bands centered near
$m/\varepsilon^2$.  Its correlation length can shrink to $O(\varepsilon^2)$.
Stage II therefore uses the smaller radii in
\eqref{eq:perturbation_radii}.  This is the reason NeuralMD can preserve the
smooth amplitude flow while stabilizing the high-frequency remainder flow.

\subsection{Training implementation}
\label{sec:train}

The training process is summarized in Algorithm~\ref{alg:neuralmd}.  Compared
with the compact protocol above, the algorithm records three implementation stages.
These are causal stage-I training, causal stage-II training, and scalar remainder selection.

\begin{algorithm}[!htb]
\caption{NeuralMD training}
\label{alg:neuralmd}
\begin{algorithmic}[1]
\REQUIRE Collocation sets $P_f,P_b$, independent selector set
$\widehat P_f$, parameters $\varepsilon,\lambda$, temporal scales
$\{R_l,k_l\}_{l=1}^{L_s}$, gate parameters
$\gamma^0,\eta_\gamma,\Delta_{\max}$.
\ENSURE Trained modulation network $z_{\theta_1^\ast}$, dynamic remainder
network $q_{\theta_2^\ast}$, and selected remainder weight $\widehat\omega$.
\STATE Construct the hard initial-data lift
$z_{\theta_1}$ in~\eqref{eq:z_hard_lift}.
\FOR{stage-I iterations}
  \STATE Sample $(\bm{x},t)\in P_f$ and build
  $\operatorname{Diff\text{-}Aug}(\bm{x},t)$ with stage-I radii
  $R_l^z=c_lR_0$.
  \STATE Compute the NLSW residual loss~\eqref{eq:Lz} using the multiscale
  representation~\eqref{eqn:multiscale}.
  \STATE Estimate the front correlation $C_{\rm front}$ and normalized front
  residual $\overline R_{\rm front}$.
  \STATE Update the causal front by~\eqref{eq:gamma_update}.
  \STATE Update $\theta_1$ by Adam/L-BFGS using the gated objective.
\ENDFOR
\STATE Freeze $\theta_1^\ast$ and construct
$r_{\theta_2}=r_{\rm ic}+q_{\theta_2}$ in~\eqref{eq:r_split}.
\FOR{stage-II iterations}
  \STATE Sample $(\bm{x},t)\in P_f$ and build
  $\operatorname{Diff\text{-}Aug}(\bm{x},t)$ with stage-II radii
  $R_l^q=c_l\varepsilon^2R_0$.
  \STATE Compute the dynamic-remainder residual loss~\eqref{eq:Lr}.
  \STATE Estimate the front diagnostics and update $\gamma$ by
  \eqref{eq:gamma_update}.
  \STATE Update $\theta_2$ by Adam/L-BFGS using the gated objective.
\ENDFOR
\STATE Freeze both pretrained networks and form
$u_\omega=u_{\rm base}+\omega q_{\theta_2^\ast}$.
\STATE Select
$\widehat\omega\in[0,1]$ by minimizing
\eqref{eq:L_blend_empirical}, or the cost-aware objective
\eqref{eq:cost_aware_selector}, on $\widehat P_f$ with a grid search followed
by safeguarded one-dimensional refinement.
\STATE Report $\widehat\omega$, initial/boundary errors, the selector objective
curve, front correlations, and the validation residual.
\RETURN $z_{\theta_1^\ast}$, $q_{\theta_2^\ast}$, and $\widehat\omega$.
\end{algorithmic}
\end{algorithm}

\subsection{Interpretability of NeuralMD}

NeuralMD uses an MLP-based architecture grounded in the Universal Approximation Theorem (UAT) \citep{hornik1989mlp}.
The UAT states that sufficiently wide neural networks approximate continuous functions $f:\mathbb R^n\to\mathbb R$.
For any $\epsilon > 0$, there exists a suitable hidden layer width $N(\epsilon)$ such that
\begin{equation}
\label{eq:auto-075}
f(x) \approx \sum_{i=1}^{N(\epsilon)} a_i \sigma(w_i x + b_i),
\end{equation}
where $\sigma(\cdot)$ is the activation function (we use tanh, which helps stabilize training and capture the asymptotic behavior of the solution). 
For deep networks, the structure is expressed as a composition of linear mappings and nonlinear activations, as
\begin{equation}
\label{eq:auto-076}
\text{MLP}(x) = (W_L \circ \sigma_{L-1} \circ W_{L-1} \circ \cdots \circ \sigma_1 \circ W_1)(x).
\end{equation}
However, MLPs concentrate nonlinearity at nodes.
Their weight matrices act as high-dimensional tensors, producing highly coupled features.
These features are difficult to interpret through input-output relationships.
In NeuralMD, this "black-box" mapping obscures the physical meaning of oscillation removal.

To address this, we introduce Kolmogorov--Arnold Networks (KANs) \citep{liu2024kan} as an interpretable alternative. 
KANs are based on the Kolmogorov--Arnold representation theorem (KAT) \citep{schmidt2021kg}.
The theorem expresses any continuous function $f(x_1,\ldots,x_n)$ as finite sums of one-dimensional functions:
\begin{equation}
\label{eq:auto-077}
f(x_1, \ldots, x_n) = \sum_{q=1}^{2n+1} \Phi_q \left( \sum_{p=1}^{n} \varphi_{q,p}(x_p) \right),
\end{equation}
where $\Phi_q$ and $\varphi_{q,p}$ are learnable continuous one-dimensional functions.
Unlike MLPs, KANs apply nonlinearity on edges rather than nodes.
Each edge corresponds to a learnable one-dimensional function $\varphi_{l,j,i}(\cdot)$.
Nodes then perform linear summation:
\begin{equation}
\label{eq:auto-078}
x_{l+1,j} = \sum_{i=1}^{n_l} \varphi_{l,j,i}(x_{l,i}), \quad x_{l+1} = \Phi_l x_l.
\end{equation}
Thus, the entire network can be viewed as a composition of functions:
\begin{equation}
\label{eq:auto-079}
\text{KANs}(x) = \left( \Phi_{L-1} \circ \Phi_{L-2} \circ \cdots \circ \Phi_0 \right)(x),
\end{equation}
allowing each edge to correspond to an explicit one-dimensional function, thereby improving interpretability.

To visualize and symbolize edge functions, KANs use B-splines for parameterization:
\begin{equation}
\label{eq:auto-080}
\varphi(x) = w_b b(x) + w_s \, \text{spline}(x), \quad \text{spline}(x) = \sum_{i} c_i B_i(x),
\end{equation}
where $b(x) = \text{silu}(x)$.
The control points $c_i$ shape the edge-function curve.
They can explicitly learn operator terms in the equations.
During training, sparse and entropy regularization terms are introduced:
\begin{equation*}
\mathcal{L}^z(\theta_1) = \lambda_{\mathrm{Res}} \mathcal{L}_{\mathrm{Res}}^z(\theta_1) + \lambda_{\mathrm{Ic}} \mathcal{L}_{\mathrm{Ic}}^z(\theta_1) + \lambda \left( \mu_1 \sum_{l} \| \Phi_l \|_1 + \mu_2 \sum_{l} S(\Phi_l) \right),
\end{equation*}
where $S(\Phi_l)$ is the entropy regularization term.
It measures the information distribution in edge functions of the $l$-th KAN layer:
\begin{equation*}
S(\Phi_l) = -\sum_{i,j} \frac{\left|\varphi_{l,i,j}\right|}{\sum_{i',j'} \left|\varphi_{l,i',j'}\right|} \log\left( \frac{\left|\varphi_{l,i,j}\right|}{\sum_{i',j'} \left|\varphi_{l,i',j'}\right|} + \delta \right).
\end{equation*}
When each learnable basis function $\Phi_{l,i,j} \in C^{k+1}$, a B-spline approximation exists.
For order $k$, the approximation error is bounded by
\begin{equation*}
\| f - \left( \Phi_{L-1}^G \circ \cdots \circ \Phi_0^G \right) x \|_{C^m} \le C G^{-(k+1-m)}, \quad 0 \le m \le k.
\end{equation*}
where $G$ is the spline density, and $k=3$ for cubic splines.

We integrate KANs into NeuralMD while retaining operator-space approximation capacity.
This provides more structured and symbolic physical interpretability.

For reference, the key displayed relations in this subsection are \eqref{eq:auto-075}, \eqref{eq:auto-076}, \eqref{eq:auto-077}, \eqref{eq:auto-078}, \eqref{eq:auto-079}, \eqref{eq:auto-080}.

\section{Numerical experiments}
\label{sec:experiment}
In this section, we demonstrate the numerical performance of the proposed NeuralMD method. 
We focus on spectral bias and propagation failure in the temporally oscillatory NKGE.
The goal is not to tune the best possible baseline model.
We report the main benchmark results and diagnostic comparisons with other collocation-based methods. 
All experiments are conducted on an Nvidia A100-SXM4-80GB GPU.
Code, data-generation scripts, and raw logs will be released with the camera-ready version.
The reported experiments should be interpreted with the implementation details below.

\subsection{Benchmarks}
For all benchmarks, we use three random time regions, denoted ${\text{\#scale}}=3$.
Their sizes are $\{R_1,R_2,R_3\}=\{10^{-2},5\times10^{-2},9\times10^{-2}\}$. 
These regions are designed to represent different temporal scales that are characteristic of the underlying problem. 
In each region, we perturb the time coordinate with different perturbation counts.
We set $\{k_1,k_2,k_3\}=\{3,5,7\}$ for multiscale configurations. 
The representation feature dimension is set to $d_{\text{model}} = 64$, which provides a balanced trade-off between model capacity and computational efficiency.
For 1D, 2D, and 3D NKGE tests, we follow the protocol of \citet{zhao2023former}.
Training combines Adam with L-BFGS \citep{liu1989limited,kingma2014adam}.
Adam is used for the first 500 iterations to obtain fast initial convergence.
L-BFGS is then used for 500 iterations to refine the solution. 
All models are trained for 1,000 iterations in PyTorch \citep{paszke2019torch}.
To assess performance, we evaluate relative L1 error (rMAE) and relative root mean square error (rRMSE). 
These metrics provide a robust measure of the model's accuracy by comparing the predicted solutions against the ground truth. 
The rMAE quantifies average absolute error relative to the true solution.
The rRMSE measures relative root-mean-square error.
Both metrics are computed over all evaluation points.
This enables comparison across the tested regimes.

\subsection{Baselines}
In addition to vanilla PINNs \citep{raissi2019phy}, we compare NeuralMD with ten PINN variants.
These include QRes, FLS, CausalPINNs, PirateNet, RoPINNs, ProPINNs, PINNsFormer, SetPINNs, PINNsMamba, MSPINNs, and PhasePINNs.
They follow conventional PINN training, where different collocation points are optimized independently. 
PINNsFormer \citep{zhao2023former}, SetPINNs \citep{nagda2024set}, and PINNsMamba \citep{xu2025mamba} are based on the Transformer backbone to capture spatiotemporal correlation among PDEs.
PirateNet and PINNsFormer are previous state-of-the-art models. 
We also integrate R3 sampling \citep{wang2023r3} and loss reweighting \citep{wang2022why}.
We include RoPINN \citep{wu2024ropinn} to test orthogonality with recent optimization methods.

\subsection{Model Configuration}
In our experiments, we compare NeuralMD with twelve baselines. 
We use the following implementation details for these baselines:
\begin{itemize}
\item For vanilla PINN \citep{raissi2019phy}, QRes \citep{bu2021quad}, FLS \citep{wong2022learn}, and PirateNet \citep{wang2024pirat}, we follow the PyTorch implementation of these models provided in the official version. 
Specifically, vanilla PINN, QRes, and FLS all use 9 layers with 64 hidden channels for the feedforward layer.
\item As for SetPINNs \citep{nagda2024set}, PINNsFormer \citep{zhao2023former}, and PINNsMamba \citep{xu2025mamba}, we implement with only 1 encoder layer, which contains 64 hidden channels for the attention mechanism and 128 hidden channels for the feedforward layer.
\item For RoPINNs~\citep{wu2024ropinn} and ProPINNs~\citep{wu2025propinn}, we use their official code to obtain comparative results. 
For RoPINNs, the initial region size is $10^{-4}$.
We use 10 past iterations and sample one point per region per iteration.  
Since these points are pre-selected, their values are not updated after sampling.
Thus we do not apply region sampling to them in our experiments.  
For ProPINNs, we perturb the input dimension $(d+1)$ across 3 scales with sizes 0.03, 0.05, and 0.07.  
We use uniform sampling with fixed nodes to sample the expanded regions, and the hidden dimension is set to 64.
\item For \textbf{NeuralMD}, unless otherwise stated, we randomly perturb three different time scales to construct multiscale time regions with sizes $1\times10^{-2}$, $5\times10^{-2}$, and $9\times10^{-2}$, and sample 3, 5, and 7 points, respectively.  
For time-regime perturbation sampling, spatial points remain unchanged.
Time points are sampled using a "tanh" causal gate with $\gamma\in[-0.5,1.5]$.  
The time-region mixing layer consists of two linear layers with an activation in between, applied only to the time region dimension. 
This layer first projects the three time scales to 8 and then to 1 after the activation layer.  
The final projection layer consists of three linear layers with inner activations, and the hidden dimension is set to 64.
\end{itemize}

\subsection{Metrics}
The model's performance is evaluated using two common metrics: the relative L1 error (rMAE) and the relative Root Mean Square Error (rRMSE). For the ground truth $u$ and the predicted solution $u_\theta$, these metrics are defined as
\begin{equation*}
    \text{rMAE: }\frac{\sum_{i=1}^n \left| u_{\theta}(\bm{x}_i) - u(\bm{x}_i) \right|}{\sum_{i=1}^n \left| u(\bm{x}_i) \right|}, \quad
    \text{rRMSE: }\sqrt{\frac{\sum_{i=1}^n \left| u_{\theta}(\bm{x}_i) - u(\bm{x}_i) \right|^2}{\sum_{i=1}^n \left| u(\bm{x}_i) \right|^2}}.
\end{equation*}
where $\{\bm{x}_i\}_{i=1}^n$ denotes the set of evaluation points, which is kept separate from the training collocation set.

\subsection{1D nonlinear Klein-Gordon equation in whole regime}
We evaluate NeuralMD against 12 baselines across representative regimes: the relativistic regime ($\varepsilon = 0.8$), the transition regime ($\varepsilon = 0.5$), and the nonrelativistic limit regime ($\varepsilon = 0.1, 0.01$).
We take $d=1$ and $\lambda=1$ in the NKGE (\ref{eqn:KG}).
Take the initial data as
\begin{equation*}
 \phi_1(x)=\frac{e^{-x^2}}{\sqrt{\pi}},\quad \phi_2(x)=\frac{1}{2}\mathrm{sech}(x^2)\sin(x), \qquad x\in{\mathbb R},
\end{equation*}
\begin{equation*}
\phi_1(x)= \frac{3\sin(x)}{e^{x^2/2}+e^{-x^2/2}}, \quad \phi_2(x)= \frac{2e^{-x^2}}{\sqrt{\pi}},\quad x\in {\mathbb R}.
\end{equation*}

As shown in Table~\ref{tab:compare_results}, for $\varepsilon = 0.8$ (low-frequency temporal oscillation), several models, including PINNs, QRes, PirateNet, MSPINNs, PhasePINNs, and CausalPINNs, have substantially larger errors than the best-performing methods. 
Here time oscillation is not the primary cause of failure.
The lack of temporal causality is the dominant factor. 
This suggests that many models are affected by propagation failure even when the temporal oscillation is not severe.  
MSPINNs and PhasePINNs stretch high-frequency oscillation into the low-frequency domain.
They still fail to handle time oscillation effectively, even in low-frequency regions. 
These models do not establish temporal correlation, leading to early propagation failure.
The CausalPINNs method, which attempts to establish temporal causality, also has relatively large errors under the configuration used here. 
In contrast, models such as FLS, which use Fourier feature scaling without time correlation, do not experience failure for low-frequency oscillation problems. 
However, they have larger errors when handling intermediate-frequency temporal oscillation at $\varepsilon = 0.5$.
Several propagation-failure models perform competitively in the low-frequency region.
These include PINNsFormer, SetPINNs, RoPINNs, and PINNsMamba. 
However, they suffer accuracy degradation when applied to intermediate-frequency temporal oscillation at $\varepsilon = 0.5$.
The ProPINNs model, with gradient correction across spatiotemporal regions, maintains competitive accuracy in the intermediate-frequency regime. 
However, as $\varepsilon$ decreases to 0.1, ProPINNs experiences severe failure, and all models fail with no predictive accuracy.

In contrast, NeuralMD compensates the remainder amplitude in the relativistic regime.
It uses random time perturbations and multiscale mixing to alleviate propagation failure. 
In the nonrelativistic limit regime, NeuralMD separates the oscillatory carrier.
It then uses WKB reconstruction to recover the high-frequency temporal structure. 
In the transition regime, NeuralMD balances phase separation with remainder-amplitude compensation, reducing errors associated with both spectral bias and propagation failure. 
NeuralMD achieves the lowest errors among the tested methods in these experiments, especially in the nonrelativistic limit regime ($\varepsilon = 0.1, 0.01$).

\begin{figure}[!htb]
    \centering
    \includegraphics[width=0.95\textwidth]{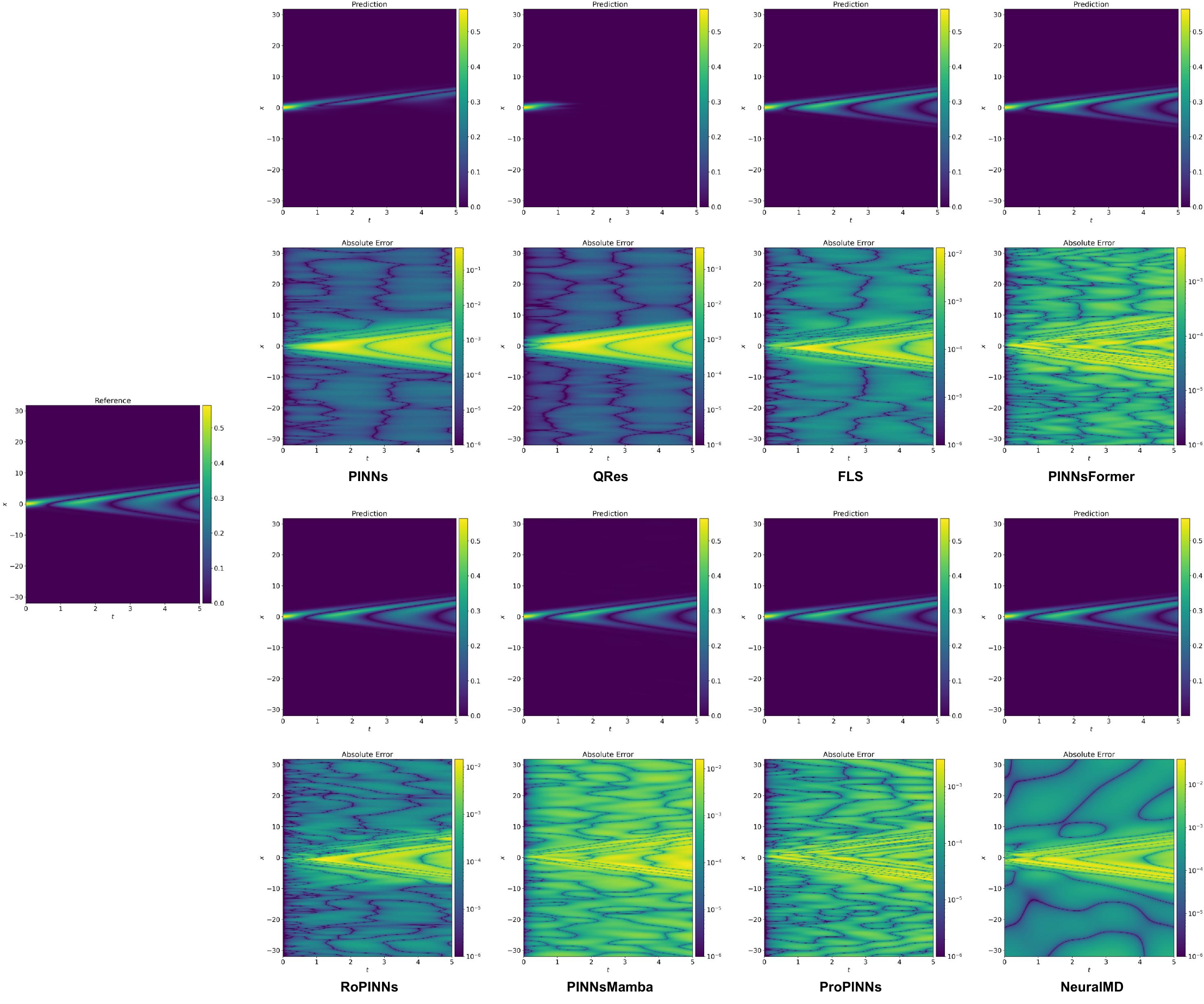}
    \vspace{-8pt}
    \caption{The prediction solution of NeuralMD and baselines for $\varepsilon=0.8,T=5.0$.}
    \label{fig:1d_eps_0.8}
\end{figure}

\begin{figure}[!htb]
    \centering
    \includegraphics[width=0.95\textwidth]{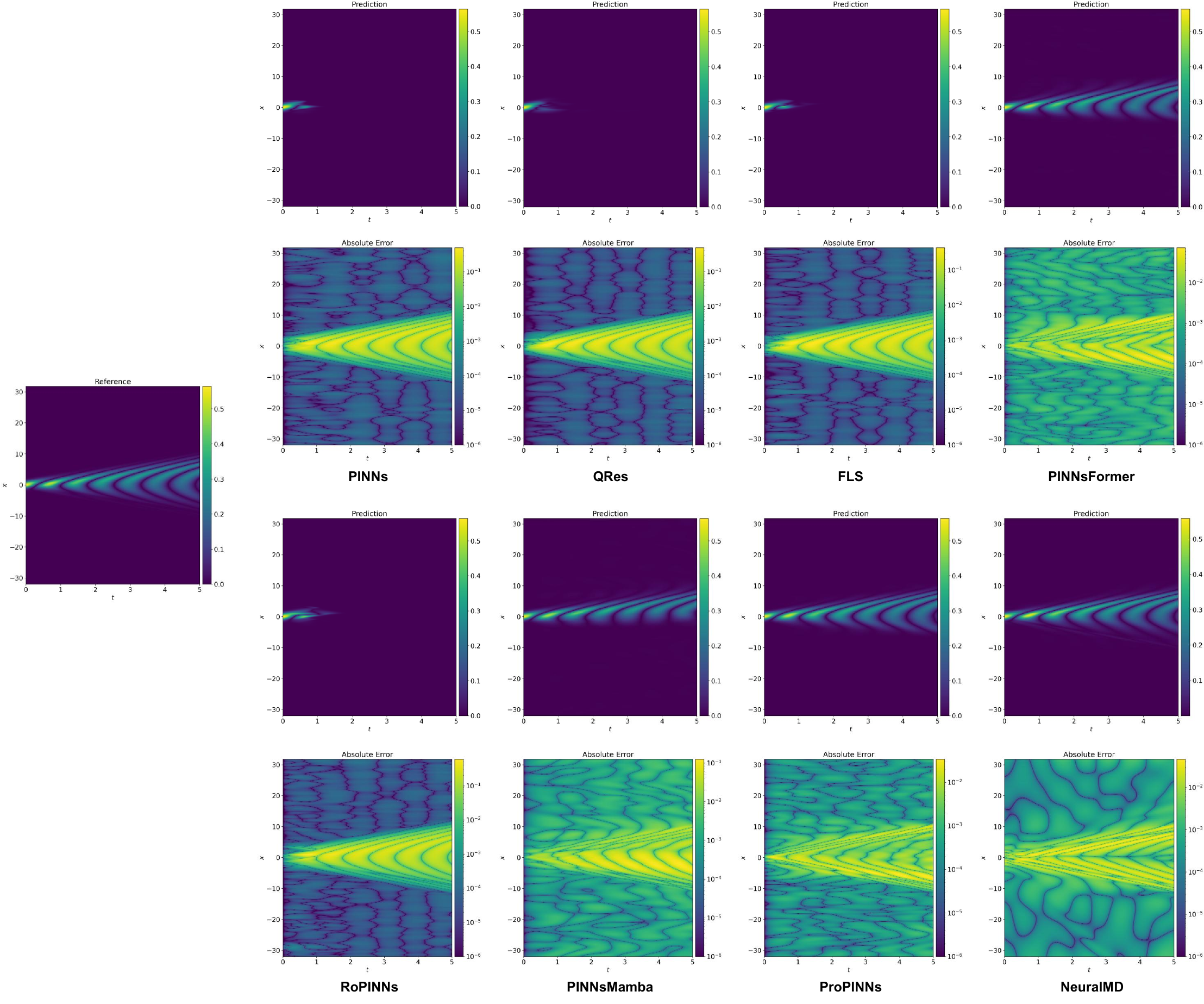}
    \vspace{-8pt}
    \caption{The prediction solution of NeuralMD and baselines for $\varepsilon=0.5,T=5.0$.}
    \label{fig:1d_eps_0.5}
\end{figure}

\begin{figure}[!htb]
    \centering
    \includegraphics[width=0.95\textwidth]{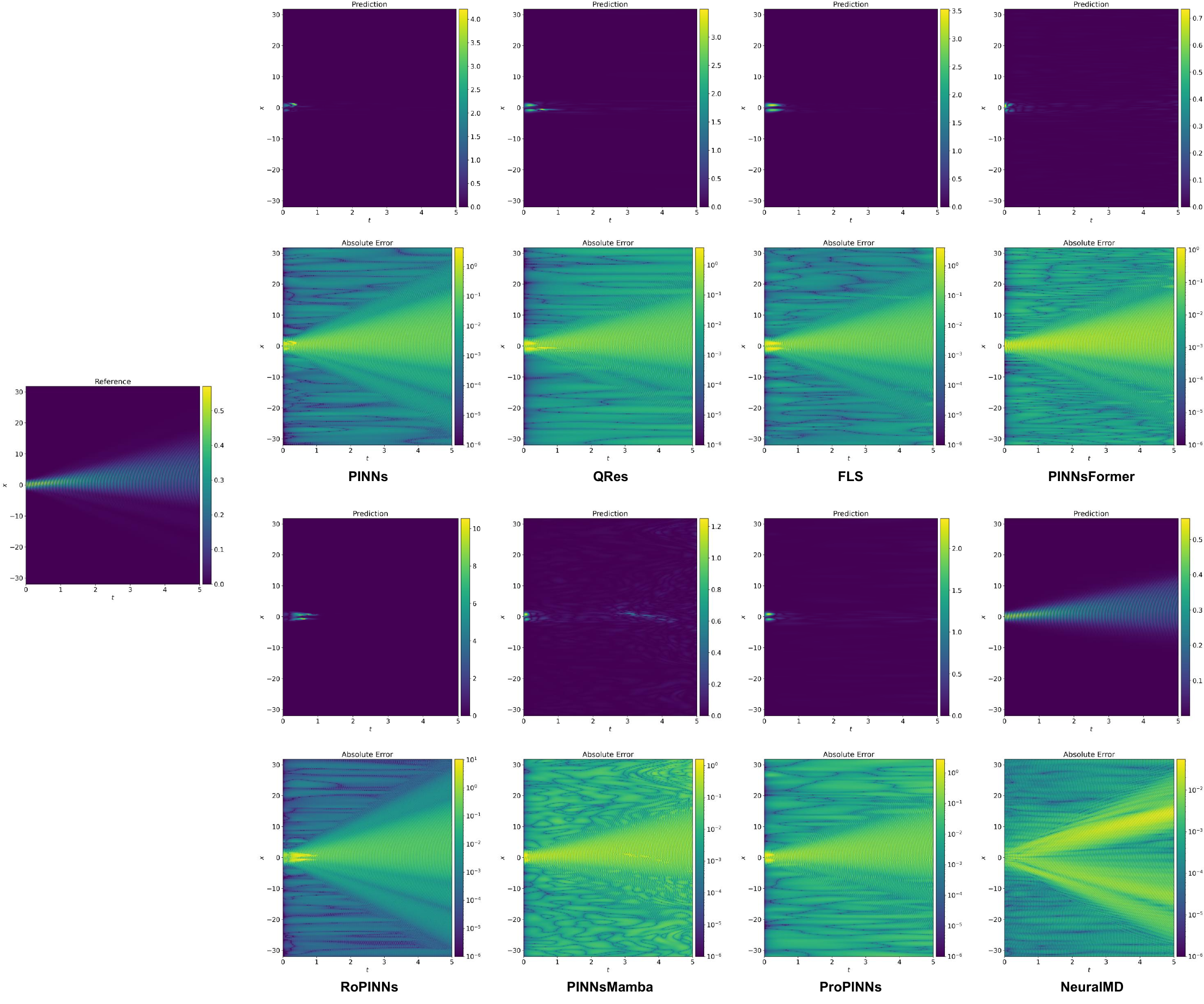}
    \vspace{-8pt}
    \caption{The prediction solution of NeuralMD and baselines for $\varepsilon=0.1,T=5.0$.}
    \label{fig:1d_eps_0.1}
\end{figure}

To compare model performance, Figure~\ref{fig:1d_eps_0.8}--\ref{fig:1d_eps_0.1} shows reference solutions, predictions, and absolute errors.
For $\varepsilon = 0.8$ (low-frequency temporal oscillation), Figure~\ref{fig:1d_eps_0.8} shows that several baseline methods (e.g., PINNs, QRes, PirateNet) produce predictions with visible errors, while NeuralMD has smaller localized discrepancies near wave peaks. 
For $\varepsilon = 0.5$ (intermediate-frequency temporal oscillation), Figure~\ref{fig:1d_eps_0.5} shows that several baseline methods fail to track the oscillatory wave structure over the full time interval.  NeuralMD maintains lower visual error across the plotted spatiotemporal domain.
For $\varepsilon = 0.1$, the tested baselines show large errors.
NeuralMD captures high-frequency oscillations more accurately via decomposition and WKB reconstruction.

\begin{table*}[t]
	\caption{Performance comparison of different PINNs architectures on the whole regime NKGE. 
    Both rMAE and rRMSE are recorded. Smaller values indicate better performance. 
    For clarity, the best result is in bold and the second best is underlined.  
    }
	\label{tab:compare_results}
	\vspace{-10pt}
	\vskip 0.15in
	\centering
	\begin{small}
	\renewcommand{\multirowsetup}{\centering}
	\setlength{\tabcolsep}{6pt}
	\begin{tabular}{l|cccccccc}
		\toprule
        \multirow{3}{*}{Model} & \multicolumn{2}{c}{$\varepsilon=0.8$} & \multicolumn{2}{c}{$\varepsilon=0.5$} & \multicolumn{2}{c}{$\varepsilon=0.1$} & \multicolumn{2}{c}{$\varepsilon=0.01$}  \\
        \cmidrule(lr){2-3}\cmidrule(lr){4-5}\cmidrule(lr){6-7}\cmidrule(lr){8-9}
		& $\mathrm{rMAE}$ & $\mathrm{rRMSE}$ & $\mathrm{rMAE}$ & $\mathrm{rRMSE}$ & $\mathrm{rMAE}$ & $\mathrm{rRMSE}$ & $\mathrm{rMAE}$ & $\mathrm{rRMSE}$ \\
		\midrule
        PINNs & 0.811 & 0.809 & 0.883 & 0.938 & 1.407 & 1.540 & 644.918 & 811.071 \\
        QRes & 0.761 & 0.904 & 0.863 & 0.945 & 1.943 & 1.821 & 1128.042 & 1215.761  \\
        FLS & 0.022 & 0.021 & 0.891 & 0.940 & 2.086 & 2.422 & \underline{22.665} & \underline{31.897}  \\
        PirateNet & 0.692 & 0.789 & 0.851 & 0.923 & 2.789 & 2.941 & 1449.254 & 1678.351  \\
        MSPINNs & 0.766 & 0.835 & 0.854 & 0.921 & 1.342 & 1.458 & 156.662 & 214.652  \\
        PhasePINNs & 0.673 & 0.754 & 0.784 & 0.855 & 1.124 & 1.246 & 121.890 & 189.706  \\
        CausalPINNs & 0.534 & 0.521 & 1.382 & 1.198 & 25.628 & 17.085 & 2479.044 & 2153.672  \\
        PINNsFormer & 0.022 & 0.021 & 0.124 & 0.146 & 1.028 & \underline{1.012} & 30.362 & 46.366  \\
        SetPINNs & 0.006 & 0.008 & 0.107 & 0.124 &  \underline{1.026} & 1.020 & 30.561 & 45.286  \\
        RoPINNs & 0.022 & 0.021 & 1.014 & 0.953 & 4.184 & 5.094 & 59.345 & 69.012  \\
        ProPINNs & \underline{0.004} & \underline{0.005} & \underline{0.050} & \underline{0.053} & 1.417 & 1.277 & 449.188 & 412.412  \\
        PINNsMamba & 0.027 & 0.037 & 0.235 & 0.264 & 1.163 & 1.042 & 28.562 & 35.267  \\
        \midrule
        \textbf{NeuralMD} & \textbf{0.002} & \textbf{0.003} & \textbf{0.008} & \textbf{0.010} & \textbf{0.005} & \textbf{0.007} & \textbf{0.006} & \textbf{0.008}  \\
        Improvement & 50.0\% & 40.0\% & 84.0\% & 81.1\% & 99.5\% & 99.3\% & 99.9\% & 99.9\% \\
		\bottomrule
        \end{tabular}
    \end{small}
\vspace{-10pt}
\end{table*}

\subsection{2D nonlinear Klein-Gordon equation in whole regime}
In this section, we explore the wave interactions in 2D.
We take $d=2$ and $\lambda=1$ in the NKGE (\ref{eqn:KG}) and choose the initial data as
\begin{equation*}
\begin{split}
 &\phi_1(x,y)=\exp{(-(x+2)^2-y^2)}+\exp{(-(x-2)^2-y^2)},\\
 & \phi_2(x,y)=\exp{(-x^2-y^2)},\qquad (x,y)\in{\mathbb R}^2.
\end{split}
\end{equation*}
The problem is solved numerically on a bounded computational domain
$\Omega=(-16,16)\times(-16,16)$ with the periodic boundary condition.

\begin{figure}[!htb]
    \centering
    \includegraphics[width=0.9\textwidth]{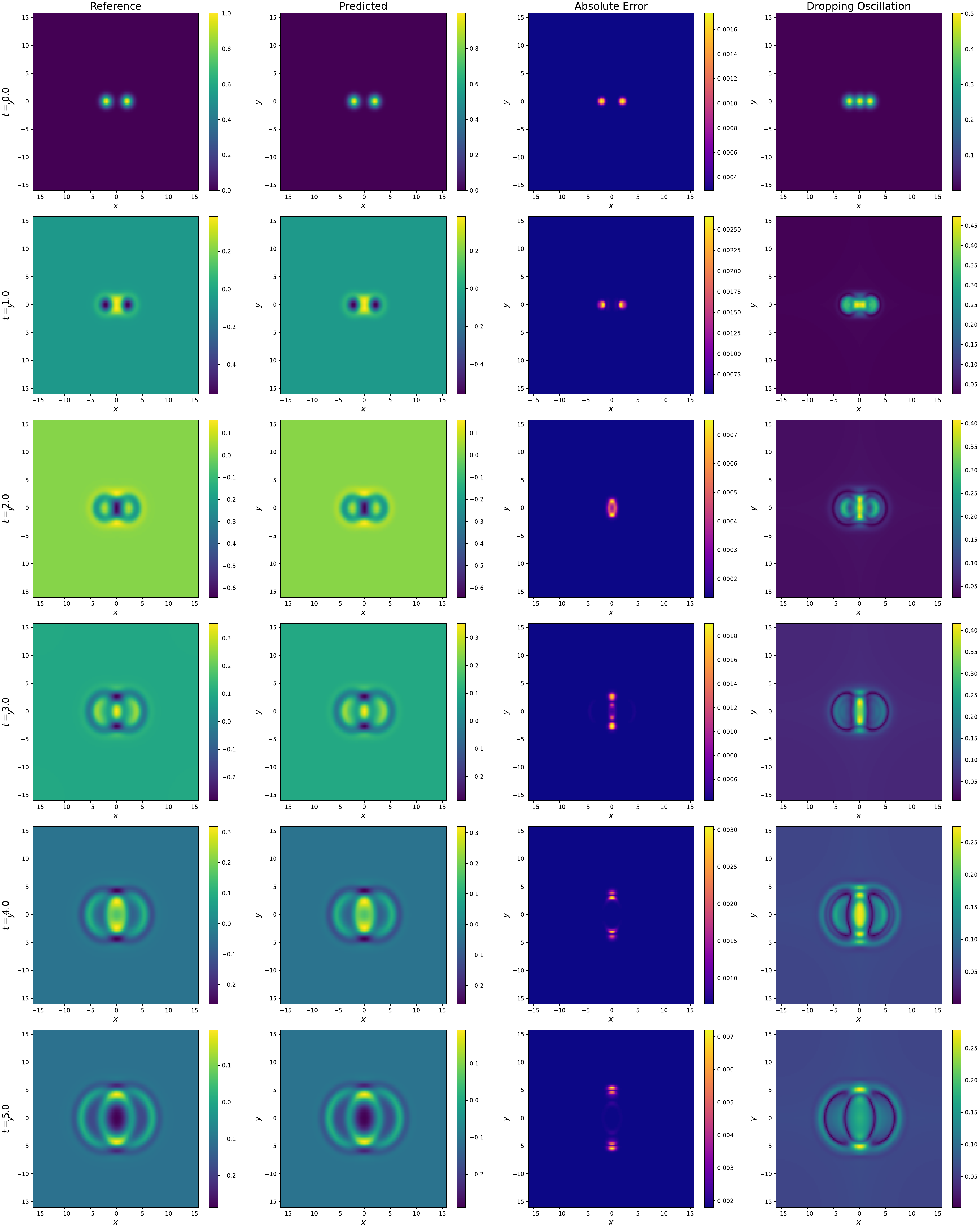}
    \vspace{-8pt}
    \caption{The prediction solution of NeuralMD in 2D for $\varepsilon=0.8,T=5.0$.}
    \label{fig:2d_eps_0.8}
\end{figure}

\begin{figure}[!htb]
    \centering
    \includegraphics[width=0.9\textwidth]{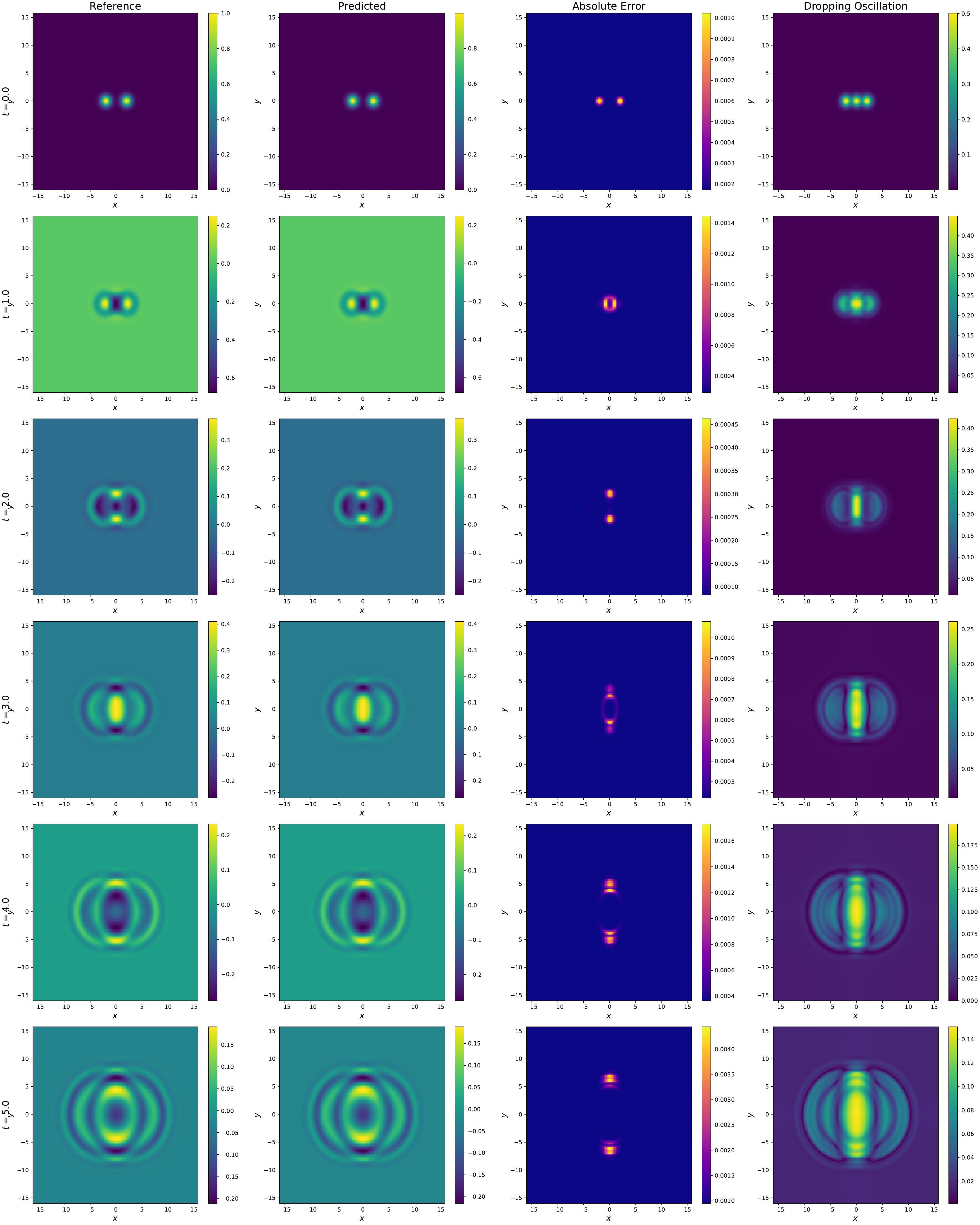}
    \vspace{-8pt}
    \caption{The prediction solution of NeuralMD in 2D for $\varepsilon=0.5,T=5.0$.}
    \label{fig:2d_eps_0.5}
\end{figure}

\begin{figure}[!htb]
    \centering
    \includegraphics[width=0.9\textwidth]{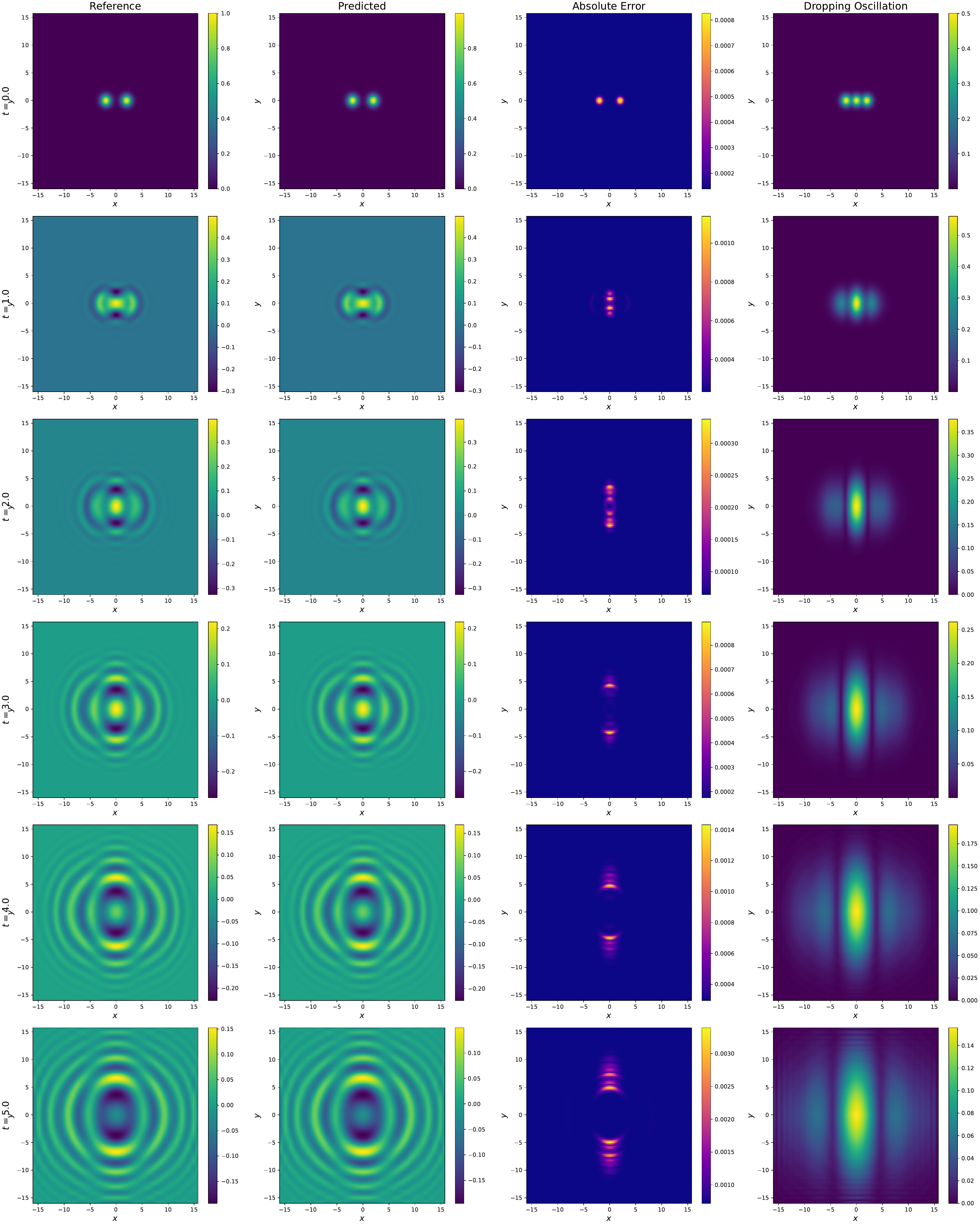}
    \vspace{-8pt}
    \caption{The prediction solution of NeuralMD in 2D for $\varepsilon=0.1,T=5.0$.}
    \label{fig:2d_eps_0.1}
\end{figure}

Figure~\ref{fig:2d_eps_0.8} shows the 2D NeuralMD prediction for $\varepsilon=0.8$.
The reference solution contains two interacting Gaussian wave packets.
NeuralMD captures both wave amplitude and phase.
The error is mainly localized near interaction regions. 
For $\varepsilon=0.5$, Figure~\ref{fig:2d_eps_0.5} shows accurate prediction under stronger temporal oscillation.
The wave structure remains well preserved, and the absolute error stays low.
For $\varepsilon=0.1$, Figure~\ref{fig:2d_eps_0.1} shows qualitative high-frequency reconstruction through WKB expansion.

\subsection{3D nonlinear Klein-Gordon equation in whole regime}
In this section, we explore the wave interactions in 3D.
We take $d=3$ and $\lambda=1$ in the NKGE (\ref{eqn:KG}) and choose the initial data as
\begin{equation*}
\begin{split}
 &\phi_1(x,y,z)=2\exp{(-x^2-2y^2-3z^2)},\\
 & \phi_2(x,y,z)=\exp{(-(x+0.5)^2-y^2-z^2)},\quad (x,y,z)\in{\mathbb R}^3.
\end{split}
\end{equation*}

\begin{figure}[!htb]
    \centering
    \includegraphics[width=1.0\textwidth]{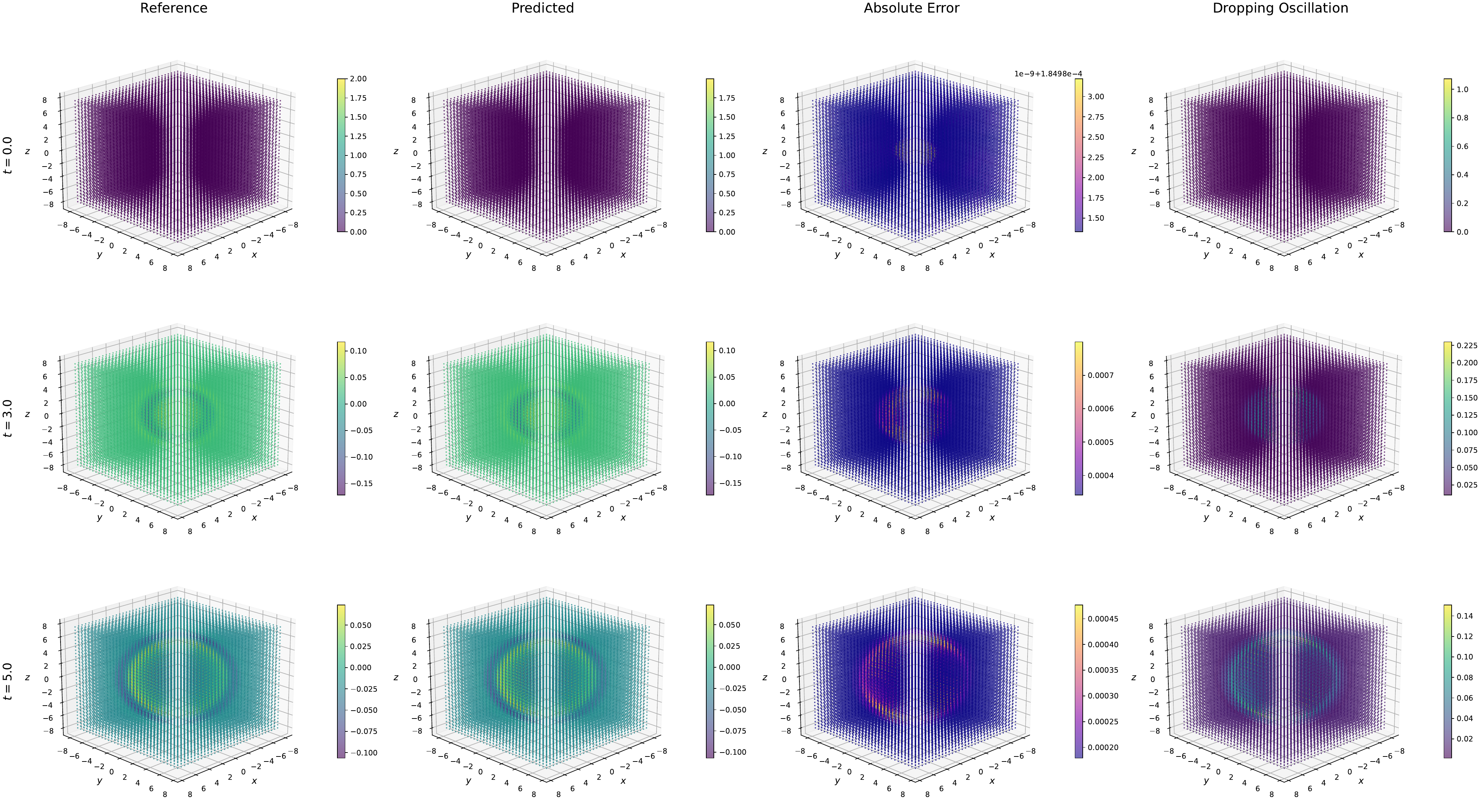}
    \vspace{-8pt}
    \caption{The prediction solution of NeuralMD in 3D for $\varepsilon=0.8,T=5.0$.}
    \label{fig:3d_eps_0.8}
\end{figure}

\begin{figure}[!htb]
    \centering
    \includegraphics[width=1.0\textwidth]{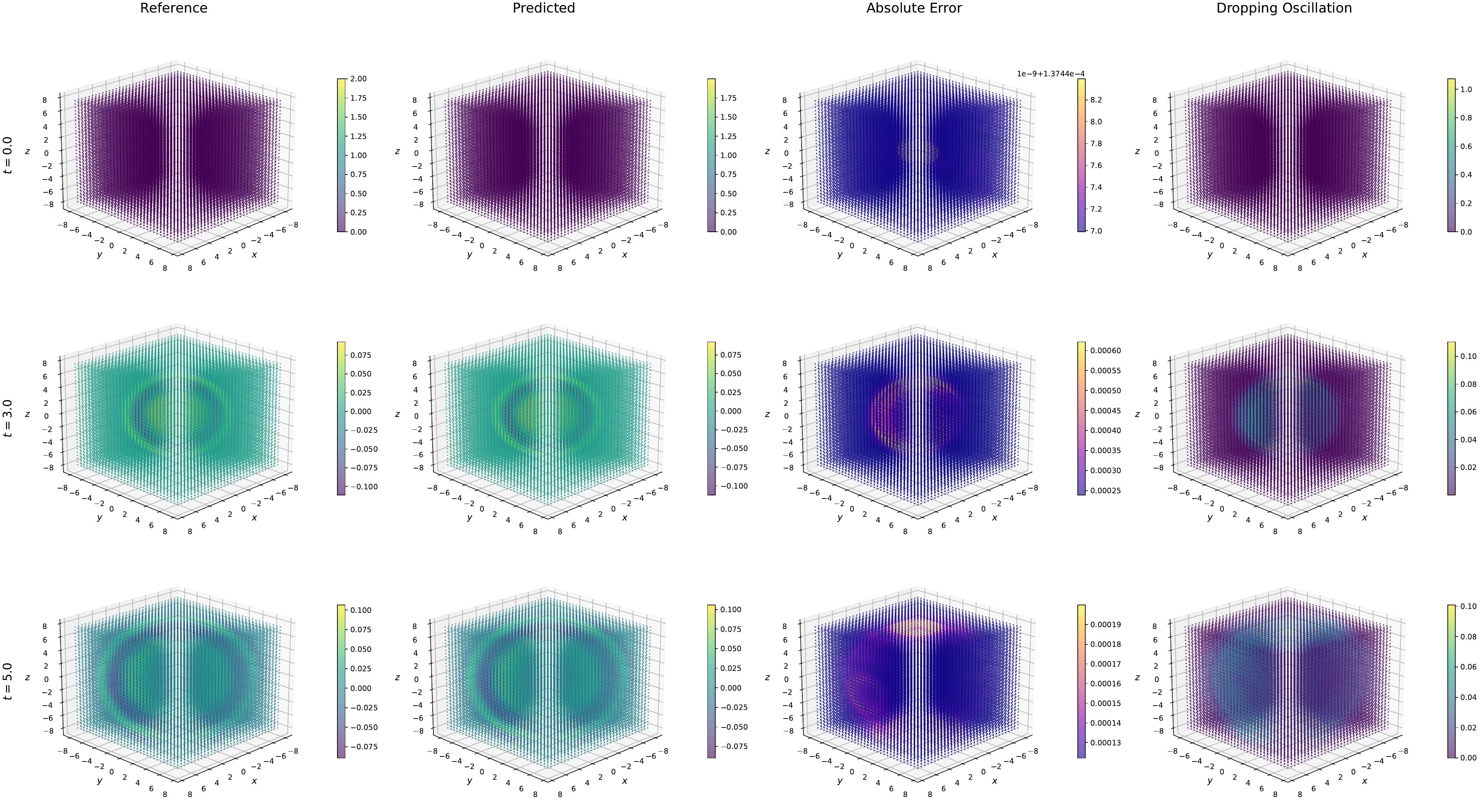}
    \vspace{-8pt}
    \caption{The prediction solution of NeuralMD in 3D for $\varepsilon=0.5,T=5.0$.}
    \label{fig:3d_eps_0.5}
\end{figure}

\begin{figure}[!htb]
    \centering
    \includegraphics[width=1.0\textwidth]{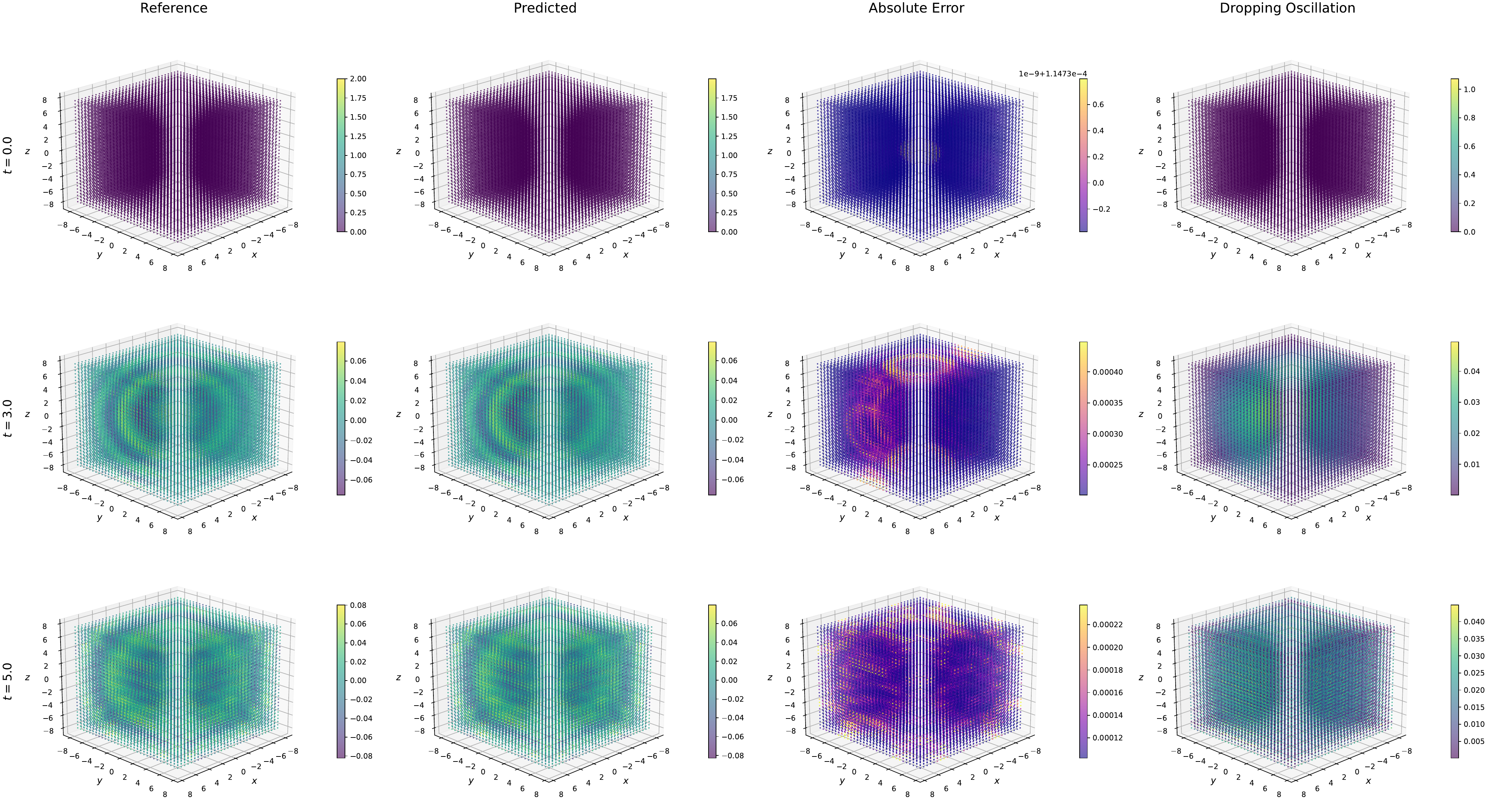}
    \vspace{-8pt}
    \caption{The prediction solution of NeuralMD in 3D for $\varepsilon=0.1,T=5.0$.}
    \label{fig:3d_eps_0.1}
\end{figure}

Figure~\ref{fig:3d_eps_0.8} shows the 3D NeuralMD prediction for $\varepsilon=0.8$.
The visualization shows anisotropic spreading of a Gaussian wave packet.
NeuralMD closely matches the reference solution across cross-sectional views.
For $\varepsilon=0.5$, Figure~\ref{fig:3d_eps_0.5} shows that accuracy remains stable.
The 3D wave structure is well preserved under stronger temporal oscillation.
For $\varepsilon=0.1$, Figure~\ref{fig:3d_eps_0.1} shows qualitative high-frequency reconstruction.
This reconstruction is obtained through the multiscale decomposition.

\subsection{Phase separation and reconstruction}
To validate carrier separation, we test $\varepsilon=0.8$, $\varepsilon=0.5$, and $\varepsilon=0.1$.

\begin{figure}[!htb]
    \centering
    \includegraphics[width=0.95\textwidth]{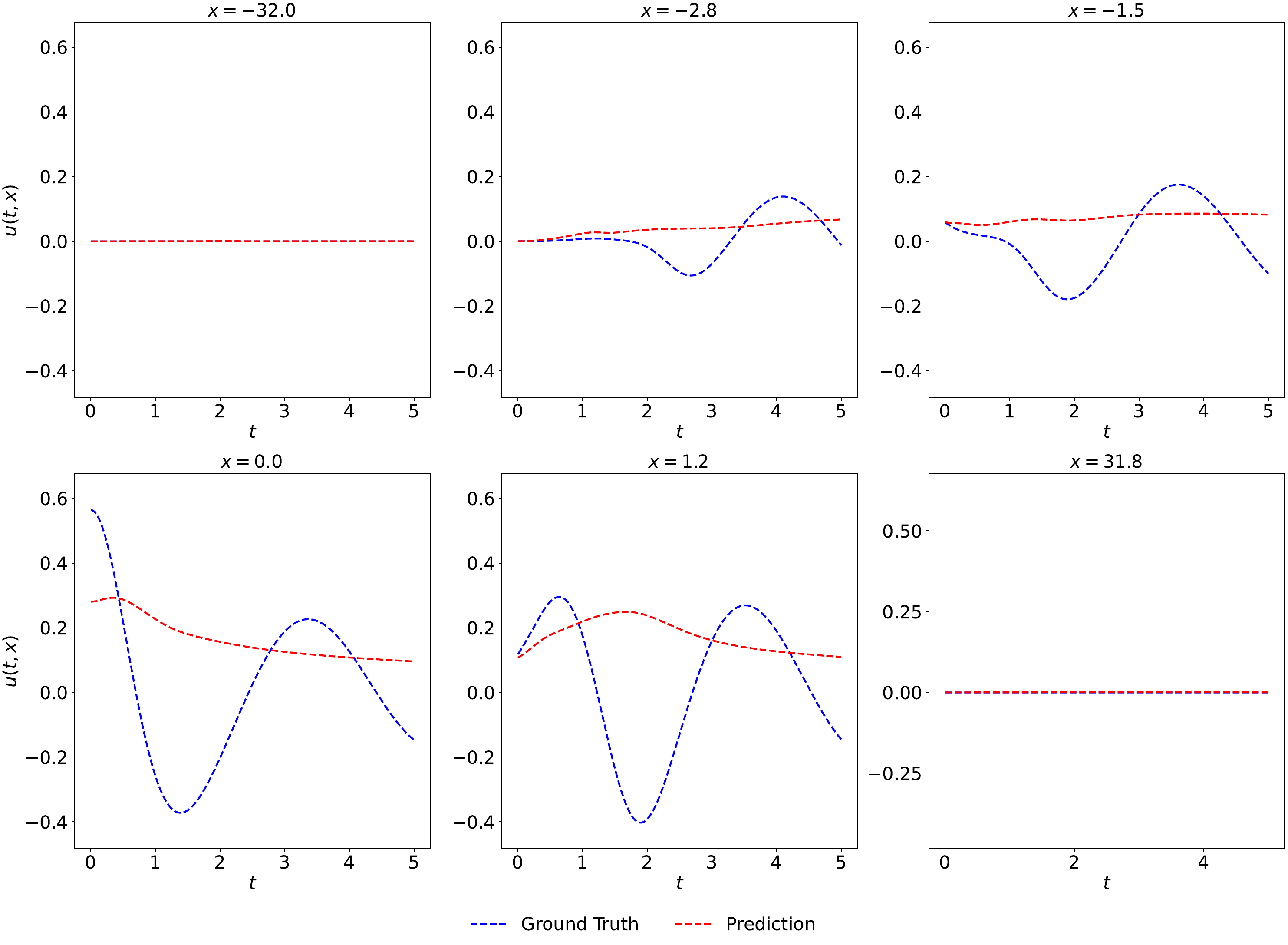}
    \vspace{-8pt}
    \caption{Carrier-separation results of NeuralMD for $\varepsilon=0.8,T=5.0$.}
    \label{fig:drop1d_eps_0.8}
\end{figure}

\begin{figure}[!htb]
    \centering
    \includegraphics[width=0.95\textwidth]{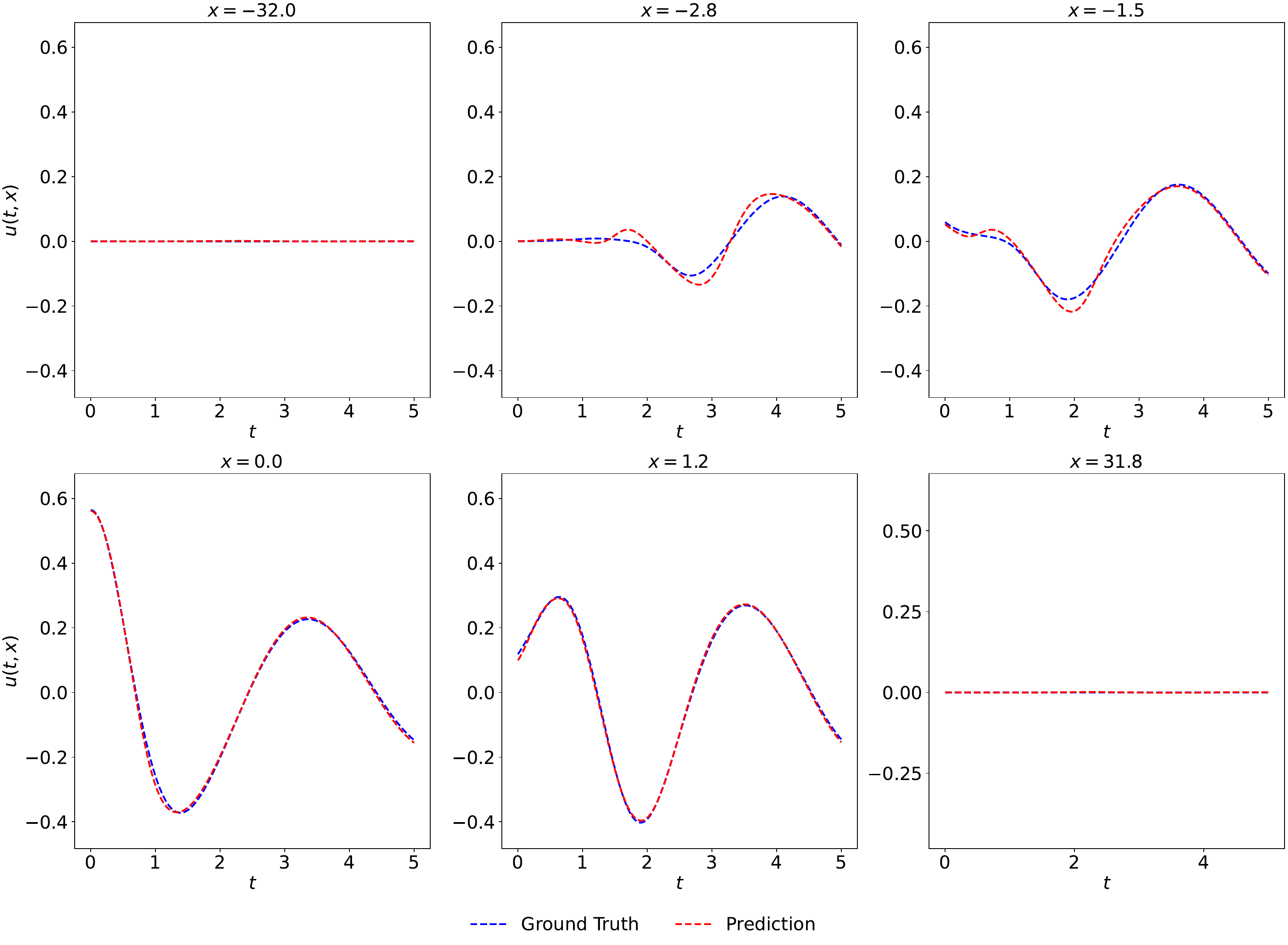}
    \vspace{-8pt}
    \caption{Results of reconstructing time oscillation of NeuralMD for $\varepsilon=0.8,T=5.0$.}
    \label{fig:recont1d_eps_0.8}
\end{figure}

For $\varepsilon=0.8$, Figure~\ref{fig:drop1d_eps_0.8} shows carrier separation at $x=-2,0,2$.
Blue curves are original oscillatory solutions.
Orange dashed curves are smoothed envelope predictions after carrier separation.  
At this stage, the impact of time oscillation is minimal, and the process is primarily dominated by the remainder amplitude compensation.  
A notable observation is that while time oscillation remains at each position, the oscillatory behavior is modulated into a smooth curve.  
Figure~\ref{fig:recont1d_eps_0.8} shows WKB reconstruction with remainder terms.
The reconstructed solutions closely match the reference slices.
The $L^2$ relative error (rRMSE) is $0.064$.

\begin{figure}[!htb]
    \centering
    \includegraphics[width=0.95\textwidth]{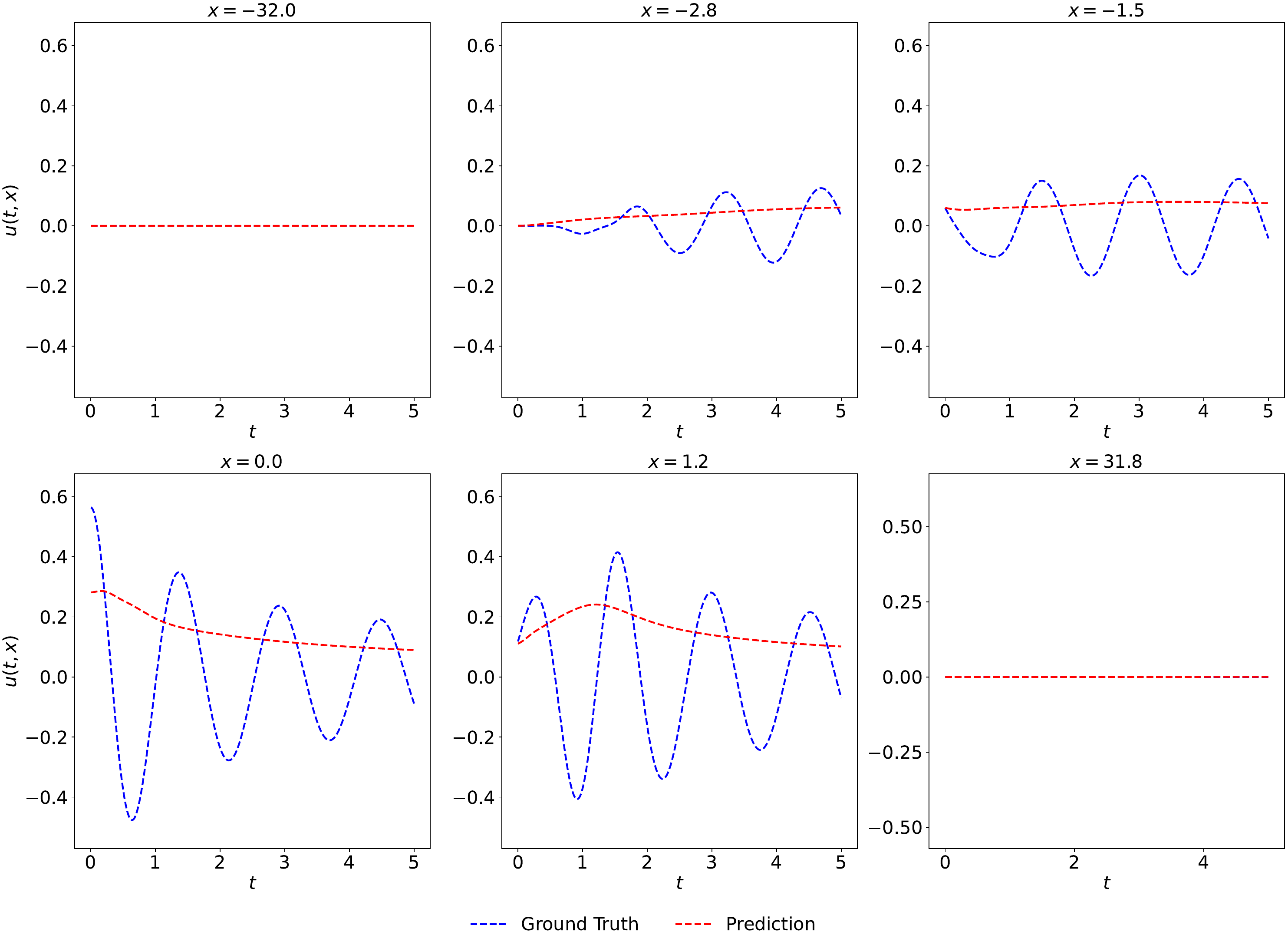}
    \vspace{-8pt}
    \caption{Carrier-separation results of NeuralMD for $\varepsilon=0.5,T=5.0$.}
    \label{fig:drop1d_eps_0.5}
\end{figure}

\begin{figure}[!htb]
    \centering
    \includegraphics[width=0.95\textwidth]{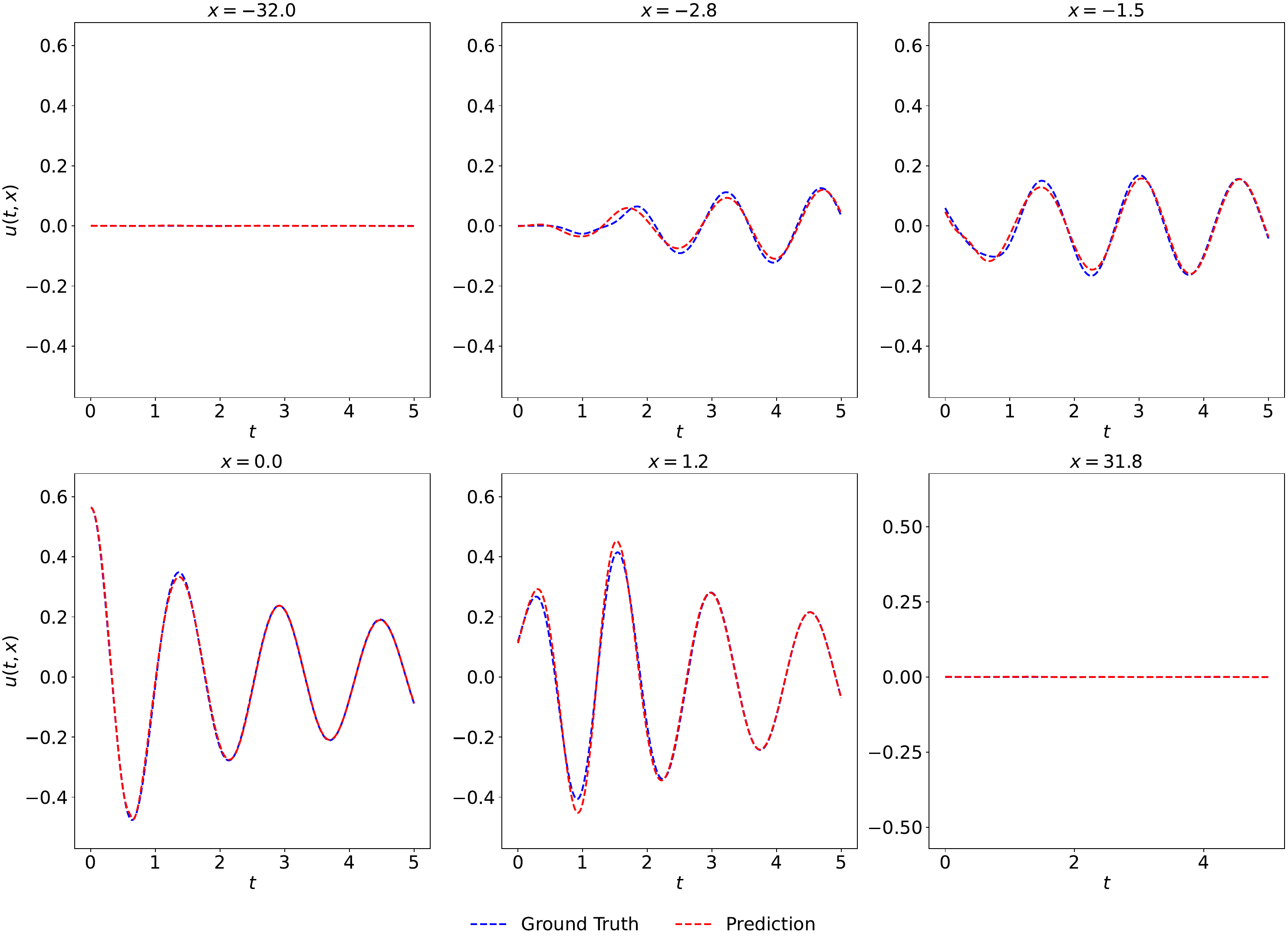}
    \vspace{-8pt}
    \caption{Results of reconstructing time oscillation of NeuralMD for $\varepsilon=0.5,T=5.0$.}
    \label{fig:recont1d_eps_0.5}
\end{figure}

For $\varepsilon=0.5$, Figure~\ref{fig:drop1d_eps_0.5} shows time evolution at $x=-2,0,2$.
The oscillatory solutions are represented through smoother modulated envelopes.  
Figure~\ref{fig:recont1d_eps_0.5} shows accurate WKB reconstruction across all positions.
The method achieves an $L^2$ relative error of $0.099$.

\begin{figure}[!htb]
    \centering
    \includegraphics[width=0.95\textwidth]{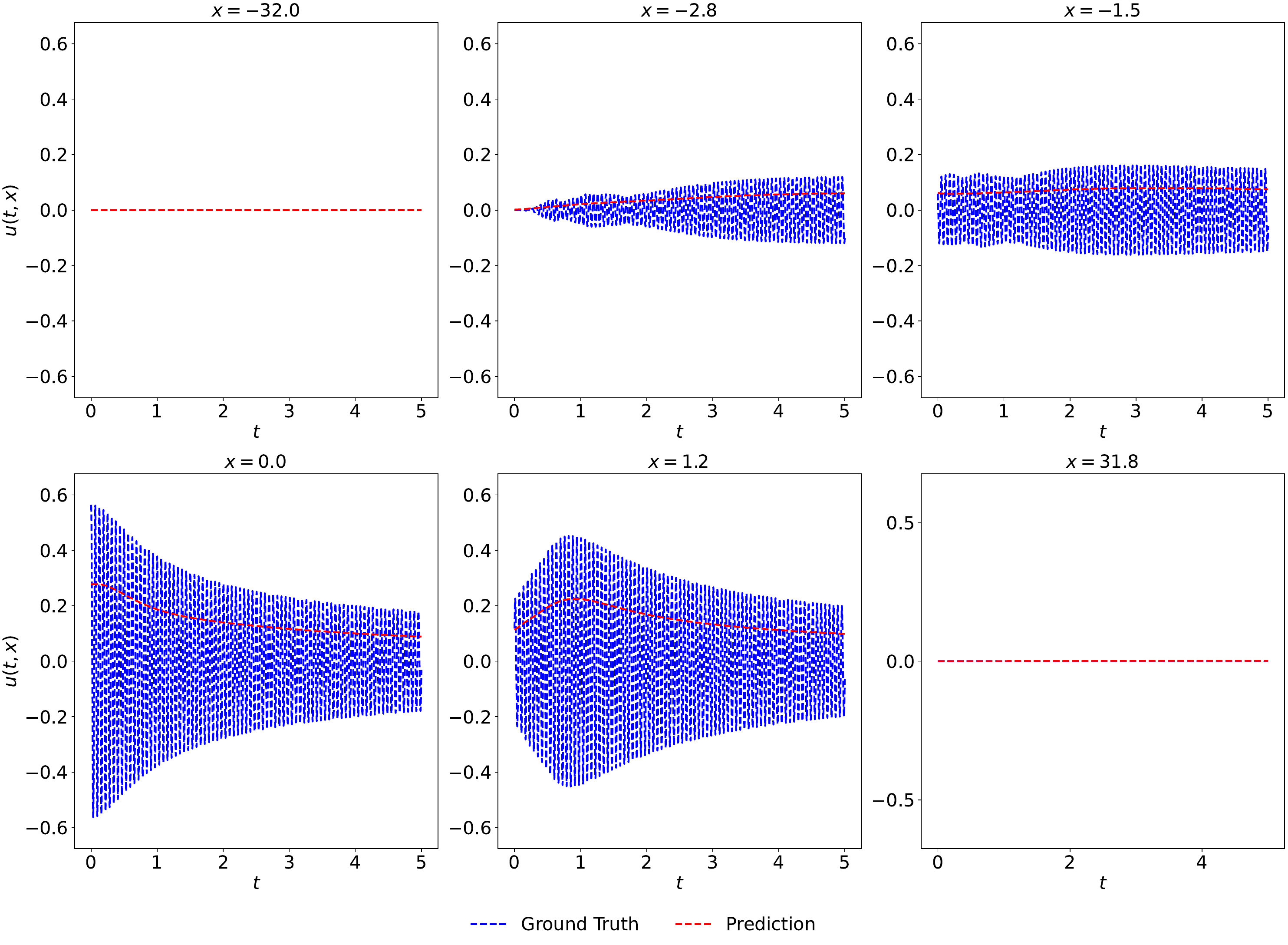}
    \vspace{-8pt}
    \caption{Carrier-separation results of NeuralMD for $\varepsilon=0.1,T=5.0$.}
    \label{fig:drop1d_eps_0.1}
\end{figure}

\begin{figure}[!htb]
    \centering
    \includegraphics[width=0.95\textwidth]{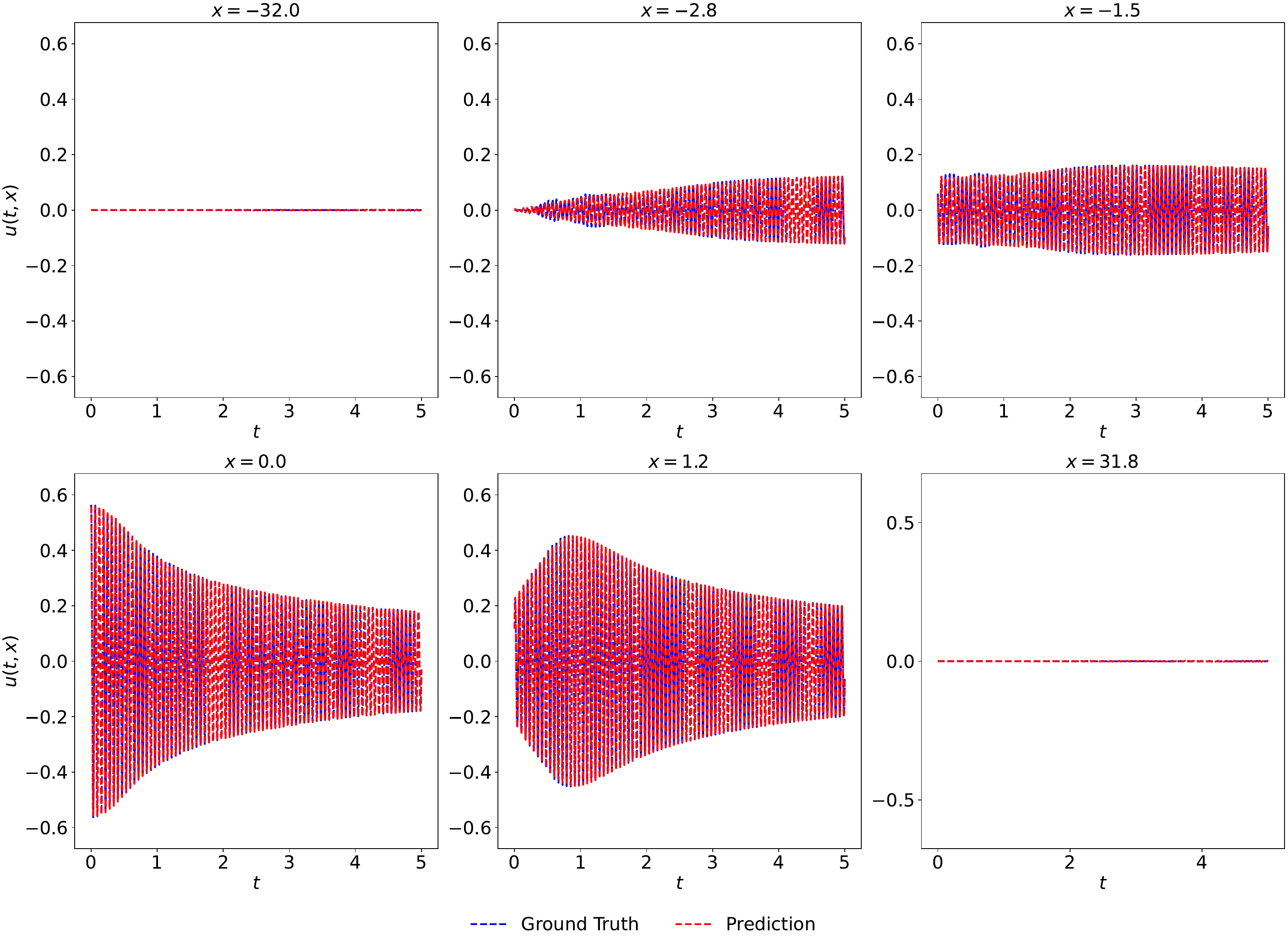}
    \vspace{-8pt}
    \caption{Results of reconstructing time oscillation of NeuralMD for $\varepsilon=0.1,T=5.0$.}
    \label{fig:recont1d_eps_0.1}
\end{figure}

For $\varepsilon=0.1$, Figure~\ref{fig:drop1d_eps_0.1} shows high-frequency temporal oscillation.
NeuralMD factors the high-frequency carrier through multiscale decomposition.
This transforms oscillatory signals into smoother neural representations.  
NeuralMD reconstructs temporal oscillation through the WKB formula; see Figure~\ref{fig:recont1d_eps_0.1}.
The WKB reconstruction with remainder terms achieves rRMSE $0.069$ in this test.

\begin{figure}[!htb]
    \centering
    \includegraphics[width=0.85\textwidth]{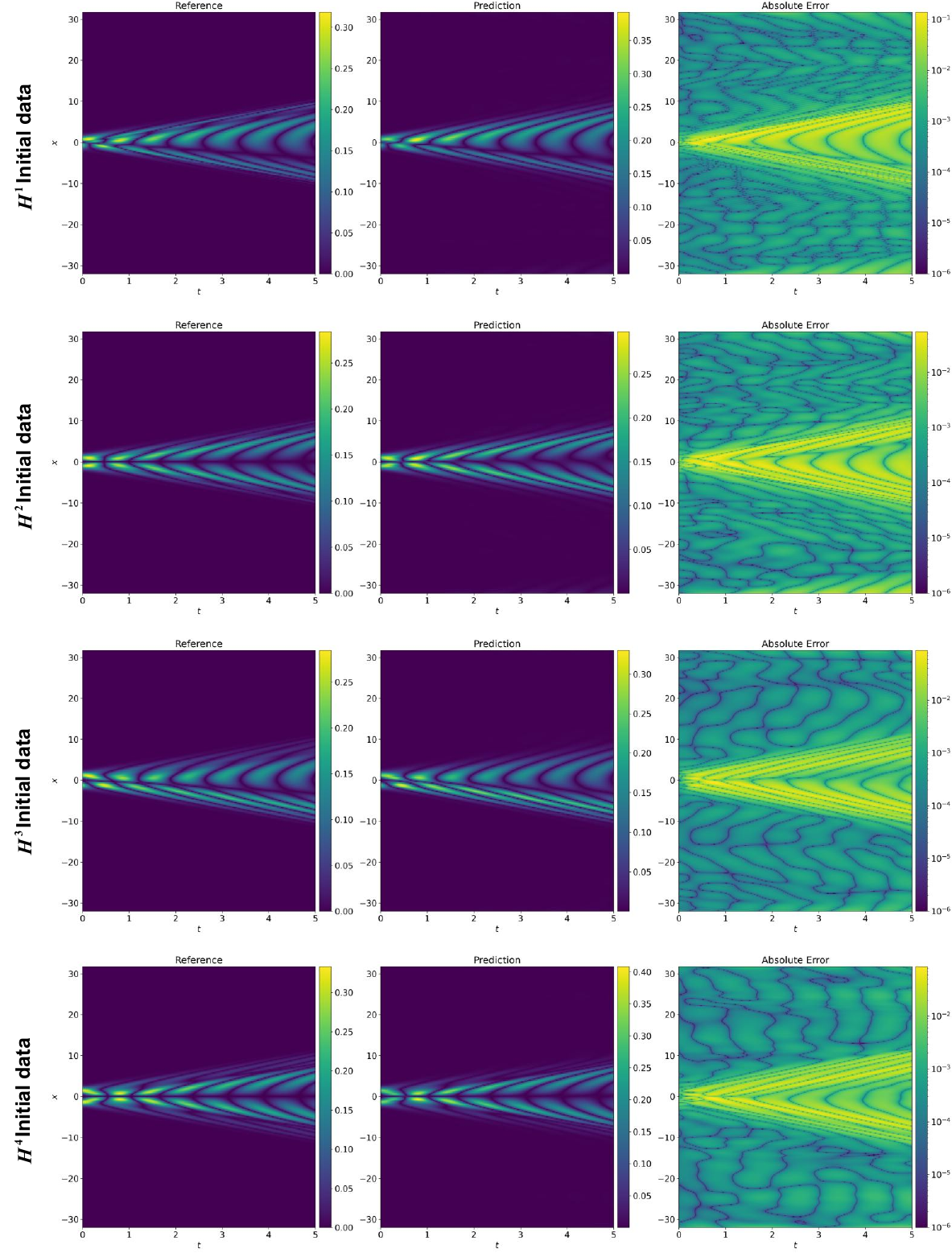}
    \vspace{-8pt}
    \caption{The prediction solution of NeuralMD for $\varepsilon=0.5$ under initial data with different regularity.}
    \label{fig:H_initial_data}
\end{figure}

Furthermore, we evaluate carrier separation under initial data with varying regularities, as illustrated in Figure~\ref{fig:H_initial_data}. 
For $H^1$ to $H^4$ initial data, NeuralMD maintains stable reconstruction errors.
These data represent solutions with different Sobolev regularity. 
These results suggest robustness across the tested levels of solution regularity.

\subsection{Efficiency comparison}
\begin{figure}[!htb]
    \centering
    \includegraphics[width=1.0\textwidth]{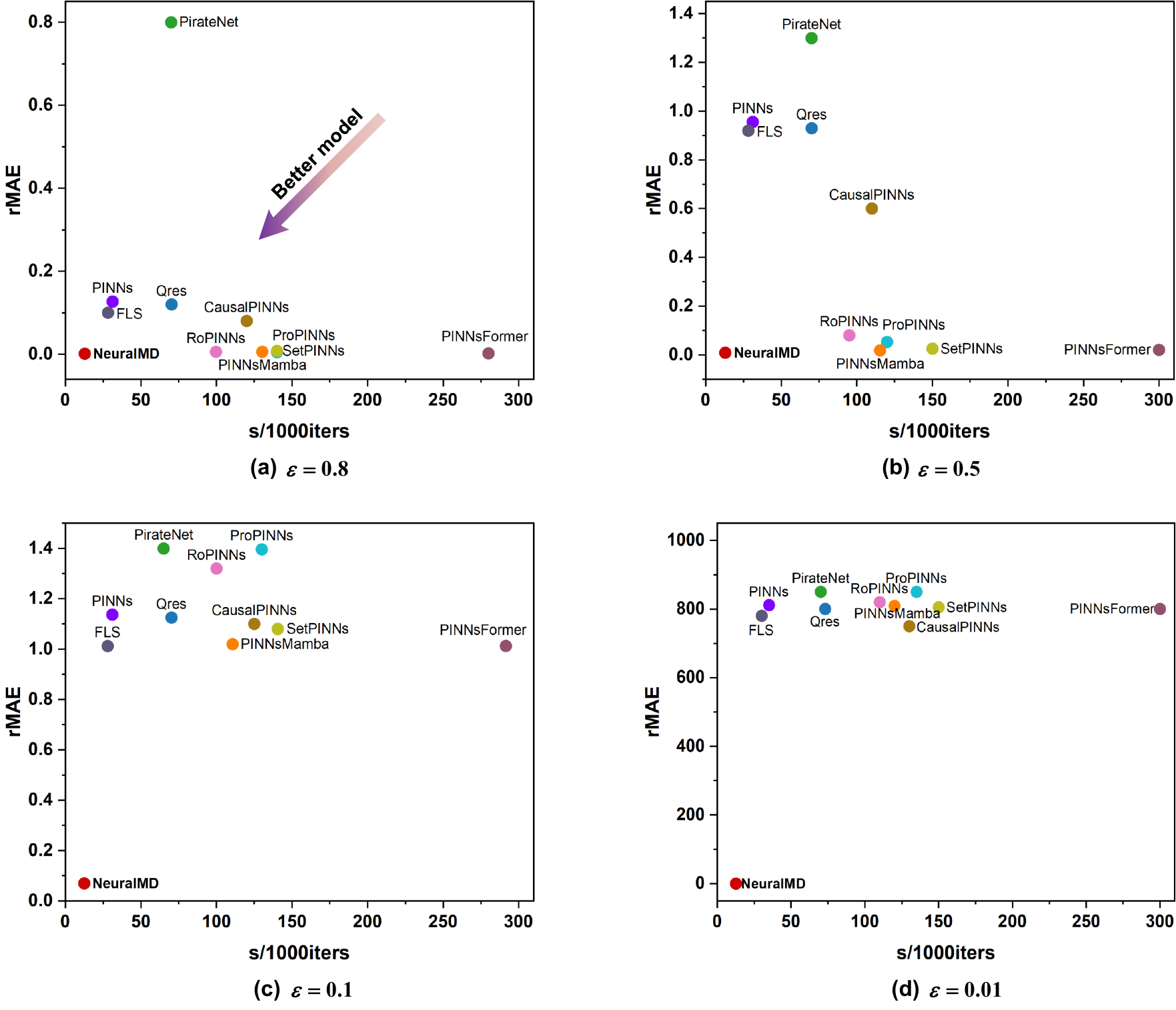}
    \vspace{-8pt}
    \caption{Efficiency comparison for a temporally oscillatory problem under different $\varepsilon$.}
    \label{fig:efficient}
\end{figure}

To assess practicability, Figure~\ref{fig:efficient} compares efficiency across representative $\varepsilon$ regimes. 
The left panel shows rMAE versus training time for each method.
In these experiments, NeuralMD reaches lower error within the measured training budget.
As $\varepsilon\to0$, several baseline errors increase rapidly.
For $\varepsilon=0.01$, the vanilla PINN rMAE reaches 800.
This suggests spectral bias and propagation failure under extreme temporal oscillation.
In contrast to the baselines, NeuralMD maintains errors below 0.01 in these runs. The error curves for NeuralMD show rapid convergence and stable final accuracy on the tested cases.
The right panel shows that NeuralMD is about 2--3$\times$ faster than Transformer-based baselines.
These include PINNsFormer~\citep{zhao2023former}, SetPINN~\citep{nagda2024set}, and PINNMamba~\citep{xu2025mamba}.  
With lightweight projection layers, NeuralMD remains comparable to single-point PINN variants.
These include QRes, FLS, CausalPINNs, and PirateNet.
Overall, NeuralMD achieves a favorable performance-efficiency trade-off on the tested temporally oscillatory PDE benchmarks.

\section{Conclusions}
\label{sec:conclusion}
In this work, we proposed NeuralMD for temporally oscillatory nonlinear Klein--Gordon equations.
The method targets representative values of $\varepsilon\in(0,1]$.
NeuralMD adopts a two-stage pretraining strategy.
The first stage learns the mildly oscillatory NLSW modulation.
The second stage learns a small-amplitude remainder equation.
This remainder corrects residual oscillatory components. 
The full NKGE solution is reconstructed through a WKB expansion.
An error-based criterion decides whether the remainder contribution is retained.

To alleviate propagation failure arising when training the modulated NLSW and oscillatory remainder networks, we introduced a gated gradient correlation metric. 
We also propose gated residual sampling to couple temporal collocation points.
This strengthens the causal structure of the training dynamics.
We further extend NeuralMD to an interpretable edge-activation variant.
Nodes perform only summation, while edges carry learnable activation functions. 
In this formulation, the mechanisms of carrier separation and dynamic remainder-amplitude compensation become structurally interpretable at the level of the learned operators.
Experiments cover different dimensions, initial-data regularities, and long-time prediction settings.
They indicate lower spectral-bias and propagation-failure errors than tested collocation baselines.

Several directions merit further investigation. 
First, joint training over $\varepsilon\in(0,1]$ is a natural extension.
Such training may face severe gradient conflicts across regimes. 
Second, NeuralMD could be extended to semi-classical nonlinear Schr\"odinger equations.
There, simultaneous temporal and spatial oscillation may benefit from carrier separation.
Finally, the WKB cascade may support higher-order MTI schemes with learned operators.


\acks{The authors gratefully acknowledge financial support from the National Natural Science Foundation of China (12202157). 
We also thank the Exploration Foundation of the Key Laboratory of CNC Equipment Reliability, Ministry of Education.
We thank the National Key Laboratory of Automotive Chassis Integration and Bionics, School of Mechanical and Aerospace Engineering, Jilin University.}

\newpage

\appendix
\section{Stability and Monte Carlo diagnostics for NeuralMD}
\label{app:stability_diagnostics}
The theoretical results in the main text are deliberately restricted.
They can be verified on the scalar remainder-selection path or through
local descent inequalities.  For completeness, this appendix records the
standard diagnostics we use for gated multiscale temporal sampling.  These are
not claimed as global optimization guarantees for the nonconvex neural
training problem.

Let $\mathcal L_{\rm gate}$ denote the gated multiscale empirical objective.
It averages pointwise losses over temporal perturbations.
The weights are nonnegative, and the gate satisfies $0\le h\le1$.
Assume the pointwise loss is $L$-Lipschitz.
Also assume the stochastic-gradient update is uniformly stable.
Then the expected generalization gap has the form
\begin{equation*}
  \mathbb E[\mathcal E_{\rm gen}]
  \le
  \frac{2L^2}{|\mathcal S|}
  \sum_{t=0}^{T-1}\alpha_t ,
\end{equation*}
up to constants that depend on the precise sampling and replacement model.
The gate affects only the effective Lipschitz and variance constants.
It does not imply a negative or unconditional improvement term.

Similarly, assume $\mathcal L_{\rm gate}$ is lower bounded.
Assume it has an $H_{\rm gate}$-Lipschitz gradient.
Assume also that the stochastic gated gradient has bounded second moment.
Then SGD with $\alpha_t=(t+1)^{-1/2}$ gives the standard diagnostic
\begin{equation*}
  \mathbb E\|\nabla_\theta\mathcal L_{\rm gate}(\theta_{I_T})\|^2
  \le
  \frac{\mathcal L_{\rm gate}(\theta_0)-\inf_\theta\mathcal L_{\rm gate}(\theta)}
       {\sum_{t=0}^{T-1}\alpha_t}
  +
  \frac{H_{\rm gate}\sigma_g^2}{2}
  \frac{\sum_{t=0}^{T-1}\alpha_t^2}
       {\sum_{t=0}^{T-1}\alpha_t},
\end{equation*}
where $\Pr(I_T=t)=\alpha_t/\sum_{k=0}^{T-1}\alpha_k$.  With this step size the
right-hand side is $O(\log(T)/\sqrt{T})$.  This is a diagnostic for a randomly
selected iterate, not a guarantee for the final iterate.

Finally, take $M$ independent temporal perturbation samples.
Let their gradients be $\widehat g_1,\ldots,\widehat g_M$.
Assume their covariance trace is bounded by $\sigma_g^2$.
The averaged estimator
$\bar g_M=M^{-1}\sum_{j=1}^M\widehat g_j$ satisfies
\begin{equation*}
  \mathbb E\|\bar g_M-\mathbb E\widehat g\|^2
  \le \frac{\sigma_g^2}{M}.
\end{equation*}
This elementary variance reduction is the only Monte Carlo claim used to
interpret the number of temporal perturbations.

\section{Corrected generalization analysis for NeuralMD}
\label{app:corrected_generalization}
This section gives a stability-based generalization statement for the gated
multiscale temporal objective.  The result is conditional on standard
algorithmic stability assumptions.  It is not a global optimization guarantee
for the nonconvex neural network.
All results here apply to the idealized expected gated objective below.
They also apply to its finite Monte Carlo estimator under stated sampling assumptions.
They do not prove convergence of the final trained NeuralMD iterate.
They also do not imply a PDE solution-error bound without main-text compatibility assumptions.

\begin{table}[t]
\caption{Notation used in Appendix~\ref{app:corrected_generalization}--\ref{app:corrected_gradient_estimation}.}
\label{tab:appendix_notation}
\centering
\begin{tabular}{ll}
\toprule
Symbol & Meaning \\
\midrule
$s=(\bm{x},t)$ & space-time sample point from $Q_T$ \\
$\vartheta$ & random temporal perturbation index \\
$t_\vartheta$ & perturbed time generated from $(t,\vartheta)$ \\
$n$ & number of training samples in $\mathcal S$ \\
$K$ & number of optimization iterations \\
$M$ & number of Monte Carlo temporal perturbations \\
$I_K$ & random optimization-iterate index \\
$\bar\ell,\mathcal J,\widehat{\mathcal J}_{\mathcal S}$ & gated pointwise, population, and empirical objectives \\
$L_g,\beta_g,H_g$ & effective Lipschitz, smoothness, and gradient-Lipschitz constants \\
\bottomrule
\end{tabular}
\end{table}

Let $s=(\bm{x},t)$ denote a space-time collocation point sampled from $\mu$ on
$Q_T=\Omega\times[0,T]$.  Given $s$, let $\vartheta$ denote the random temporal
perturbation produced by the multiscale rule.
First, a scale $l$ is sampled with probability $\lambda_l$.
Then $t_\vartheta=t_\vartheta(t,\vartheta)$ is sampled in the corresponding neighborhood.
Boundary clipping is applied when needed.
The gate satisfies $0\le h(t_\vartheta)\le1$.  For a
pointwise residual or data loss $\ell_{\rm pt}(\theta;\bm{x},t_\vartheta)$, define
\begin{equation}
\label{eq:app_gated_loss}
  \bar\ell(\theta;s)
  =
  \mathbb E_{\vartheta\mid s}
  \bigl[
    h(t_\vartheta)\ell_{\rm pt}(\theta;\bm{x},t_\vartheta)
  \bigr].
\end{equation}
The population and empirical objectives are
\begin{equation}
\label{eq:app_pop_emp_loss}
  \mathcal J(\theta)=\mathbb E_{s\sim\mu}\bar\ell(\theta;s),
  \qquad
  \widehat{\mathcal J}_{\mathcal S}(\theta)
  =
  \frac{1}{n}\sum_{i=1}^{n}\bar\ell(\theta;s_i),
  \qquad
  \mathcal S=\{s_i\}_{i=1}^{n}.
\end{equation}

\begin{definition}[Uniform stability for the gated objective]
\label{def:app_uniform_stability}
A randomized training algorithm $\mathcal A$ is $\beta_n$-uniformly stable for
$\bar\ell$ under the following condition.
For any two datasets $\mathcal S$ and $\mathcal S'$ differing in one sample,
\begin{equation*}
  \sup_s
  \mathbb E_{\mathcal A}
  \left|
    \bar\ell(\mathcal A(\mathcal S);s)
    -
    \bar\ell(\mathcal A(\mathcal S');s)
  \right|
  \le \beta_n .
\end{equation*}
\end{definition}

\begin{theorem}[Expected generalization from stability]
\label{thm:app_stability_generalization}
Assume $\bar\ell(\mathcal A(\mathcal S);s)$ is integrable for
$s\sim\mu$ and for the algorithmic randomness.  If $\mathcal A$ is
$\beta_n$-uniformly stable for $\bar\ell$, then
\begin{equation*}
  \left|
  \mathbb E_{\mathcal S,\mathcal A}
  \left[
    \mathcal J(\mathcal A(\mathcal S))
    -
    \widehat{\mathcal J}_{\mathcal S}(\mathcal A(\mathcal S))
  \right]
  \right|
  \le \beta_n .
\end{equation*}
\end{theorem}

\begin{proof}
Let $\mathcal S^{(i)}$ be obtained from $\mathcal S$ by replacing $s_i$ with an
independent copy $s_i'$.  By exchangeability,
\begin{equation*}
\begin{split}
&\mathbb E_{\mathcal S,\mathcal A}
\left[
  \mathcal J(\mathcal A(\mathcal S))
  -
  \widehat{\mathcal J}_{\mathcal S}(\mathcal A(\mathcal S))
\right]  \\
&=
\frac{1}{n}\sum_{i=1}^{n}
\mathbb E
\left[
  \bar\ell(\mathcal A(\mathcal S);s_i')
  -
  \bar\ell(\mathcal A(\mathcal S);s_i)
\right].
\end{split}
\end{equation*}
Replacing $\mathcal S$ by $\mathcal S^{(i)}$ in the first term preserves the
joint distribution.  Hence each summand has absolute value at most $\beta_n$ by
Definition~\ref{def:app_uniform_stability}.  Averaging over $i$ proves the claim.
\end{proof}

\begin{theorem}[Stability of convex gated SGD]
\label{thm:app_convex_stability}
Assume that $\bar\ell(\theta;s)$ is convex in $\theta$, $L_g$-Lipschitz in
$\theta$, and $\beta_g$-smooth in $\theta$ for every $s$.  Consider SGD with
sample replacement,
\begin{equation*}
  \theta_{k+1}
  =
  \theta_k-\alpha_k\nabla_\theta\bar\ell(\theta_k;s_{i_k}),
  \qquad
  i_k\sim {\rm Unif}\{1,\ldots,n\},
\end{equation*}
with identical initialization on neighboring datasets and
$0\le\alpha_k\le2/\beta_g$.  Then the algorithm is
\begin{equation*}
  \beta_n
  \le
  \frac{2L_g^2}{n}\sum_{k=0}^{K-1}\alpha_k
\end{equation*}
uniformly stable for $\bar\ell$.
\end{theorem}

\begin{proof}
Let $\theta_k$ and $\theta_k'$ be coupled SGD iterates on neighboring datasets.
For a shared sample, the update map
$G_{\alpha,s}(\theta)=\theta-\alpha\nabla_\theta\bar\ell(\theta;s)$ is
nonexpansive.  This follows from convexity, $\beta_g$-smoothness, and
$\alpha\le2/\beta_g$.  When the differing sample is selected, both gradients
are bounded by $L_g$.  The distance can then increase by at most $2\alpha_kL_g$.
Therefore
\begin{equation*}
  \mathbb E\|\theta_{k+1}-\theta_{k+1}'\|
  \le
  \mathbb E\|\theta_k-\theta_k'\|
  +
  \frac{2\alpha_kL_g}{n}.
\end{equation*}
Since the initializations coincide, summing over $k$ gives
\begin{equation*}
  \mathbb E\|\theta_K-\theta_K'\|
  \le
  \frac{2L_g}{n}\sum_{k=0}^{K-1}\alpha_k .
\end{equation*}
The $L_g$-Lipschitz property of $\bar\ell$ gives the stated stability bound.
\end{proof}

\begin{theorem}[Conditional nonconvex stability]
\label{thm:app_nonconvex_stability}
Assume that $\bar\ell(\theta;s)$ is $L_g$-Lipschitz and $\beta_g$-smooth in
$\theta$ for every $s$, but not necessarily convex.  Under the same
sample-replacement coupling as in Theorem~\ref{thm:app_convex_stability}, SGD
satisfies
\begin{equation}
\label{eq:app_nonconvex_stability}
  \beta_n
  \le
  \frac{2L_g^2}{n}
  \sum_{k=0}^{K-1}
  \alpha_k
  \prod_{j=k+1}^{K-1}(1+\alpha_j\beta_g).
\end{equation}
\end{theorem}

\begin{proof}
For a shared sample, $\beta_g$-smoothness implies
\begin{equation*}
  \|G_{\alpha,s}(\theta)-G_{\alpha,s}(\theta')\|
  \le
  (1+\alpha\beta_g)\|\theta-\theta'\|.
\end{equation*}
For the differing sample, the additional increase is at most
$2\alpha_kL_g/n$ in expectation.  Thus
\begin{equation*}
  \mathbb E\|\theta_{k+1}-\theta_{k+1}'\|
  \le
  (1+\alpha_k\beta_g)\mathbb E\|\theta_k-\theta_k'\|
  +
  \frac{2\alpha_kL_g}{n}.
\end{equation*}
Unrolling the recursion from the common initialization and multiplying by
$L_g$ yields~\eqref{eq:app_nonconvex_stability}.
\end{proof}

The constants $L_g$ and $\beta_g$ are effective constants for the gated loss.
Because $0\le h\le1$, and the scale weights are convex, temporal averaging
does not increase these constants.  This holds when pointwise losses satisfy
the same bounds uniformly.  Overlap between temporal neighborhoods
may reduce empirical variance, but it does not justify a negative
generalization term without additional assumptions.  The nonconvex bound is
informative only when the product factor in
\eqref{eq:app_nonconvex_stability} is controlled.  Control may come from
stepsizes, the training horizon, or a local smoothness estimate.
Otherwise, it is a diagnostic inequality rather than a practical bound.

\section{Corrected convergence rate of NeuralMD}
\label{app:corrected_convergence}
The following theorem is a stationarity guarantee for the population objective
$\mathcal J$ in~\eqref{eq:app_pop_emp_loss}.  It applies to a randomly selected
iterate and does not claim convergence of the final iterate or of the PDE
solution error.

\begin{theorem}[Random-iterate stationarity]
\label{thm:app_random_iterate_stationarity}
Assume that $\mathcal J$ is lower bounded by $\mathcal J_{\inf}$ and has an
$H_g$-Lipschitz gradient.  Let $G_k$ be an unbiased stochastic gradient,
\begin{equation*}
  \mathbb E[G_k\mid\theta_k]=\nabla_\theta\mathcal J(\theta_k),
  \qquad
  \mathbb E[\|G_k\|^2\mid\theta_k]\le M_g^2 .
\end{equation*}
For SGD, $\theta_{k+1}=\theta_k-\alpha_kG_k$, define
\begin{equation*}
  \Pr(I_K=k)=\frac{\alpha_k}{\sum_{j=0}^{K-1}\alpha_j}.
\end{equation*}
Then
\begin{equation}
\label{eq:app_stationarity_general}
  \mathbb E\|\nabla_\theta\mathcal J(\theta_{I_K})\|^2
  \le
  \frac{
    \mathcal J(\theta_0)-\mathcal J_{\inf}
    +\frac{H_gM_g^2}{2}\sum_{k=0}^{K-1}\alpha_k^2
  }{
    \sum_{k=0}^{K-1}\alpha_k
  } .
\end{equation}
With $\alpha_k=(k+1)^{-1/2}$, this gives $O(\log K/\sqrt K)$.  With the
finite-horizon choice $\alpha_k=\alpha/\sqrt K$, where $\alpha>0$ is fixed, it
gives $O(K^{-1/2})$.
\end{theorem}

\begin{proof}
By $H_g$-smoothness,
\begin{equation*}
\mathcal J(\theta_{k+1})
\le
\mathcal J(\theta_k)
-\alpha_k\langle\nabla\mathcal J(\theta_k),G_k\rangle
+\frac{H_g\alpha_k^2}{2}\|G_k\|^2 .
\end{equation*}
Taking conditional expectation and using unbiasedness and the bounded second
moment gives
\begin{equation*}
\mathbb E[\mathcal J(\theta_{k+1})]
\le
\mathbb E[\mathcal J(\theta_k)]
-\alpha_k\mathbb E\|\nabla\mathcal J(\theta_k)\|^2
+\frac{H_gM_g^2}{2}\alpha_k^2 .
\end{equation*}
Rearranging, summing over $k=0,\ldots,K-1$, and using the lower bound
$\mathcal J_{\inf}$ yields
\begin{equation*}
\sum_{k=0}^{K-1}\alpha_k
\mathbb E\|\nabla\mathcal J(\theta_k)\|^2
\le
\mathcal J(\theta_0)-\mathcal J_{\inf}
+\frac{H_gM_g^2}{2}\sum_{k=0}^{K-1}\alpha_k^2 .
\end{equation*}
Dividing by $\sum_k\alpha_k$ gives~\eqref{eq:app_stationarity_general}.
\end{proof}

\section{Corrected gradient estimation error of NeuralMD}
\label{app:corrected_gradient_estimation}
The stochastic gradient used by the gated multiscale sampling is a Monte Carlo
estimate of the gradient of~\eqref{eq:app_gated_loss}.  The next theorem
records the exact unbiasedness and variance relation.

\begin{theorem}[Monte Carlo gradient estimation]
\label{thm:app_mc_gradient}
Assume differentiation can be interchanged with the conditional expectation over
$\vartheta$.  Also assume that
\begin{equation*}
  g(\theta;s,\vartheta)
  =
  h(t_\vartheta)\nabla_\theta
  \ell_{\rm pt}(\theta;\bm{x},t_\vartheta)
\end{equation*}
has a finite second moment.  Then
\begin{equation*}
  \nabla_\theta\bar\ell(\theta;s)
  =
  \mathbb E_{\vartheta\mid s}g(\theta;s,\vartheta).
\end{equation*}
For $M$ independent perturbations $\vartheta_1,\ldots,\vartheta_M$ and
\begin{equation*}
  \widehat g_M(\theta;s)
  =
  \frac{1}{M}\sum_{j=1}^{M}g(\theta;s,\vartheta_j),
\end{equation*}
we have
\begin{equation}
\label{eq:app_mc_variance}
  \mathbb E_{\vartheta\mid s}\widehat g_M(\theta;s)
  =
  \nabla_\theta\bar\ell(\theta;s),
  \qquad
  \mathbb E_{\vartheta\mid s}
  \|\widehat g_M(\theta;s)-\nabla_\theta\bar\ell(\theta;s)\|^2
  =
  \frac{1}{M}\operatorname{tr}
  \operatorname{Cov}_{\vartheta\mid s}(g(\theta;s,\vartheta)).
\end{equation}
\end{theorem}

The interchange condition holds under a standard dominated-differentiation condition.
For example, suppose $\ell_{\rm pt}$ is differentiable near the iterates.
Assume there is an integrable envelope $m_s(\vartheta)$ such that
$\|\nabla_\theta \ell_{\rm pt}(\theta;\bm{x},t_\vartheta)\|
\le m_s(\vartheta)$ in that neighborhood.  For unbounded residual losses, use
uniform-integrability or truncated-loss conditions instead.  These conditions can depend on $\varepsilon$ because
the NKGE residual contains stiff coefficients.

\begin{proof}
The first identity follows directly from
Definition~\eqref{eq:app_gated_loss} and the assumed interchange of
differentiation and expectation.  For the second identity, write
$\Sigma_s=\operatorname{Cov}_{\vartheta\mid s}(g(\theta;s,\vartheta))$.  Independence of
$\vartheta_1,\ldots,\vartheta_M$ gives
\begin{equation*}
  \operatorname{Cov}_{\vartheta\mid s}(\widehat g_M)
  =
  \frac{1}{M}\Sigma_s .
\end{equation*}
The mean-squared Euclidean error equals the trace of this covariance, which
proves~\eqref{eq:app_mc_variance}.
\end{proof}

\begin{theorem}[Stationarity with gradient estimation error]
\label{thm:app_stationarity_grad_error}
Assume that $\mathcal J$ is lower bounded by $\mathcal J_{\inf}$ and has an
$H_g$-Lipschitz gradient.  Suppose the update uses
\begin{equation*}
  \widehat G_k
  =
  \nabla\mathcal J(\theta_k)+b_k+\eta_k,
  \qquad
  \mathbb E[\eta_k\mid\theta_k]=0,
\end{equation*}
where $\|b_k\|\le B_k$ and
$\mathbb E[\|\eta_k\|^2\mid\theta_k]\le \sigma_k^2/M$.  If
$0<\alpha_k\le1/(2H_g)$ and
$\Pr(I_K=k)=\alpha_k/\sum_{j=0}^{K-1}\alpha_j$, then
\begin{equation}
\label{eq:app_grad_error_stationarity}
\begin{split}
  \mathbb E\|\nabla\mathcal J(\theta_{I_K})\|^2
  \le&
  \frac{4(\mathcal J(\theta_0)-\mathcal J_{\inf})}
       {\sum_{k=0}^{K-1}\alpha_k}
  +
  \frac{6\sum_{k=0}^{K-1}\alpha_k B_k^2}
       {\sum_{k=0}^{K-1}\alpha_k}  \\
  &+
  \frac{2H_g\sum_{k=0}^{K-1}\alpha_k^2\sigma_k^2/M}
       {\sum_{k=0}^{K-1}\alpha_k}.
\end{split}
\end{equation}
\end{theorem}

\begin{proof}
Let $g_k=\nabla\mathcal J(\theta_k)$.  Smoothness gives
\begin{equation*}
\mathcal J(\theta_{k+1})
\le
\mathcal J(\theta_k)
-\alpha_k\langle g_k,\widehat G_k\rangle
+\frac{H_g\alpha_k^2}{2}\|\widehat G_k\|^2 .
\end{equation*}
Taking conditional expectation gives the next bound.
We use $\mathbb E[\eta_k\mid\theta_k]=0$.
We also use Young's inequality
$|\langle g_k,b_k\rangle|\le \|g_k\|^2/4+B_k^2$.
Finally, use $\|g_k+b_k\|^2\le2\|g_k\|^2+2B_k^2$ and
$\alpha_k\le1/(2H_g)$.
\begin{equation*}
\mathbb E[\mathcal J(\theta_{k+1})\mid\theta_k]
\le
\mathcal J(\theta_k)
-\frac{\alpha_k}{4}\|g_k\|^2
+\frac{3\alpha_k}{2}B_k^2
+\frac{H_g\alpha_k^2}{2}\frac{\sigma_k^2}{M}.
\end{equation*}
After taking full expectation, summing over $k$, and using the lower bound
$\mathcal J_{\inf}$, division by $\sum_k\alpha_k$ gives
\eqref{eq:app_grad_error_stationarity}.
\end{proof}

The term involving $B_k$ shows the effect of persistent gradient bias.
It leads only to convergence to a stationarity neighborhood.
In the unbiased Monte Carlo case of Theorem~\ref{thm:app_mc_gradient}, $B_k=0$.
Increasing the number of temporal perturbations then reduces variance by $1/M$.
The constants may depend on the stiffness parameter $\varepsilon$
through the residual derivatives; the result does not remove that dependence.

\section{Additional results}
This section reports additional results that supplement the main text.
They include hyperparameter analysis, inverse experiments, and optimizer comparisons.

\section{Inverse problem experiments}
As we stated in the implementations, NeuralMD can also be applied to tasks with data loss for the inverse problem.
We further test NeuralMD on inverse problems with sparse measurements.
In these tasks, part of the solution is unknown and must be recovered.
Specifically, we recover the NKGE initial condition from limited spatiotemporal observations.
This setting is challenging because the model must predict forward evolution.
It must also reconstruct the unknown initial data.

We test 1D NKGE with varying $\varepsilon$ values.
Only $5\%$ of the spatiotemporal domain is observed.
The goal is to infer the initial condition $\phi_1(x)$ and $\phi_2(x)$ from the available measurements.
We compare NeuralMD with vanilla PINNs and other baseline methods using the same experimental setup.
The preliminary results suggest that NeuralMD can reconstruct the initial conditions with lower relative errors than the tested baselines in these settings.
This behavior is consistent with the multiscale decomposition strategy, which separates the oscillatory dynamics before solving the inverse problem.

Furthermore, we investigate the robustness of NeuralMD to noise in the observations.
We add Gaussian noise with varying standard deviations ($\sigma = 0.01, 0.05, 0.1$) to the measurements and evaluate the recovery performance.
NeuralMD maintains stable reconstruction errors in the tested noisy-observation settings.
This behavior may relate to gated temporal mixing.
The mechanism averages information across scales and can reduce observation noise.

\section{Impact of optimizers}
In this work, we adhere to established benchmarks and use Adam+L-BFGS for standard evaluations.
This section examines the influence of different optimizers on the performance of NeuralMD.

We evaluate how optimization strategies affect convergence and final accuracy.
We compare Adam alone, L-BFGS alone, and the default Adam+L-BFGS hybrid.

Adam provides rapid initial convergence through adaptive learning rates.
It scales gradients using first and second moments.
This is useful early in training, when careful step-size adaptation is needed.
However, we observe that Adam alone may converge to local minima prematurely, resulting in suboptimal final accuracy for highly oscillatory PDEs.

The L-BFGS optimizer, as a quasi-Newton method, approximates the Hessian matrix to achieve second-order optimization convergence.
This approach often yields higher accuracy in the final stages of training by performing precise gradient-based updates.
Nevertheless, L-BFGS requires substantial memory to store curvature information and may converge slowly if the initial point is far from the optimum.

The hybrid Adam+L-BFGS strategy leverages the strengths of both optimizers.
We first apply Adam for 500 iterations to reach a favorable parameter region.
L-BFGS then runs for 500 iterations to refine the solution.
This two-stage approach gives the best performance among the optimizer choices tested for NeuralMD in our experiments ($\varepsilon = 0.8, 0.5, 0.1, 0.01$).

We further investigate the impact of learning rate scheduling on NeuralMD's performance.
Cosine annealing and step decay were also tested.
The default constant learning rate remains competitive in our settings.

\section{Hyperparameter Analysis}
This section evaluates NeuralMD under several hyperparameter configurations.
We vary perturbation counts $(k_1,k_2,k_3)$, region sizes $(R_1,R_2,R_3)$, and \#scales.

\subsection{Number of Random Time Perturbations}
The number of random time perturbations $\{k_1, k_2, k_3\}$ controls the diversity of temporal sampling within each scale.
For scale $l\in\{1,2,3\}$, we sample $k_l$ perturbed time points.
They are $\{t+\delta_t\}_{j=1}^{k_l}$ with $\delta_t\sim\mathcal U(-R_l,R_l)$.
Increasing $k_l$ enhances the Monte Carlo approximation quality of the gated temporal loss, reducing the variance of gradient estimates.
Formally, the gated multiscale time-region loss can be written as
\begin{equation*}
\mathcal{L}^{\text{gms}}(u_\theta, \mathbf{x}) = \sum_{l=1}^{\#\text{scale}} \lambda_l \cdot \frac{1}{k_l} \sum_{j=1}^{k_l} h(t + \delta_{t,j}) \mathcal{L}_{\text{pt}}(u_\theta, \mathbf{x}, t + \delta_{t,j}),
\end{equation*}
where $\delta_{t,j} \sim \mathcal{U}(-R_l, R_l)$ and $h(\cdot)$ denotes the gating function.
As $k_l \to \infty$, the empirical average converges to the true expectation $\mathbb{E}_{\delta_t}[\cdot]$.
Our diagnostics indicate that increasing $k_l$ from 3 to 7 improves performance.
Returns diminish beyond $k_l=7$ because computational cost increases.

\subsection{Perturbation Region Size}
The perturbation region size $\{R_1, R_2, R_3\}$ determines the temporal window for random perturbations at each scale.
These values are related to the local temporal correlation scale.
The NKGE carrier frequency scales as $O(\varepsilon^{-2})$ in the nonrelativistic limit.
Thus stage-II perturbations should be capped at an $O(\varepsilon^2)$ scale.
For $\varepsilon\approx1$, larger neighborhoods can be used.
For $\varepsilon\ll1$, the radius must avoid averaging unrelated carrier phases.
We set $R_1 = 1 \times 10^{-2}$, $R_2 = 5 \times 10^{-2}$, and $R_3 = 9 \times 10^{-2}$ in our experiments.
The choice of perturbation region sizes follows the principle that $R_l$ should be proportional to the local oscillation period, i.e., $R_l \sim \varepsilon_l / \omega$ where $\omega$ denotes the characteristic frequency.

\subsection{Number of Scales}
The number of temporal scales $\#\text{scale}$ controls the granularity of multiscale decomposition.
We denote the scales as $\{\Omega_{t_1}, \Omega_{t_2}, \ldots, \Omega_{t_{\#\text{scale}}}\}$ where $\Omega_{t_l} = (t - R_l, t + R_l)$.
Incorporating additional scales enables the model to capture oscillatory structures across a broader range of frequencies.
The mixing weights $\{\lambda_l\}_{l=1}^{\#\text{scale}}$ are learned through the time-region mixing layer, which can be expressed as
\begin{equation*}
\tilde{u}_\theta = \text{MLP}\left(\sum_{l=1}^{\#\text{scale}} \lambda_l \cdot u_\theta(\Omega_{t_l})\right),
\end{equation*}
where the mixing weights satisfy $\sum_{l=1}^{\#\text{scale}} \lambda_l = 1$ and $\lambda_l \geq 0$.
Our diagnostic ablation indicates that increasing $\#\text{scale}$ from 1 to 3 improves performance in the tested NKGE setting.
Further increasing the number of scales may lead to overfitting and increased computational overhead.

\begin{itemize}
  \item \emph{Increasing the number of perturbations improves the model's performance by reducing gradient estimation variance.}
  \item \emph{The perturbation region size is determined by the NKGE property and should scale with the oscillation period.}
  \item \emph{Incorporating additional scales can enhance the model's performance by capturing multi-frequency oscillatory structures.}
\end{itemize}

\vskip 0.2in

\bibliography{references}

\end{document}